\documentclass[11pt]{article}

\usepackage{amssymb}
\usepackage{amsmath}
\usepackage{theorem}
\usepackage{epsfig}
\usepackage{verbatim}
\usepackage{graphicx}
\usepackage[all]{xy}

\usepackage{times}
\usepackage[T1]{fontenc}

\usepackage{color}

\numberwithin{equation}{section}

\newtheorem{thm}{Theorem}[section]

\newtheorem{rem}[thm]{Remark}

\newtheorem{lemma}[thm]{Lemma}

\newcommand{\R}{\mathbb{R}}
\newcommand{\C}{\mathbb{C}}
\newcommand{\Z}{\mathbb{Z}}

\newcommand{\T}{\mathbb{T}}

\newcommand{\supp}{{\rm supp}\thinspace}

\newcommand{\dx}{\partial_x}

\newcommand{\la}{\Lambda}

\newcommand{\al}{\alpha}
\newcommand{\ep}{\varepsilon}

\newcommand{\un}{\underline{z}}

\newcommand{\pa}{\partial}

\newcommand{\hz}{\hat{\zeta}}
\newcommand{\hh}{\hat{h}}
\newcommand{\hw}{\hat{w}}
\newcommand{\hze}{\hat{z}}
\newcommand{\that}{\hat{t}}

\begin{document}

\title{Breakdown of smoothness for the Muskat problem}
\date{27-1-2012}
\author{\'Angel Castro, Diego C\'ordoba, Charles Fefferman,\\
 and Francisco Gancedo.}

\maketitle
 \abstract{In this paper we show that there exist  analytic initial data in the stable regime for the Muskat problem such that the solution turns to the unstable regime and later breaks down, i.e., no longer belongs to $C^4$.}
 \tableofcontents
\section{Introduction}

The Muskat equation governs the motion of an interface separating two fluids in a porous medium (e.g. oil and water in sand). See \cite{Muskat}. The same equation governs an interface separating two fluids trapped between two closely spaced parallel vertical plates (a "Helle Shaw cell"). See \cite{H-S}.

In this paper, we exhibit a solution of the Muskat equation for which the interface is initially smooth and stable, but later enters an unstable regime, and still later develops a singularity.

Let us briefly recall the Muskat equation and its derivation from physical laws. Imagine the plane $\R^2$ filled with two incompressible fluids that cannot mix. At time $t$, let $\Omega^1(t)$ be the region occupied by FLUID 1, and let $\Omega^2(t)$ be the region occupied by FLUID 2. These two regions are separated by their common boundary $\pa\Omega^1(t)=\pa\Omega^2(t)$ (the "interface"), which we represent as a parametrized curve
$$\{ z(\alpha,t) \,:\, \alpha\in \R\},$$
where $$z(\alpha,t)=(z_1(\alpha,t),z_2(\alpha,t))\in\R^2.$$

We suppose that all points $(x_1,x_2)\in \R^2$ with large positive $x_2$ belong to $\Omega^1(t)$ and that all $(x_1,x_2)$ with large negative $x_2$ belong to $\Omega^2(t)$.

Let $p^i(x,t)$ and $u^i(x,t)$ denote (respectively) the pressure and fluid velocity in region $\Omega^i(t)$, for $i=1,2$. Then the $p^i$ and $u^i$ satisfy the following conditions.
\begin{enumerate}
\item $\text{div} u^i=0$ in $\Omega^i(t)$, ($i=1,2$), i.e., the fluids are incompressible.
\item $\frac{\mu_i}{\kappa}u^i=-\nabla p^i-g(0,\rho_i)$ in $\Omega^i(t)$ ($i=1,2$); this is Darcy's law \cite{Muskat} for fluid flow in a porous medium. (Here $\kappa$ is a physical constant associated to the porous medium, $\mu_i$ and $\rho_i$ are physical constants associated to FLUID i, and $g$ is the acceleration due to gravity. In particular, $\rho_i$ is the density of FLUID i and $\mu_i$ is the viscosity of FLUID i.)

 Moreover, on the interface,

 \item $p^1=p^2$

 and

 \item $u^1-u^2$ is tangent to the interface.

 Finally, the interface moves with the fluid, i.e.,

 \item $\pa_t z(\alpha,t)=u^1(z(\alpha,t),t)+c_1(\alpha,t)\pa_\alpha z(\alpha,t)$.

 Here the function $c_1(\alpha,t)$ may be chosen arbitrarily; it affects only the parametrization of the interface, which has no intrinsic physical meaning. (We could have used $u^2$ in 5, in place of $u^1$.)

See, e.g., \cite{DY} where the Muskat problem is interpreted in terms of weak solutions.
\end{enumerate}

We may study two scenarios, namely
\begin{description}
\item ASYMPTOTICALLY FLAT INTERFACES: $\lim_{\alpha\to \pm \infty} [z(\alpha,t)-(\alpha,0)]=0$

and

\item PERIODIC INTERFACES: $z(\alpha,t)-(\alpha,0)$ is $2\pi-$periodic in $\alpha$ for fixed $t$.
\end{description}

We restrict attention to the case of two fluids with the same viscosity but different densities. In that case, $\mu^1=\mu^2$ but $\rho_1\neq\rho_2$ in 2.

Equations 1,...,5 take pressure and gravity into account but neglect surface tension.

Once we know the interface $z(\alpha,t)$, we can derive explicit formulas for $u^1$, $u^2$, $p^1$ and $p^2$ from 1,...,4 and elementary potential theory. Substituting these formulas into 5, and making a convenient choice of the function $c_1(\alpha,t)$, we obtain the Muskat equation for $z(\alpha,t)$.

For asymptotically flat interfaces, the Muskat equation takes the form
\begin{equation}
\pa_t z_\mu(\alpha,t)=\frac{\rho_2-\rho_1}{2\pi}P.V.\int_{-\infty}^\infty\frac{z_1(\alpha,t)-z_1(\beta,t)}{|z(\alpha,t)-z(\beta,t)|^2}
[\pa_\alpha z_\mu(\alpha,t)-\pa_\beta z_\mu(\beta,t)]d\beta\quad (\mu=1,2),\label{6}
\end{equation}
where "$P.V.$" denotes a principal value integral (the integral does not converge absolutely at infinity).

For periodic interfaces, the Muskat equation is
\begin{equation}
\pa_t z_\mu(\alpha,t)=\frac{\rho_2-\rho_1}{2\pi}\int_{\R/2\pi\Z}\frac{\sin(z_1(\alpha,t)-z_1(\beta,t))[\pa_\alpha z_\mu(\alpha,t)-\pa_\beta z_\mu(\beta,t)]}
{\cosh(z_2(\alpha,t)-z_2(\beta,t))-\cos(z_1(\alpha,t)-z_1(\beta,t))}
d\beta\quad (\mu=1,2).\label{7}
\end{equation}

Given an initial interface $z^0(\alpha)=(z_1^0(\alpha),z_2^0(\alpha))$ ($\alpha\in\R$), we try to solve \eqref{6} or \eqref{7}, with the initial condition
\begin{equation}\label{8}
z(\alpha,0)=z^0(\alpha)\quad (\text{all $\alpha\in\R$}).
\end{equation}
There is no significant difference between the two flavors \eqref{6}, \eqref{7} of the Muskat equation. Except for the introduction, we will work in the periodic setting. For the moment, we work with asymptotically flat interfaces, since the relevant formulas are then a bit simpler.

\begin{figure}[h!]
\centering
\includegraphics[width=1\textwidth]{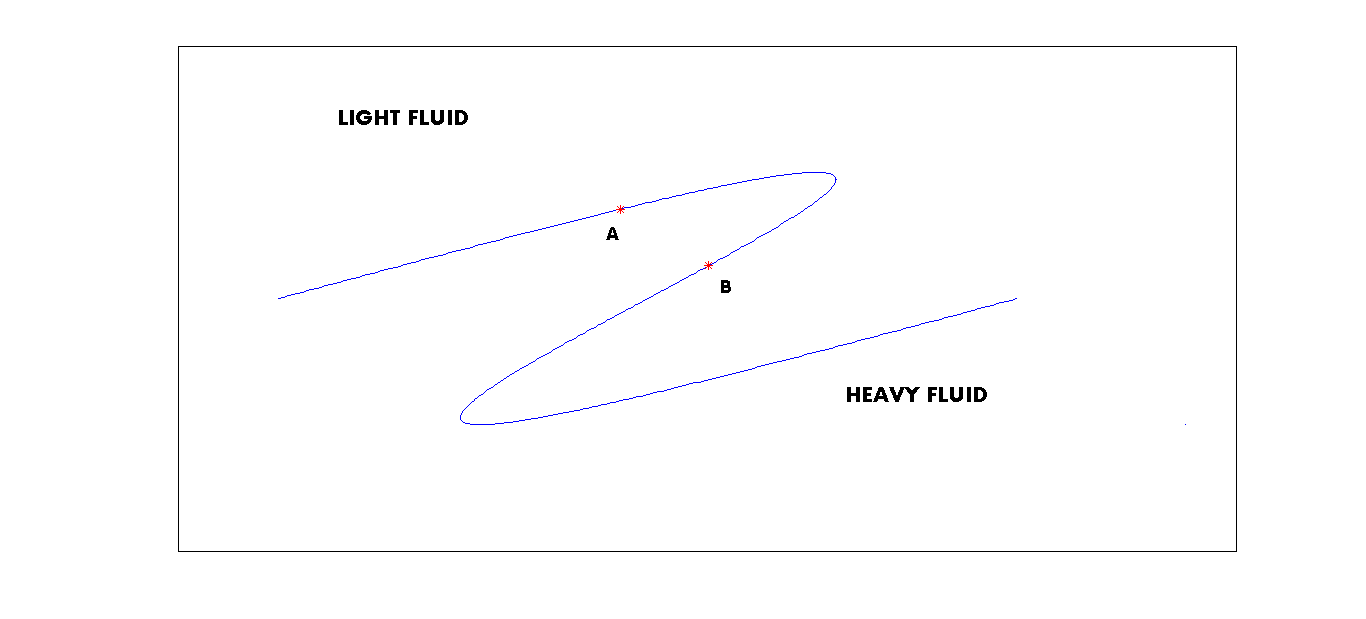}\caption{Non-graph interface.}
\label{X}
\end{figure}

Suppose first that the interface is the graph of a function, $x_2=f(x_1)$. Then it is well-known that the initial-value problem \eqref{6}, \eqref{8} is linearly stable for positive time, if and only if the heavy fluid lies below the interface and the light fluid lies above it (i.e., $\rho_2>\rho_1$). If instead the heavy fluid lies above the light fluid ($\rho_1>\rho_2$), then the initial-value problem \eqref{6}, \eqref{8} is linearly stable for negative time. See \cite{DY}. If the interface is not the graph of a function, as in Figure \ref{X}, then the heavy fluid lies below the light fluid locally near point A, and the reverse holds locally near point B.

In that case, the initial-value problem \eqref{6}, \eqref{8} is linearly unstable, both for positive time (thanks to point B), and for negative time (thanks to point A).

\begin{figure}[h!]
\centering
\includegraphics[width=1\textwidth]{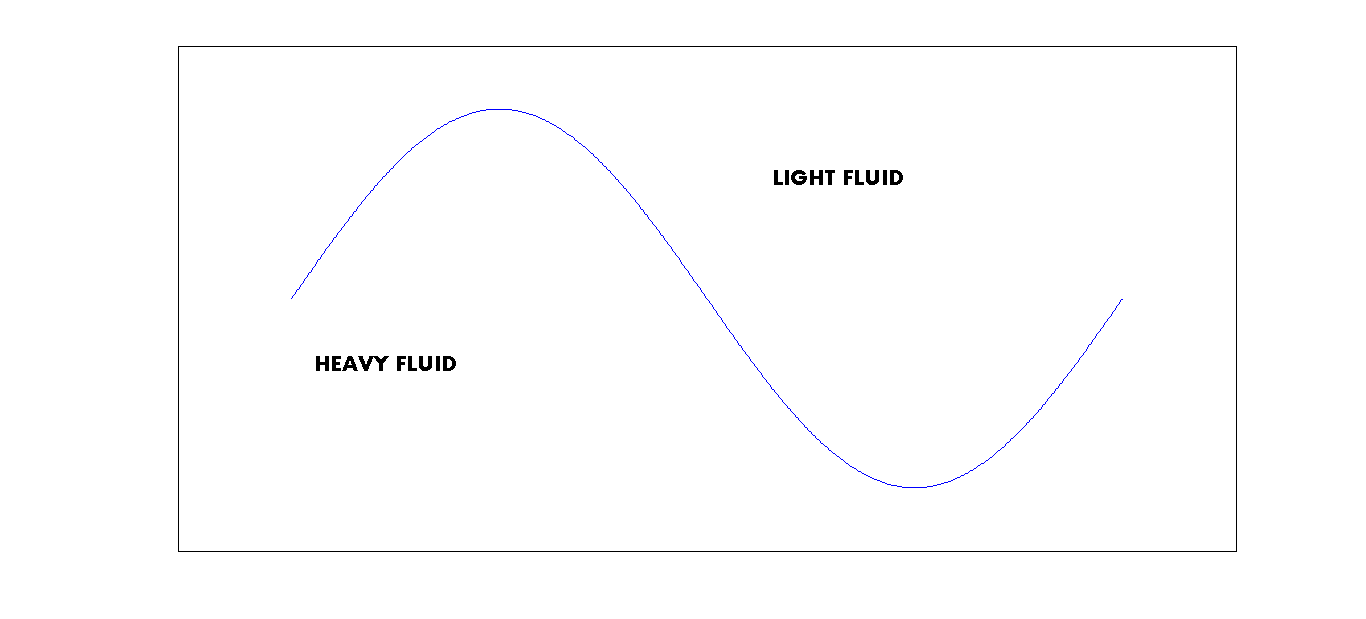}\caption{Initially stable interface.}
\label{Ya}
\end{figure}

\begin{figure}[h!]
\centering
\includegraphics[width=1\textwidth]{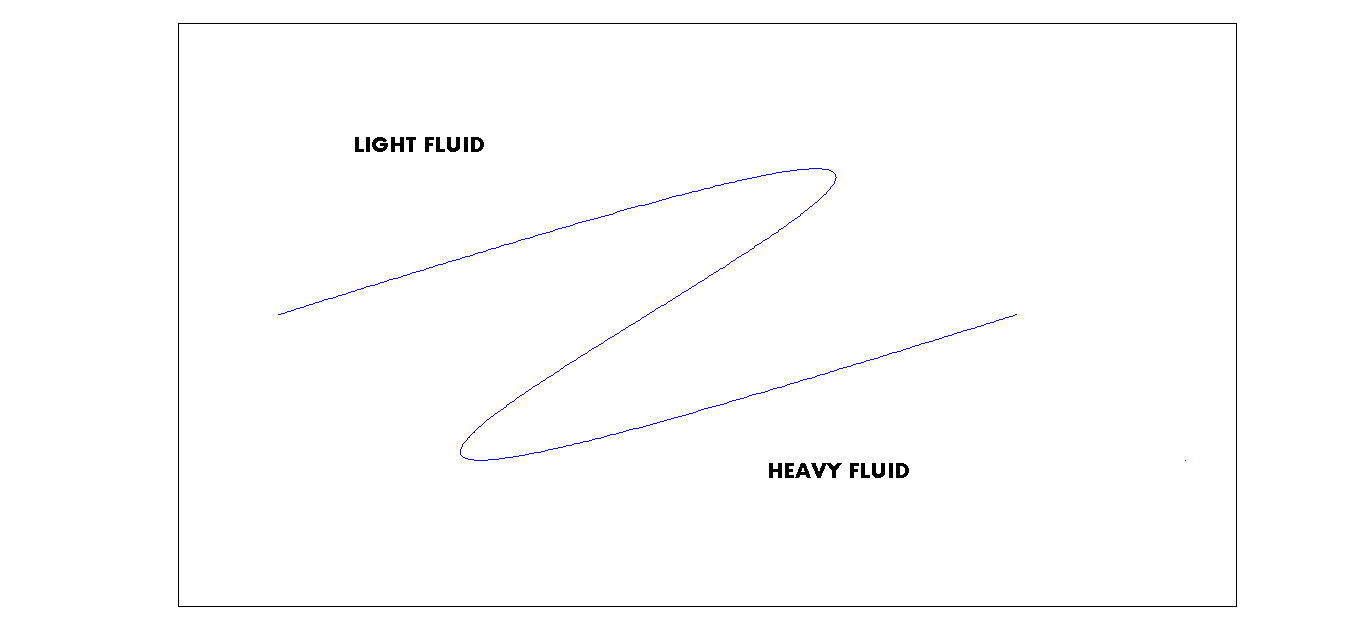}\caption{Turnover.}
\label{Yb}
\end{figure}

Therefore, one asks whether an interface may begin in a (presumably) stable configuration as in Figure \ref{Ya}, and
then "turn over" to a (presumably) highly unstable configuration as in Figure \ref{Yb}.

\begin{figure}[h!]
\centering
\includegraphics[width=1\textwidth,height=8cm]{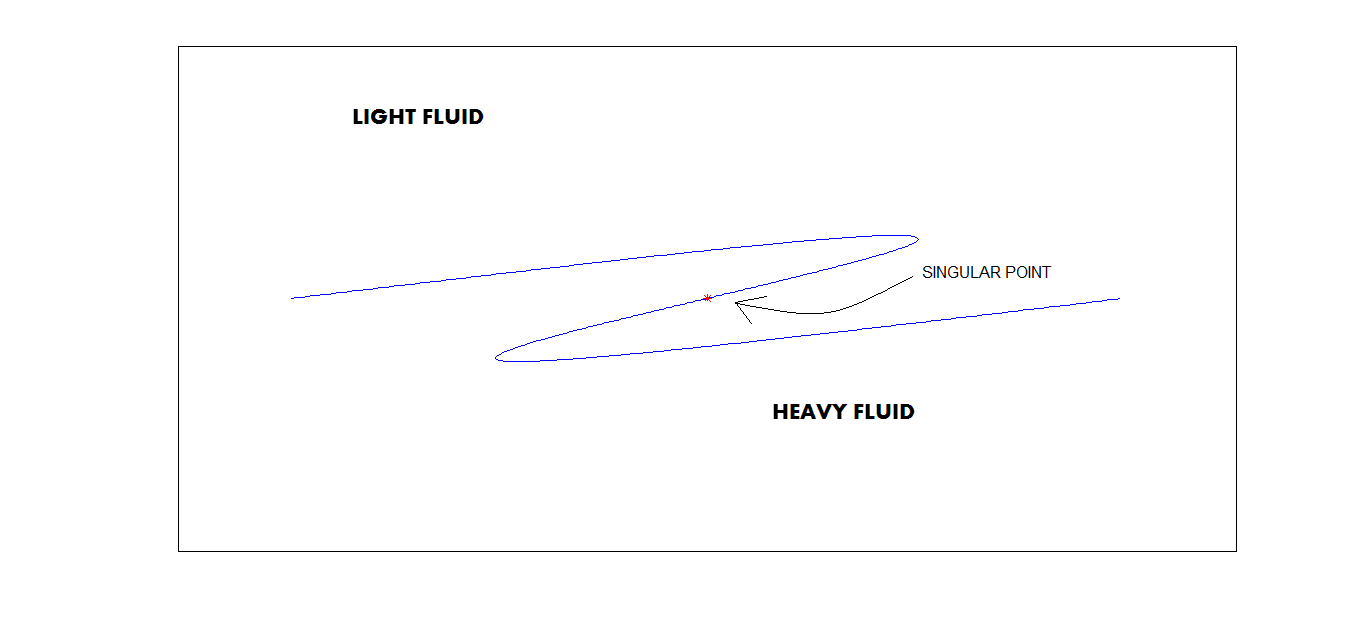}\caption{Breakdown.}
\label{Yc}
\end{figure}

For "very small" initial data, Constantin-Pugh \cite{Peter}, Yi \cite{Yi2},  Siegel-Caflish-Howison \cite{SCH}, C\'ordoba-Gancedo \cite{DY} and Escher-Matioc \cite{Esch2}, prove  that such a "turnover" cannot occur. It is shown in \cite{ccgs} that a turnover cannot occur if the initial interface has the form
$$\{(x_1,x_2)\,:\, x_2=f(x_1)\}$$
with $|f'|<1$ everywhere.

Nevertheless, Castro et al.  \cite{turning}, \cite{CCFGL}  show that a turnover from Figure \ref{Ya} to Figure \ref{Yb} can occur. Moreover, the Muskat solution exists and remains smooth for some time after the turnover, despite the instability of interfaces that look like Figure \ref{Yb}.

The reason why the Muskat solution continues to exist in a highly unstable regime is that, by the time the turnover occurs, the interface has become real-analytic. In fact, starting from a smooth, stable configuration as in Figure \ref{Ya}, the Muskat solution becomes instantly real-analytic. For $t>0$, the function $z(\alpha,t)$ defining the interface $\partial\Omega^{1}(t)=\partial\Omega^{2}(t)=\{z(\alpha,t): \alpha\in \R\}$ continues analytically to a strip $\{\alpha\in\C: |\Im \alpha|<h(t)\}$, with $h(t)>ct$ for small $t$. As $t$ approaches the turnover time, the function $h(t)$ may decrease, but it remains strictly positive up to the moment when the turnover occurs (i.e. up to the appearance of a vertical tangent to the interface).

Therefore, the variants \cite{Nirenberg}, \cite{Nishida} of the Cauchy-Kowalewski theorem  show that a real-analytic Muskat solution continues to exist for a short time after the turnover. The interface becomes (presumably) more and more unstable as the turnover progresses. For details see \cite{CCFGL}.


In \cite{DY} (following the scheme of \cite{SCH} for a scenario of two fluids with different viscosities and the same density) the construction of global-in-time solutions leads to ill-posedness (in Sobolev spaces in the Hadamard sense) for the unstable regime by going backwards in time. In \cite{Otto1} and \cite{Otto2} it was proposed, in the unstable regime,  a relaxation of the phase distribution  that leads to a well-posed problem and captures a mixing profile. From a completely different approach, a convex integration framework, unstable weak solutions are constructed for porous media in \cite{Laszlo, Shvydkoy, CFG}.

It is natural to believe that at some finite time $T$, the Muskat solution will break down. However, the technique in \cite{turning}, \cite{CCFGL} is not strong enough to prove that such a breakdown occurs.

In this paper, we show that a breakdown occurs, by proving the following result.

\begin{thm}\label{conclusion} (Main Theorem). There exists a solution $z(x,t)$ of the Muskat equation (\ref{7}), defined for $t\in[t_0, t_2]$, $\alpha\in\R /2\pi\Z$, such that the  following hold:
\begin{description}
\item A. At time $t_0$, the interface is a graph.
\item B. At some time $t_1\in (t_0, t_2),$ the interface is no longer a graph.
\item C. For each time $t\in[t_0,t_2)$, the interface is real analytic.
\item D. At time $t_2$, the interface is $C^3$ smooth but not $C^4$ smooth.
\end{description}
Thus, at time $t_2$, there is no longer a $C^4$ solution of the Muskat equation.
\end{thm}

Note that (D) includes the assertion that the curve $\alpha \rightarrow z(\alpha, t_2)$ cannot be made $C^4$-smooth by reparametrization. It is not enough merely to show that $z(\alpha, t_2)\notin C^4_{loc}$. In our proof of the Main Theorem, the failure of $C^4$ smoothness in (D) will occur at a single point, as shown in Figure \ref{Yc}.

To start to understand why Muskat solutions behave as described above, let us linearize equation \eqref{6} about the trivial solution $z(\alpha,t)=(\alpha,0),$ a stationary horizontal interface. Perturbing $z(\alpha,t)$ to $z(\alpha,t) + \delta z(\alpha,t)$, and keeping only terms up to first order in $\delta z$ in \eqref{6}, we obtain the linear equation
\begin{equation}\label{9}
\pa_t [\delta z_{\mu}(\alpha)]=\frac{\rho_2-\rho_1}{2\pi}P.V.\int_{-\infty}^{\infty} \frac{\pa_{\alpha}[\delta z_{\mu}(\alpha,t)] - \pa_\beta [\delta z_\mu(\beta,t)]}{\alpha - \beta}d\beta
\end{equation}
with $\mu = 1,2$.

The operator
\begin{equation}
f \rightarrow \Lambda f(x):= - \frac{1}{\pi} P.V. \int_{-\infty}^{\infty} \frac{f'(x) -f'(y)}{x-y}dy
\end{equation}
is the square root of the Laplacian in one dimension, i.e.,
$\widehat{\Lambda f}(\xi) =|\xi|\widehat{f}(\xi)$ (all $\xi\in \R$), where $\,\,\widehat{}$ denotes the Fourier transform.

The equation \eqref{9} takes the form
\begin{align}
\pa_t F(x,t) =&\frac{(\rho_1 -\rho_2)}{2} \Lambda F(x,t)\label{10}\\
F(x,0) =&f(x) \quad\quad \text{(given)}\nonumber.
\end{align}

The operator $\Lambda$ has a well-known connection with harmonic functions (see \cite{Eli}). In fact, let $u(x+iy)$ denote the Poisson integral of $f(x)$, defined for $y>0$. Similarly, for each fixed $t$, let $U(x+iy,t)$ denote the Poisson integral of $F(x,t)$. Then \eqref{10} is equivalent to the initial-value problem
\begin{align*}
\pa_t U(x+iy,t)=&-\frac{(\rho_1 -\rho_2)}{2} \pa_y U(x+iy,t)\\
U(x+iy,0)=&u(x+iy),
\end{align*}
which (formally) admits the solution $U(x+iy,t)=u(x+i[y- \frac{(\rho_1 - \rho_2)}{2}t])$.

In particular setting $y=0$, we find that
\begin{equation}\label{11}
F(x,t)=u(x-i\frac{(\rho_1-\rho_2)}{2}t).
\end{equation}
This is the formal solution of our linearized Muskat equation \eqref{10}.

 It is now obvious that the sign of $\rho_1-\rho_2$ plays a crucial role, and that solutions are real-analytic in the stable case. Indeed the harmonic function $u(x+iy)$ is real-analytic in the upper half-plane $y>0$ but undefined in the lower half-plane $y<0$ and unless $f$ is real-analytic. Therefore, if $(\rho_2-\rho_1)t>0$, then \eqref{11} yields a real-analytic function of $x$; more precisely, $F(x,t)$ continues analytically to the strip $|\Im x|<\frac{(\rho_2-\rho_1)}{2}t$. On the other hand, if $(\rho_2-\rho_1)t<0$, then \eqref{11} makes no sense unless the initial datum $f(x)$ continues analytically to the strip $|\Im x|<|\frac{(\rho_2-\rho_1)}{2}t|$.

 Note that the same ``Rayleigh-Taylor" condition that predicts stability as time flows forward also predicts instability as time flows backward.

  More generally, if we linearize the Muskat equation about any given solution $z(\al,t)=(z_1(\al,t),z_2(\al,t))$, we arrive at a linearized Muskat equation of the form
\begin{align*}
F_t(x,t)&=-\sigma(x,t)\Lambda F(x,t)+\text{HARMLESS TERMS},\\
F(x,0)&=f(x),
\end{align*}
where
$$
\sigma(x,t)=\frac{(\rho_2-\rho_1)}{2}\frac{ \dx z_1(x,t)}{(\dx z_1(x,t))^2+(\dx z_2(x,t))^2}\qquad \text{the ``Rayleigh-Taylor function."}
$$

We therefore expect that our Muskat solution will be real-analytic for $t>0$ provided $\sigma(x,t)>0$ everywhere, and highly unstable or non-existent for $t>0$ if $\sigma(x,t)<0$ everywhere.

 If we pass from initial conditions at time zero to solutions at times $t<0$, then the roles of $\sigma>0$ and $\sigma<0$ are reversed.

 Observe that (assuming $\rho_2>\rho_1$) the function $\sigma(x,t)$ is positive if and only if the interface $(z_1(\al,t),z_2(\al,t))$ is the graph of a function. In particular the significance of a turnover is that the Rayleigh-Taylor condition $\sigma>0$ for linear stability breaks down.

  In \cite{CCFGL}, we made a rigorous analysis of the full nonlinear problem and established analytic continuation of Muskat solutions to the time-varying strip $\{|\Im x|<h(t)\}$. This allowed us to construct  real-analytic Muskat solutions $z^0(x,t)$ defined for $t\in [-T,T]$ that ``turn over" at time 0. Here, we prove our Main Theorem \ref{conclusion} by constructing Muskat solutions $z(x,t)$ analytic on a carefully chosen time-varying domain $\Omega(t)=\{|\Im x|<h(\Re x,t)\}$.

  We take as initial datum the interface $z(x,\tau)$ at some small positive time $\tau$, and solve the Muskat equation backwards in time until we reach time $-\tau^2$. We take our initial datum $z(x,\tau)$ to be a small perturbation of $z^0(x,\tau)$, where $z^0$ is the real-analytic ``turnover" solution constructed in \cite{CCFGL}.

  Thus, we obtain Muskat solutions $z$ analytic on $\Omega(t)$ and close to $z^0$ in a suitable Sobolev norm.

\begin{figure}
\centering
\includegraphics[width=0.7\textwidth]{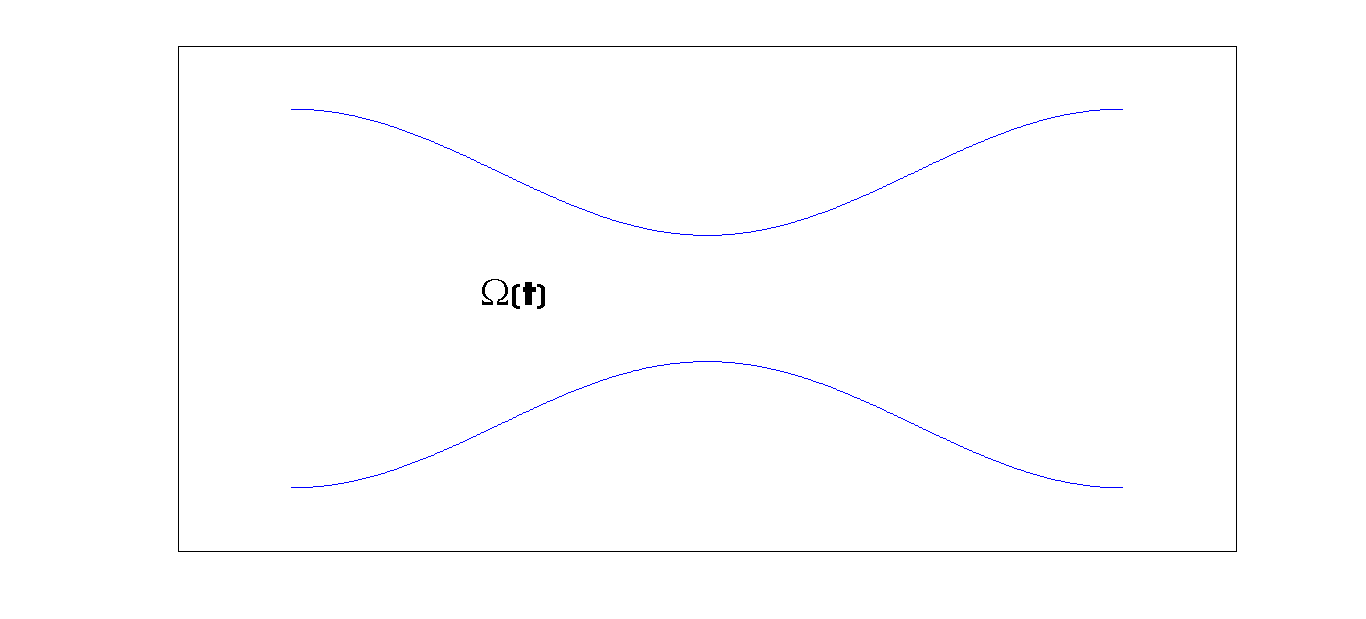}\caption{Domain $\Omega(t)$ at time  $t<\tau$.}
\label{Q(a)}
\end{figure}

\begin{figure}
\centering
\includegraphics[width=0.7\textwidth]{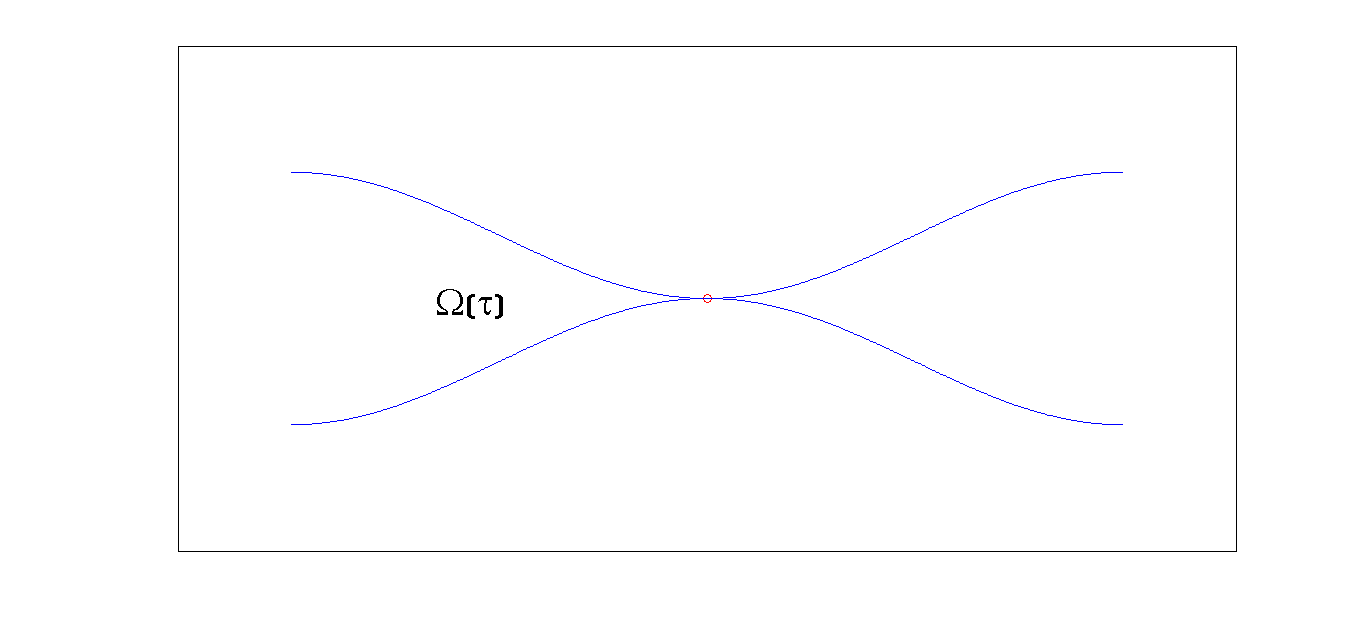}\caption{Domain $\Omega(t)$ at time $t=\tau$.}
\label{Q(b)}
\end{figure}

  For all times $t\in [-\tau^2,\tau)$, we take the domain $\Omega(t)$ to be a neighborhood of the real axis, as in Figure \ref{Q(a)}. However, at time $t=\tau$, we take $\Omega(t)$ to be a domain that pinches to a point as shown in Figure \ref{Q(b)}.


In particular an initial datum $z(x,\tau)$ may continue analytically to $\Omega(\tau)$, yet fail to belong to $C^4$ as a function of $x\in\R$. We take $z(x,\tau)$ to belong to $C^3$ but not $C^4$. Solving the Muskat equation (backwards in time) on $\Omega(t)$ $(t\in [-\tau^{2},\tau])$ for such an initial datum, we obtain a Muskat solution as in the statement of our Main Theorem \ref{conclusion}. In particular, assertions (A) and (B) of the Main Theorem hold because $z$ is a small perturbation of $z^0$, which is known to be a graph at time $t_0=-\tau^{2}$, but not at time $t_1=\tau^{2}.$ (Recall that $z^0$ turns over at time zero.) Assertion (C) holds because $z(x,t)$ continues analytically to $\Omega(t)$, with $\Omega(t)$ as in Figure \ref{Q(a)} for $t$ less than $t_2=\tau$. Finally, assertion (D) holds because we picked our initial datum $z(\cdot,\tau)$ to be in $C^3$ but not in $C^4$.

This concludes our summary of the proof of our main theorem. We warn the reader that the summary is somewhat oversimplified. The correct version of our proofs appears in the sections below.

\section{Preliminary estimates}\label{seccion1}
This section is devoted to proving crucial estimates which will be used in sections 4 and 5.

\subsection{Basic Properties of $\Lambda$}
Let $\T=\R/2\pi\Z$ and let $f\in C^2(\T)$. For $x\in \T$ we define
\begin{equation}\label{ladefi}
\Lambda f(x)=-\frac{1}{\pi}\int_{y\in\T}\frac{1}{2}\cot\left(\frac{x-y}{2}\right)(f'(x)-f'(y))dy.
\end{equation}

For $x\in \T$, let
\begin{equation}\label{distance}
||x||=\text{distance}(x,2\pi\Z).
\end{equation}
Then a carefully justified integration by parts yields
\begin{equation}\label{lambda}
\la f(x)=\frac{1}{4\pi}\lim_{\ep\to 0}\int\limits_{\substack{y\in\T\\||x-y||>\ep}}\csc^2\left(\frac{x-y}{2}\right)(f(x)-f(y))dy
\end{equation}
where the limit converges uniformly in $x$.
Consequently,
$$\int_{x\in \T}\overline{f}(x)\la f(x)dx=\frac{1}{4\pi}\lim_{\ep \to 0}\iint\limits_{\substack{x,y\in\T\\||x-y||>\ep}}\overline{f}(x)\csc^2\left(\frac{x-y}{2}\right)(f(x)-f(y))dydx$$
$$=\frac{1}{8\pi}\lim_{\ep \to 0}\iint\limits_{\substack{x,y\in\T\\||x-y||>\ep}}(\overline{f}(x)-\overline{f}(y))\csc^2\left(\frac{x-y}{2}\right)(f(x)-f(y))dydx$$
i.e.,
\begin{equation}\label{flaf}
\int_{x\in \T}\overline{f}(x)\la f(x)dx=\frac{1}{8\pi}\int\int_{x,y\in\T}\csc^2\left(\frac{x-y}{2}\right)|f(x)-f(y)|^2dydx.\end{equation}
More generally, let $A\in C^2(\T)$, and suppose that $|A''| \leq C$ on $T$. We write $C'$, $C''$, etc. to denote constants depending only on $C$.

Then  the expression (\ref{lambda}) gives
$$\int_{x\in\T} A(x) \overline{f}(x)\la f(x)dx=\frac{1}{4\pi}\lim_{\ep\to 0}\iint\limits_{\substack{x,y\in \T\\||x-y||>\ep}}A(x)\overline{f}(x)\csc^2\left(\frac{x-y}{2}\right)(f(x)-f(y))dydx$$
$$=\frac{1}{8\pi}\lim_{\ep\to 0}\iint\limits_{\substack{x,y\in \T\\||x-y||>\ep}}(A(x)\overline{f}(x)-A(y)\overline{f}(y))\csc^2\left(\frac{x-y}{2}\right)(f(x)-f(y))dydx$$
$$=\frac{1}{8\pi}\lim_{\ep\to 0}\iint\limits_{\substack{x,y\in \T\\||x-y||>\ep}}(A(x)(\overline{f}(x)-\overline{f}(y))+(A(x)-A(y))\overline{f}(y))\csc^2\left(\frac{x-y}{2}\right)(f(x)-f(y))dydx$$
$$=\frac{1}{8\pi}\iint_{x,y\in\T}A(x)\csc^2\left(\frac{x-y}{2}\right)|f(x)-f(y)|^2dydx.$$
\begin{equation}\label{flafgeneral}
+\frac{1}{8\pi}\int_{y\in\T}\overline{f}(y)\left\{\int_{x\in \T}(A(x)-A(y))\csc^2\left(\frac{x-y}{2}\right)(f(x)-f(y))dx\right\}dy\end{equation}
The second term on the right in \eqref{flafgeneral} is bounded by $C||f||_{L^2(\T)}^2$, by Calder\'on's theorem on the commutator integral (see \cite{Eli}). Therefore \eqref{flafgeneral} gives
$$\int_{x\in \T}A(x)\overline{f}(x)\la f(x) dx$$
\begin{equation}\label{six}=\frac{1}{8\pi}\iint\limits_{x,y\in
\T}A(x)\csc^2\left(\frac{x-y}{2}\right)|f(x)-f(y)|^2dxdy+Error,\end{equation}
where
$$|Error|\leq C'||f||^2_{L^2(\T)}$$
provided $|A'|\leq C$ on $\T$.

\begin{rem}Note that the formula (\ref{six}) holds for
complex-valued $A$ and $f$.
\end{rem}

\subsection{Contour Integrals Analogous to $\Lambda$}\label{seccion3}

Let $h$ be a positive function in $C^2(\T)$ ($h>0$), such that $|h|$,
$|h'|$, $|h''|$ and $|h'''|\leq C$. In this section, $C'$, $C''$ etc. denote constants depending only on $C$.

Let $\Gamma_+=\{x+ih(x)\,:\, x\in \T\}$ and
$\Gamma_-=\{x-ih(x)\,:\, x\in \T\}$. Here, both $\Gamma_+$ and
$\Gamma_-$ are contours in $\T+i\R$, oriented so that x increases as $x\pm ih(x)$ moves along $\Gamma_{\pm}$.

Let $F(z)$ be $2\pi-$periodic in $z$, and analytic in the region between $\Gamma_+$ and
$\Gamma_-$, and suppose $F$ extends smoothly to the closure of
that region. We write $f_+(x)=F(x+ih(x))$ and $f_-(x)=F(x-ih(x))$.

We study the contour integral

\begin{equation}\label{dosuno}\Lambda_{\Gamma_+}F(z)=-\frac{1}{\pi}\int\limits_{w\in\Gamma_+}\frac{1}{2}\cot\left(\frac{z-w}{2}\right)(F'(z)-F'(w))dw,\end{equation}
for $z\in\Gamma_+$.

We write $z=x+ih(x)$ and $w=u+ih(u)$. Hence, in the expression
(\ref{dosuno}) we have $F'(z)=(1+ih'(x))^{-1}f'_+(x)$ and
$F'(w)=(1+ih'(u))^{-1}f'_+(u)$. Also, $dw=(1+ih'(u))du$.

Next, we show the following lemma:
\begin{lemma}\label{dosdos}
The function
$$A(x,u)=\frac{1}{2}\cot\left(\frac{(x+ih(x))-(u+ih(u))}{2}\right)-\frac{1}{2}(1+ih'(u))^{-1}\cot\left(\frac{x-u}{2}\right)$$
is a $C^1$ function of $(x,u)\in \T\times\T$ with $C^1$ norm $\leq
C'$.
\end{lemma}
Proof: Since $A(x,u)$ is $2\pi-$periodic in $u$ for fixed $x$, we
may suppose that $|x-u|\leq\pi$. In the region $c\leq|x-u|\leq\pi$
(for small $c$ to be chosen below), we have that $\cot((x-u)/2)$
is $C^1$ with $C^1$ norm $\leq C'$, and same holds for
$\cot(((x+ih(x))-(u+ih(u)))/2)$ since
$$ dist((x+ih(x))-(u+ih(u)),2\pi\Z)\geq dist(x-u,2\pi\Z)\geq c.$$
In the region $|x-u|<c$ for small enough $c$, we have
$$|(x+ih(x))-(u+ih(u))|<\frac{1}{10},$$
and therefore
$$\frac{1}{2}\cot\left(\frac{(x+ih(x))-(u+ih(u))}{2}\right)-\frac{1}{(x+ih(x))-(u+ih(u))}=g(x,u),$$
where the $C^1$ norm of $g(x,u)$ is bounded a priori in the region
we are considering. Since also $\cot((x-u)/2)/2$ differs from
$1/(x-u)$ in $|x-u|<c$ by an error whose $C^1$ norm is bounded a
priori, lemma (\ref{dosdos}) follows easily.\\

Substituting the above results into the expression (\ref{dosuno}),
we learn that

$$\la_{\Gamma_+}F(x+ih(x))=-\frac{1}{\pi}\int_{u\in\T}\left\{(1+ih'(u))^{-1}\frac{1}{2}\cot\left(\frac{x-u}{2}\right)
+A(x,u)\right\}$$
$$\times\left\{(1+ih'(x))^{-1}f'_+(x)-(1+ih'(u))^{-1}f'_+(u)\right\}(1+ih'(u))du$$
$$=-\frac{1}{\pi}\int_{u\in\T}\frac{1}{2}\cot\left(\frac{x-u}{2}\right)
\left\{(1+ih'(x))^{-1}f'_+(x)-(1+ih'(u))^{-1}f'_+(u)\right\}du$$
$$-\frac{1}{\pi}\int_{u\in\T}(1+ih'(u))A(x,u)
\left\{(1+ih'(x))^{-1}f'_+(x)-(1+ih'(u))^{-1}f'_+(u)\right\}du$$
$$=(1+ih'(x))^{-1}\left\{-\frac{1}{\pi}\int_{u\in\T}\frac{1}{2}\cot\left(\frac{x-u}{2}\right)(f'_+(x)-f'_+(u))dy\right\}$$
$$-\frac{1}{2\pi}\int_{u\in\T}\left\{((1+ih'(x))^{-1}-(1+ih'(u))^{-1})\cot\left(\frac{x-u}{2}\right)\right \}f'_+(u)du$$
$$-\frac{1}{\pi}(1+ih'(x))^{-1}f'_+(x)\int_{u\in\T}(1+ih'(u))A(x,u)du+\frac{1}{\pi}\int_{u\in\T}A(x,u)f'_+(u)du.$$
On the right hand side, the first term is $(1+ih'(x))^{-1})\la
f_+(x)$; the expression in curly brackets in the second term has
$C^1$ norm bounded a priori, thus allowing us to integrate by
parts. We can perform a similar argument with the last term.
Therefore
$$\la_{\Gamma_+}F(x+ih(x))=(1+ih'(x))^{-1}\la f_+(x)$$
\begin{equation}\label{dostres}-\frac{1}{\pi}(1+ih'(x))^{-1}f'_+(x)
\int_{u\in\T}(1+ih'(u))A(x,u)du+Error(x),\end{equation}
with
$$|Error(x)|\leq C ||f_+||_{L^2(\T)}.$$
We study the $A-$integral in the expression (\ref{dostres}). By
definition,
$$A(x,u)=\frac{1}{2}\cot\left(\frac{z-w}{2}\right)-(1+ih'(u))^{-1}\frac{1}{2}\cot\left(\frac{x-u}{2}\right),$$
with $z=x+ih(x)$ and $w=u+ih(u)$. Moreover, $dw=(1+ih'(u))du$.
Therefore,
$$\int_{u\in \T} A(x,u)(1+ih'(u))du=$$
$$P.V.\int_{w\in
\Gamma_+}\frac{1}{2}\cot\left(\frac{z-w}{2}\right)dw-P.V.\int_{u\in
\T}\frac{1}{2}\cot\left(\frac{x-u}{2}\right)du,$$ with
$z=x+ih(x)$. Since
$$P.V.\int_{u\in
\T}\frac{1}{2}\cot\left(\frac{x-u}{2}\right)du=0,$$ we have that
\begin{equation}\label{doscuatro}\int_{u\in\T}A(x,u)(1+ih'(u))du=P.V.\int_{w\in\Gamma_+}\cot\left(\frac{z-w}{2}\right)dw=0.\end{equation}
The last equality in \eqref{doscuatro} is an exercise in contour integration.

Therefore, the expressions (\ref{dostres}) and (\ref{doscuatro})
yield the following result,
\begin{equation}\label{dossiete}
\la_{\Gamma_+}F(x+ih(x))=(1+ih'(x))^{-1}\la f_+(x)+Error_+(x),
\end{equation}
where
$$||Error_+||_{L^2(\T)}\leq C||f_+||_{L^2(\T)}.$$
Similarly,
\begin{equation}\label{dosocho}\Lambda_{\Gamma_-}F(z)=-\frac{1}{\pi}\int_{w\in\Gamma_-}
\frac{1}{2}\cot\left(\frac{z-w}{2}\right)(F'(z)-F'(w))dw,\end{equation}
for $z\in\Gamma_-$  satisfies
\begin{equation}\label{dosnueve}
\Lambda_{\Gamma_-}F(z)=
(1-ih'(x))^{-1}\la f_-(x)+Error_-(x),\end{equation}
where $$||Error_-||_{L^2(\T)}\leq C||f_-||_{L^2(\T)}.$$
Next we return to the definition (\ref{dosuno}), and shift the
contour from $\Gamma_+$ to $\Gamma_-$ (note that the integrand in
the expression (\ref{dosuno}) has no poles).

Thus,
$$\la_{\Gamma_+}F(x+ih(x))=-\frac{1}{\pi}\int\limits_{w\in\Gamma_-}\frac{1}{2}\cot\left(\frac{z-w}{2}\right)(F'(z)-F'(w))dw$$
$$=F'(z)\left\{-\frac{1}{\pi}\int_{\Gamma_-}\frac{1}{2}\cot\left(\frac{z-w}{2}\right)dw\right\}+
\frac{1}{2\pi}\int_{w\in\Gamma_-}\cot\left(\frac{z-w}{2}\right)F'(w)dw$$
for $z=x+ih(x)$.

Another exercise in contour integration shows that

$$\int_{\Gamma_-}\frac{1}{2}\cot\left(\frac{z-w}{2}\right)dw=-i\pi,$$

and therefore,
$$\la_{\Gamma_+}F(x+ih(x))=iF'(x+ih(x))+\frac{1}{2\pi}\int_{w\in\Gamma_-}\cot\left(\frac{z-w}{2}\right)F'(w)dw$$
\begin{equation}\label{dosdiez}=iF'(x+ih(x))-\frac{1}{4\pi}\int_{w\in\Gamma_-}\csc^2\left(\frac{z-w}{2}\right)F(w)dw.\end{equation}
For $z=x+ih(x)\in \Gamma_+$, $w=u-ih(u)\in \Gamma_-$, we have
$$=\left|\csc^2\left(\frac{z-w}{2}\right)\right|\leq \frac{C}{||\Re{(z-w)}||^2+|\Im{(z-w)}|^2}=\frac{C}{||x-u||^2+(h(x)+h(u))^2}$$
$$\leq \frac{C}{||x-u||^2+h(x)^2},\quad |F(w)|=|f_-(u)|,\,\,\, |dw|\leq Cdu.$$
Hence,
$$\left|\int_{w\in\Gamma_-}\csc^2\left(\frac{z-w}{2}\right)F(w)dw\right|\leq C\int_{u\in\T}\frac{|f_-(u)|}{||x-u||^2+h(x)^2}du$$
$$\leq C h^{-1}(x)M[f_-](x),$$
where $M[f_-]$ denotes the Hardy-Littlewood maximal function of $f_-$. (See \cite{Eli}, page 9).

Together with (\ref{dosdiez}), this implies:
\begin{equation}\label{dosonce}\la_{\Gamma_+}F(x+ih(x))=iF'(x+ih(x))+Erreur_+(x),\end{equation}
where
$$||h(\cdot)Erreur_+(\cdot)||_{L^2(\T)}\leq C||f_-||_{L^2(\T)}.$$
Similarly,
\begin{equation}\label{dosdoce}\la_{\Gamma_-}F(x-ih(x))=-iF'(x-ih(x))+Erreur_-(x),\end{equation}
where
$$||h(\cdot)Erreur_-(\cdot)||_{L^2(\T)}\leq C||f_+||_{L^2(\T)}.$$
From (\ref{dossiete}) and (\ref{dosonce}), we see that
\begin{equation}\label{dostrece}
iF'(x+ih(x))=(1+ih'(x))^{-1}\la f_+(x)+ Goof_+(x)\end{equation}
where
$$||h(\cdot)Goof_+(\cdot)||_{L^2(\T)}\leq C||f_+||_{L^2(\T)}+C||f_-||_{L^2(\T)}.$$
Similarly
\begin{equation}\label{doscatorce}
-iF'(x-ih(x))=(1-ih'(x))^{-1}\la f_-(x)+ Goof_-(x)\end{equation}
where
$$||h(\cdot)Goof_-(\cdot)||_{L^2(\T)}\leq C||f_+||_{L^2(\T)}+C||f_-||_{L^2(\T)}.$$

\subsection{Applying results derived above}\label{seccion4}

In this section, we are given two functions $h$, $h_t\in C^2(\T)$ (with $C^3$ norms assumed bounded a priori). We again assume that $h>0$ and we take the
contours $\Gamma_+=\{x+ih(x)\,:\,x\in\T\}$, $\Gamma_-=\{x-ih(x)\,:\,x\in\T\}$, and let $F(z)$
be a $2\pi-$periodic  analytic function on the region between $\Gamma_+$ and $\Gamma_-$, smooth on the closure
of that region. We take $f_+(x)=F(x+ih(x))$ and $f_-(x)=F(x-ih(x))$ for $x\in\T$.

Let $\theta$ be a function on $\T$, with $C^0-$norm assumed bounded a priori.

We assume that
\begin{equation}\label{tresuno}
|h_t(x)|\leq Sh(x)\end{equation}
for all $x\in supp(\theta)$, where $S$ is a given positive number.

Let $\tilde{\sigma}_1(x)$ be a function on $\T$, assumed to have a priori bounded $C^2$ norm. We want to estimate
$$X\equiv\Re\,\int_{x\in\T}\overline {f}_+(x)\theta(x)\left\{-\frac{\tilde{\sigma}_1(x)}{\pi}\int_{u\in\T}\frac{1}{2}\cot\left(\frac{(x+ih(x))-(u+ih(u))}{2}\right)\right.$$
$$\times(F'(x+ih(x))-F'(u+ih(u)))(1+ih'(u))du$$
\begin{equation}\label{tresdos}\left.
+ih_t(x)F'(x+ih(x))\frac{}{}\right\}dx,\end{equation}
from below. Applying (\ref{dossiete}) and (\ref{dostrece}), we have
$$-\frac{\tilde{\sigma}_1(x)}{\pi}\int_{u\in\T}\frac{1}{2}\cot\left(\frac{(x+ih(x))-(u+ih(u))}{2}\right)\times(F'(x+ih(x))-F'(u+ih(u)))(1+ih'(u))du$$
$$=\tilde{\sigma}_1(x)(1+ih'(x))^{-1}\Lambda f_+(x)+Err_+(x),$$
with
$$||Err_+||_{L^2(\T)}\leq C||f_+||_{L^2(\T)},$$
and
$$ih_t(x)F'(x+ih(x))=h_t(x)(1+ih'(x))^{-1}\la f_+(x)+Oops_+(x),$$
with
$$||Oops_+||_{L^2(supp(\theta))}\leq CS\left(||f_+||_{L^2(\T)}+||f_-||_{L^2(\T)}\right).$$
Therefore,
\begin{equation}\label{trestres}
X=\Re\int_{x\in\T}\overline{f}_+(x)\theta(x)(1+ih'(x))^{-1}\{\tilde{\sigma}_1(x)+h_t(x)\}\la f_+(x)dx+Discrep_+,\end{equation}
where
$$|Discrep_+|\leq C(S+1)\left(||f_+||_{L^2(\T)}^2+||f_-||_{L^2(\T)}^2\right).$$
We recall that the derivatives of $(1+ih'(x))^{-1}$ up to order two are bounded a priori. We make the following assumption:

The $C^2$ norm of $\theta(x)\{\tilde{\sigma}_1(x)+h_t(x)\}$ is bounded a priori.

Then, by (\ref{six}) we have
$$X=\frac{1}{8\pi}\iint\limits_{x,y\in\T}\Re\left\{\theta(x)(1+ih'(x))^{-1}(\tilde{\sigma}_1(x)+h_t(x))\right\}\csc^2\left(\frac{x-y}{2}\right)|f_+(x)-f_+(y)|^2dxdy$$
\begin{equation}\label{trescinco}+Errore_+\end{equation}
with
$$|Errore_+|\leq C(S+1)\left(||f_+||_{L^2(\T)}^2+||f_-||_{L^2(\T)}^2\right).$$
If
\begin{equation}\label{tresseis}
\Re\left\{\theta(x)(1+ih'(x))^{-1}(\tilde{\sigma}_1(x)+h_t(x))\right\}\geq 0,
\end{equation}
then (\ref{trescinco}) yields
\begin{equation}\label{tressiete}
X\geq -C(S+1)\left(||f_+||_{L^2(\T)}^2+||f_-||_{L^2(\T)}^2\right).\end{equation}
Thus we have proven the following result.
\begin{lemma}\label{lemma2}
Let $\theta$, $\tilde{\sigma}_1$, $h$, $h_t$ be functions on $\T$. Let $S$ be a positive number. We make the following assumptions
\begin{enumerate}
\item $h>0$ on $\T$. \item $|h_t|\leq Sh$ on $\supp(\theta)$. \item
The $C^3$ norms of $h$ and $h_t$ are bounded a priori. \item The
$C^2$ norm of $\tilde{\sigma}{_1}$ is bounded a priori. \item The
$C^2$ norm of $\theta(x)(1+ih'(x))^{-1}(\tilde{\sigma}_1(x)+h_t(x))$
is bounded a priori. \item
$\Re\left\{\theta(x)(1+ih'(x))^{-1}(\tilde{\sigma}_1(x)+h_t(x))\right\}\geq
0$ on $\T$.
\item The $C^0$-norm of $\theta$ is bounded a priori.
\end{enumerate}
Let $F(z)$ be analytic on $\{z\in\C\,:\,|\Im(z)|\leq h(\Re(z))\}$, $2\pi-$periodic and smooth on the closure of the above region. Let
\begin{equation*}
f_+(x)=F(x+ih(x))
\end{equation*}
and
\begin{equation*}
f_-(x)=F(x-ih(x)),
\end{equation*}
for $x\in\T$. Then the following inequality holds:
$$\Re\,\int_{x\in\T}\overline {f}_+(x)\theta(x)\left\{-\frac{\tilde{\sigma}_1(x)}{\pi}\int_{u\in\T}\frac{1}{2}\cot\left(\frac{(x+ih(x))-(u+ih(u))}{2}\right)\right.$$
$$\times(F'(x+ih(x))-F'(u+ih(u)))(1+ih'(u))du$$
$$
\left.+ih_t(x)F'(x+ih(x))\frac{}{}\right\}dx$$ $$\geq -C(S+1)\left(||f_+||_{L^2(\T)}^2+||f_-||_{L^2(\T)}^2\right).$$
Similarly, under assumptions analogous to 1,...,7, with $1-ih'(x)$ in place of $1+ih'(x)$ we have
$$\Re\,\int_{x\in\T}\overline {f}_-(x)\theta(x)\left\{-\frac{\tilde{\sigma}_1(x)}{\pi}\int_{u\in\T}\frac{1}{2}\cot\left(\frac{(x-ih(x))-(u-ih(u))}{2}\right)\right.$$
$$\times(F'(x-ih(x))-F'(u-ih(u)))(1-ih'(u))du$$
$$
\left.-ih_t(x)F'(x-ih(x))\frac{}{}\right\}dx$$ $$\geq-C(S+1)\left(||f_+||_{L^2(\T)}^2+||f_-||_{L^2(\T)}^2\right).$$
Here, the constant $C$ depends only on the a priori bounds assumed in 3, 4, 5 and 7 above.
\end{lemma}

\subsection{G{\aa}rding Inequalities}

First we work on $\R$. Later on, we will switch over to $\T$.

Here $\Lambda_\R$ will be the square root of $(-\Delta)$ in one
dimension, $D=\frac{1}{i}\frac{d}{dx}$ on $\R$ and
$$\Pi_\sigma\,:\, L^2(\R)\to L^2(\R)$$
for $\sigma=\pm 1$ is the orthogonal projection onto either
positive frequencies ($\sigma=1$) or negative frequencies
($\sigma=-1$). Note $\Pi_\sigma=\frac{1}{2}\left(I+i\sigma
H\right)$, where $H$ is the Hilbert transform ($\hat{H}=-i\, sign (\xi)$).

Commutator estimates:

Let $a$ be a real-valued function on $\R$. We suppose $a'\in L^\infty(\R)$, and write $a$ for the operator $f\mapsto
a\cdot f$ on $L^2(\R)$. Then
$$[a,H]D\,:\, L^2(\R)\to L^2(\R)$$
is bounded, with norm at most $C||a'||_{L^\infty(\R)}$, where $C$
is a universal constant.

That follows from Calder\'on's theorem on the commutator
integral (see \cite{Eli}).

Since $\Lambda_{\R}=iDH$, it follows that

$$[a,H]\Lambda_\R\,:\,L^2(\R) \to L^2(\R)$$
is bounded, with norm at most $C||a'||_{L^\infty(\R)}$.

Since $\Pi_\sigma=\frac{1}{2}(I+i\sigma H)$, it follows that
$$||[a,\Pi_\sigma]D||,\,||[a,\Pi_\sigma]\Lambda_{\R}||\leq
C||a'||_{L^\infty(\R)}.$$

Next, recall that
$$2\Re\left(\overline{f}\cdot\Lambda_\R f\right)\geq \Lambda_\R
\left(|f|^2\right)$$ pointwise (see \cite{AD}).

If $a\geq0$ everywhere and $\Lambda a\in L^\infty$, then
\begin{align*}\Re\left\{\int_\R a(x)\overline{f}(x)\Lambda_\R f(x)dx\right\}&\geq\frac{1}{2}\int_\R a(x)
\left(\Lambda_\R\left(|f|^2\right)\right)dx \\
=\frac{1}{2}\int_\R\left(\Lambda_\R a\right)(x)|f(x)|^2dx &
\geq-\frac{1}{2}||\Lambda_\R a||_{L^\infty}||f||_{L^2}^2.
\end{align*}
Now suppose $a$, $b$ are real, $a'$, $b'\in L^\infty$; $\Lambda_\R
a$, $\Lambda_\R b\in L^\infty$; $a\geq|b|$ pointwise.

Note that $\Lambda_{\R}=\sum_{\sigma=\pm
1}\Pi_{\sigma}\Lambda_\R\Pi_\sigma$ and $D=\sum_{\sigma=\pm
1}\Pi_\sigma D \Pi_\sigma$.

Therefore,
\begin{align*}
a\Lambda_{\R}+b D=& \sum_{\sigma=\pm
1}\left\{[a,\Pi_\sigma]\Lambda_{\R}\Pi_\sigma+[b,\Pi_\sigma]D\Pi_\sigma\right\}\\
&+
\sum_{\sigma}\Pi_\sigma(a\Lambda_{\R}+bD)\Pi_\sigma.\end{align*}

The terms in curly brackets are operators with norm at most
$C||a'||_{L^\infty}+C||b'||_{L^\infty}$ on $L^2(\R)$, with $C$ a universal constant.

Therefore,
\begin{align*}
\Re
\langle(a\Lambda_{\R}+bD)f,f\rangle=&\text{Error}+\sum_{\sigma}\Re\langle
\Pi_\sigma (a\Lambda_\R+bD)\Pi_{\sigma}f,f\rangle\\
 =& \text{Error}+\sum_{\sigma}\Re
 \langle(a\Lambda_\R+bD)\Pi_\sigma f,\Pi_\sigma f\rangle\\
 =& \text{Error}+\sum_{\sigma}\Re \langle (a+\sigma
 b)\Lambda_{\R}\Pi_\sigma f,\Pi_\sigma f\rangle,
\end{align*}
since $D\Pi_\sigma=\sigma \Lambda_{\R}\Pi_\sigma$; here
$$|\text{Error}|\leq C
(||a'||_{L^\infty}+||b'||_{L^\infty})||f||_{L^2}^2.$$ Recalling
that $a+\sigma b\geq 0$ on $\R$, we see that
$$\Re\langle (a\Lambda_{\R}+bD)f,f\rangle\geq
-C\left\{||a'||_{L^\infty}+||b'||_{L^\infty}+||\Lambda_\R
a||_{L^\infty}+||\Lambda_{\R}b||_{L^\infty}\right\}||f||_{L^2}^2,$$
where $C$ is a universal constant.

Now suppose $a,b$ real, $a\geq|b|$ on $\R$,
$|a|,|b|,|a'|,|b'|,|a''|,|b''|\leq 1$. Then
$$\Lambda_\R
a(x)=\frac{1}{\pi}P.V.\int_{\R}\frac{a(x)-a(y)}{(x-y)^2}dy.$$ Now
$$\left|\int_{|x-y|>1}\frac{a(x)-a(y)}{(x-y)^2}dy\right|\leq \int_{|x-y|>1}\frac{|a(x)|+|a(y)|}{(x-y)^2}dy
\leq C,$$ and
$$\left|P.V.\int_{|x-y|\leq
1}\frac{a(x)-a(y)}{(x-y)^2}dy\right|=\left|P.V.\int_{|x-y|\leq
1}\frac{a(x)-a(y)+a'(x)(x-y)}{(x-y)^2}dy\right|\leq C.$$ Similarly
for $b$. Therefore, $a'$, $b'$, $\Lambda_\R a$, $\Lambda_\R b$ are
bounded on $\R$.

Therefore, we obtain the following result:

Suppose $a,b$ are real functions on $\R$, such that
$|a|,|b|,|a'|,|b'|,|a''|,|b''|\leq 1$; and suppose $a\geq |b|$
everywhere on $\R$. Then

$$\Re \langle(a\Lambda_\R+bD)f,f\rangle\geq -C||f||_{L^2}^2$$
for $f\in C^2_0(\R)$, where $C$ is a universal constant.

Rescaling, we obtain the following result:
\begin{lemma}\label{5.0}
Let $0<\delta<1$, and let $a,b$ be real functions on $\R$, such that
$|a|\leq \delta$, $|b|\leq \delta$, $|a'|\leq 1$, $|b'|\leq1$,
$|a''|\leq \delta^{-1}$, $|b''|\leq \delta^{-1}$ everywhere, and
also $a\geq |b|$ everywhere. Then
$$\Re \langle(a\Lambda_\R+bD)f,f\rangle\geq -C||f||_{L^2}^2$$
for $f\in C^2_0(\R)$, where $C$ is a universal constant.
\end{lemma}
Next, we transfer our result from $\R$ to $\T$.
\begin{lemma} Let $0<\delta<1$, and let $a,b$ be real-valued
functions on $\T$. Suppose that
$$\left|\left(\frac{d}{dx}\right)^j a(x)\right|, \left|\left(\frac{d}{dx}\right)^j
b(x)\right|\leq \delta^{1-j}$$ for $j=0,1,2$ and $x\in\T$.

Suppose also that $|b|\leq a$ for all $x\in\T$. Then, for any
smooth function $f$ on $\T$, we have
\begin{equation}\label{asterisco}
\Re\left\{\int_{\T}\overline{f}(x)\left(a(x)\Lambda
f(x)+b(x)if'(x)\right)dx\right\}\geq -C
\int_{\T}|f(x)|^2dx,\end{equation} where $C$ is a universal
constant.
\end{lemma}

 Proof: Suppose first that $supp(a)$, $supp(b)\subset
[-1,1]$. We identify $\T$ with $[-\pi,\pi]$, and let $f=0$ on
$\R\setminus[-\pi,\pi]$. Thus, we may regard $f$ as a function on
$\R$. For $x\in[-1,1]$, we have
$$\Lambda f(x)-\Lambda_{\R} f(x)=\int_\R K(x-y)f(y)dy$$
for a kernel $K$, where $||K||_{L^\infty(\R)}\leq C$.

Hence, in
proving (\ref{asterisco}), we may replace $\Lambda$ by
$\Lambda_\R$. Estimate (\ref{asterisco}) now follows from lemma
(\ref{5.0}).

We have established (\ref{asterisco}) in the case where $supp(a)$,
$supp(b)\subset I$, where $I$ is any interval of length 2 in $\T$.
Conclusion (\ref{asterisco}) in the general case now follows by
using a partition of unity.

\subsection{Applications of G{\aa}rding's Inequality}
In this section, we
are given functions $h$, $h_t\in \C^{3}(\T)$ (with norms bounded
a priori). As usual, we assume that $h>0$, and we take
$\Gamma_+=\{x+ih(x)\,:\, x\in\T\}$, $\Gamma_-=\{x-ih(x)\,:\,
x\in\T\}$. We let $F$ be analytic on the region of $\T+i\R$ bounded
by $\Gamma_+$ and $\Gamma_-$ and we suppose that $F$ is smooth on
the closure of this region.
 As usual, we set $f_+(x)=F(x+ih(x))$, $f_-(x)=F(x-ih(x))$.

Let $\tilde{\sigma}_1(x)$ and $\theta(x)$ be functions on $\T$. As in  section (\ref{seccion4}),
we want a lower bound for
$$X\equiv\Re\,\int_{x\in\T}\overline {f}_+(x)\theta(x)\left\{-\frac{\tilde{\sigma}_1(x)}
{\pi}\int_{u\in\T}\frac{1}{2}\cot\left(\frac{(x+ih(x))-(u+ih(u))}{2}\right)\right.$$
$$\times(F'(x+ih(x))-F'(u+ih(u)))(1+ih'(u))du$$
\begin{equation}\label{seisuno}\left.
+ih_t(x)F'(x+ih(x))\frac{}{}\right\}dx.\end{equation}
We assume that the $C^2$ norm of $\tilde{\sigma}_1$ is bounded a priori and that $0\leq\theta\leq 1$ on $\T$. According to (\ref{dosuno}) and (\ref{dossiete}) we have
$$X=\Re \int_{x\in\T}\overline{f}_+(x)\theta(x)\left\{\tilde{\sigma}_1(x)(1+ih'(x))^{-1}\la f_+(x)+ih_t(x)F'(x+ih(x))\right\}dx+ \check{E}rr,$$
with
$$|\check{E}rr|\leq C||f_+||_{L^2(\T)}^2.$$
Also,
$$F'(x+ih(x))=(1+ih(x))^{-1}f_+'(x),$$
as noted before. Therefore,
\begin{equation}\label{seisdos}X=\Re \int_{x\in\T}\overline{f}_+(x)\theta(x)(1+ih'(x))^{-1}\left\{\tilde{\sigma}_1(x)\la f_+(x)+ih_t(x)f_+'(x)\right\}dx+ \check{E}rr,\end{equation}
with
$$|\check{E}rr|\leq C||f_+||_{L^2(\T)}^2.$$
We now suppose that:
\begin{equation}\label{seistres}\text{The $C^2$ norms of $\theta(x)\tilde{\sigma}_1(x)$ and $\theta(x)h_t(x)$ are bounded a priori.}\end{equation}
Then the same holds for $\theta(x)(1+ih'(x))^{-1}\tilde{\sigma}_1(x)$ and $\theta(x)(1+ih'(x))^{-1}h_t(x)$.

Therefore, formula (\ref{six}) gives
$$\Re\int_{x\in\T}\overline{f}_+(x)\theta(x)(1+ih'(x))^{-1}\tilde{\sigma}_1(x)\Lambda f_+dx$$
\begin{equation}\label{seiscuatro}
=\Re\int_{x\in\T}\overline{f}_+(x)\Re\left\{\theta(x)(1+ih'(x))^{-1}\tilde{\sigma}_1(x)\right\}\la f_+(x)dx+\check{\check{E}}rr,\end{equation}
with
$$|\check{\check{E}}rr|\leq C||f_+||_{L^2(\T)}^2.$$
Also, if $iB(x)$ is real, then
$$\Re\int_{x\in\T}\overline{f}_+(x)B(x)if'_+(x)dx=\int_{x\in\T}iB(x)\Re(\overline{f}_+(x)f'_+(x))$$
$$=\int_{x\in\T}iB(x)\frac{1}{2}\frac{d}{dx}|f_+(x)|^2dx=-\frac{i}{2}\int_{x\in\T}B'(x)|f_+(x)|^2dx.$$
Applying this with $B(x)=i\Im\{\theta(x)(1+ih'(x))^{-1}h_t(x)\}$, we see that
$$\Re\int_{x\in\T}\overline{f}_+(x)\theta(x)(1+ih'(x))^{-1}h_t(x)if'_+(x)dx$$
\begin{equation}\label{seiscinco}
=\Re\int_{x\in\T}\overline{f}_+(x)\Re\{\theta(x)(1+ih'(x))^{-1}h_t(x)\}if'_+(x)dx + \check{\check{\check{E}}}rr,\end{equation}
with
$$|\check{\check{\check{E}}}rr|\leq C||f_+||_{L^2(\T)}^2.$$
Putting (\ref{seiscuatro}) and (\ref{seiscinco}) into (\ref{seisdos}) we see that
\begin{equation}\label{seissiete}
X=\Re\int_{x\in\T} \overline{f}_+(x)\left\{A(x)\la f_+(x)+B(x)if'_+(x)\right\}dx+Err,\end{equation}
where
$$|Err|\leq C||f_+||_{L^2(\T)}^2$$
and
\begin{equation}\label{seisocho}A(x)=\Re\{\theta(x)(1+ih'(x))^{-1}\tilde{\sigma}_1(x)\},\end{equation}
\begin{equation}\label{seisnueve}B(x)=\Re\{\theta(x)(1+ih'(x))^{-1}h_t(x))\}.\end{equation}
We assume that
\begin{equation}
\left|\left(\frac{d}{dx}\right)^jA(x)\right|\,,\,\left|\left(\frac{d}{dx}\right)^jB(x)\right|\leq C\delta^{1-j}\end{equation}
for $0\leq j\leq 2$;
and that
\begin{equation}\label{seisdoce}
|B(x)|\leq A(x)\quad \text{for all $x\in\T$}.\end{equation}
From (\ref{seissiete}) and  G{\aa}rding's inequality we obtain the estimate
\begin{equation}\label{seistrece}
X\geq -C||f_+||_{L^2(\T)}^2.
\end{equation}
Thus, we have proven the following result:
\begin{lemma}\label{lemma3}
Let $\theta$, $h$, $h_t$, $\tilde{\sigma}_1$ be functions on $\T$, and let $0<\delta< 1$ be a number. We make the following assumptions:
\begin{enumerate}
\item $h>0$. \item The norm of $\tilde{\sigma}_1$ in
$C^{2}$  and the norms of $h$ and $h_t$ in
$C^{3}$ are bounded a priori.
\item The $C^2$ norms of $\theta(x)\tilde{\sigma}_1(x)$ and
$\theta(x)h_t(x)$ are bounded a priori.
\item $supp\,\, \theta \subset [-\delta,\delta]$.
\item
$\left|\left(\frac{d}{dx}\right)^j\Re\{\theta(x)(1+ih'(x))^{-1}\tilde{\sigma}_1(x)\}\right|$,
$\left|\left(\frac{d}{dx}\right)^j\Re\{\theta(x)(1+ih'(x))^{-1}h_t(x)\}\right|$
$\leq C\delta^{1-j}$

for $0\leq j\leq 2$, $x\in\T$. \item
$|\Re\{\theta(x)(1+ih'(x))^{-1}h_t(x)\}|\leq
\Re\{\theta(x)(1+ih'(x))^{-1}\tilde{\sigma}_1(x)\}$ for all
$x\in\T$.
\item $0\leq\theta\leq 1$.
\end{enumerate}
Then
$$\Re\,\int_{x\in\T}\overline {f}_+(x)\theta(x)\left\{-\frac{\tilde{\sigma}_1(x)}{\pi}\int_{u\in\T}\frac{1}{2}\cot\left(\frac{(x+ih(x))-(u+ih(u))}{2}\right)\right.$$
$$\times(F'(x+ih(x))-F'(u+ih(u)))(1+ih'(u))du$$
$$\left.
+ih_t(x)F'(x+ih(x))\frac{}{}\right\}dx$$
$$\geq -C'||f_+||_{L^2(\T)^2}^2,$$
where $f_+(x)=F(x+ih(x))$.

Here, $C'$ depends only on the a priori bounds in 2, 3 and on the constant $C$ in 5.
\end{lemma}
\begin{rem}An analogous result holds for $f_-(x)=F(x-ih(x))$.\end{rem}


\section[The unperturbed solution \underline{$z$}  \; and the height function $h$]{The unperturbed solutions \underline{$z$} and the height function $h$}\label{seccion33}

In this section we study the Rayleigh-Taylor function associated to the unperturbed solutions. We construct our time-varying domain of analyticity $\Omega(t)=\{|\Im x|\leq h(\Re x,t)\}$
adapted to the Rayleigh-Taylor function.

\subsection{The Rayleigh Taylor Function}\label{seccion7}

We consider the  analytic solution
$\underline{z}(x,t)=(\underline{z}_1(x,t),\underline{z}_2(x,t))$, for $t\in [-T,T]$ ($T>0$), found in \cite{CCFGL}. We will call
$\underline{z}(x,t)$ the unperturbed solution and we recall that we can
choose the origin of time so that $\underline{z}(x,t)$ satisfies the following
properties:
\begin{itemize}
\item $\un_1(x,t)-x$ and $\un_2(x,t)$ are $2\pi-$periodic functions and real for $x$ real.
\item $\un_1(x,0)$ and $\un_2(x,0)$ are odd functions.
\item $\underline{z}_1(0,0)=0,$ \item $\pa_x \underline{z}_1(0,0)=0$,\item $\pa_x^2
\underline{z}_1(0,0)=0$, \item $\pa^3_x \underline{z}_1(0,0)>\tilde{c}_1>0$, \item $\pa_t\pa_x
\underline{z}_1(0,0)<-\tilde{c}_2<0$. \item $\pa_x \underline{z}_2(0,0)>\underline{c}>0$.
\end{itemize}
Then since $|\pa_x^4 \underline{z}_1(x,0)|\leq C_3$ it follows from Taylor's
theorem that \begin{equation}\label{capital1}c_4 x^2\leq \pa_x
\underline{z}_1(x,0)\leq C_5 x^2,\end{equation} for $|x|\leq c_6$.

Similarly, since $$\pa_t\pa_x \underline{z}_1(0,0)<-\tilde{c}_2<0$$ and
$$|\pa^2_x \pa_t \underline{z}_1(x,0)|\leq C_7$$ it follows that
$$\pa_t\pa_x \underline{z}_1(x,0)<-c_8<0$$ for $|x|<c_9$.

Since also $$|\pa_t^2 \pa_x \underline{z}_1(x,t)|\leq C_{10}$$ for $|x|\leq
c_9$, $|t|\leq c_{11}$, it follows in turn that
\begin{equation}\label{capital2}\pa_x \underline{z}_1(x,0)-C_{12}t\leq \pa_x
\underline{z}_1(x,t)\leq \pa_x \underline{z}_1(x,0)-c_{13}t\end{equation} for $|x|\leq
c_9$, $0\leq t\leq c_{14}$.

Combining our estimates (\ref{capital1}) and (\ref{capital2}), we
find that
\begin{equation}\label{capital3}
c_4x^2-C_{12}t\leq \pa_x \underline{z}_1(x,t)\leq
C_5x^2-c_{13}t\end{equation} for $|x|\leq c_{15}$, $0\leq t\leq
c_{16}$.
Similarly we have that
\begin{equation}\label{capital31/2}
c'_4x^2-c'_{12}t\leq \pa_x \underline{z}_1(x,t)\leq
C'_5x^2-C'_{13}t\end{equation} for $|x|\leq c_{15}$, $-c_{16}\leq t\leq
0$.

Also, since $\pa_x^2 \underline{z}_1(0,0)=0$ and $\underline{z}_1(x,t)$ has $C^3-$norm
at most $C$ on $\{|x|\leq c, |t|\leq c\}$, it follows that
\begin{equation}\label{capital4}
|\pa_x^2 \underline{z}_1(x,t)|\leq C_{17}|t|+ C_{17}|x|
\end{equation}
for $|x|\leq c_{18}$, $|t|\leq c_{19}$.

We have also
\begin{equation}\label{capital5}
|\pa_x^3\underline{z}_1(x,t)|\leq C_{20}\end{equation} for $|x|$, $|t|\leq
c$.
Finally we notice that
\begin{equation}\label{grafo}
\pa_x \un_1(x,-t)>c_{20}t\quad \text{for all $x\in\T$ and $0<t<c_{21}$}.
\end{equation}

Now we define the Rayleigh-Taylor function associated to the
unperturbed solution as
\begin{equation}\label{rt}
\sigma_1^0(x,t)=\frac{-2\pi \pa_x \underline{z}_1(x,t)}{(\pa_x \underline{z}_1(x,t))^2+
(\pa_x \underline{z}_2(x,t))^2}
\end{equation}
Therefore $\sigma_1^0(x,t)$ satisfies the following:
\begin{enumerate}
\item $\sigma^0_1(\cdot,t)$ is analytic on
$\{x+iy\,:\,x\in\T,\,|y|\leq c_b\}$ with $|\sigma_1^0(x+iy,t)|\leq
C$, for all $x+iy$ as above and for all $t\in[-\tau,\tau]$. \item
$\sigma_1^0(x,t)$ is real for $x\in\T$, $t\in[-\tau,\tau]$. \item
$\sigma_1^0$ has a priori bounded $C^{2}$ norm as a function of
$(x,t)\in\T\times[-\tau,\tau]$. \item
$\sigma_1^0(0,0)=0$. \item $\partial_x \sigma_1^0(0,0)=0$. \item
$\partial_x^2 \sigma_1^0(0,0)=-c_2<0$. \item $\partial_t
\sigma_1^0(0,0)=c_1>0$.
\end{enumerate}
For complex $x$ ($|\Im x|<c'$), we have
$$\sigma_1^0(x,t)=\sigma_1^0(0,t)+x\partial_x\sigma_1^0(0,t)+\frac{1}{2}\partial^2_x\sigma_1^0(0,t)x^2+O(||x||^3),$$
where $||x||=||\Re x||+|\Im x|\approx distance(x,2\pi\Z)$.

Moreover,
$$\sigma_1^0(0,t)=c_1t+O(t^2),$$
and  we have
\begin{eqnarray*}
\partial_x^2\sigma_1^0(0,t)&=&-c_2+O(t)\\
\partial_x\sigma_1^0(0,t)&=&O(t).
\end{eqnarray*}
Hence, for complex $x$ such that $|x|\leq 10^{-3}$ and $|\Im x|<c'$ we have
$$\sigma_1^0(x,t)=c_1t+O(t^2)+O(xt)-\frac{c_2}{2}x^2+O(tx^2)+O(x^3).$$
Since also $|O(xt)|\leq \frac{c_2}{1000}x^2+O(t^2)$, we conclude that
$$ |\sigma_1^0(x,t)-[c_1t-\frac{c_2}{2}x^2]|\leq 10^{-2}[c_1t+\frac{c_2}{2}x^2],$$
if $|x|<c$, $t\in [0,\tau]$.

In particular,
\begin{equation}\label{sieteocho}
c_1(1-10^{-2})t-\frac{c_2}{2}(1+10^{-2})x^2\leq \sigma_1^0(x,t)\leq c_1(1+10^{-2})t-\frac{c_2}{2}(1-10^{-2})x^2\end{equation}
for real $x$, if $|x|<\overline{c}$, $t\in [0,\tau]$. Similarly, we have that
\begin{equation}\label{sieteocho1/2}
c'_1(1+10^{-2})t-\frac{c_2}{2}(1+10^{-2})x^2\leq \sigma_1^0(x,t)\leq c'_1(1-10^{-2})t-\frac{c_2}{2}(1-10^{-2})x^2\end{equation}
for real $x$, if $|x|<\overline{c}$, $t\in [-\tau,0]$.

On the other hand, for the Muskat solution in \cite{CCFGL}, \eqref{rt} gives
\begin{equation}\label{sietenueve}
\sigma_1^0(x,0)<-c'\quad \text{if $\frac{\overline{c}}{100}\leq|x|\leq \pi$, $x$ real}.
\end{equation}
Then, by taking $\tau$ small enough, we achieve
\begin{equation}\label{sietediez}
\sigma_1^0(x,t)<-\frac{c'}{2}\quad \text{if $\frac{\overline{c}}{100}\leq |x|\leq\pi$, $x$ real, $t\in [0,\tau]$}.\end{equation}
If $h(x,t)$ is a positive functions defined for $x\in\T$, $0\leq t\le
 \tau$, and if
\begin{equation}\label{sieteonce}
0<h(x,t)<\frac{1}{2}c_b\quad \text{for all such $(x,t)$ (see property 1 above)},\end{equation}
then
\begin{equation}\label{sietedoce}
|\sigma_1^0(x\pm i h(x,t),t)-\sigma_1^0(x,t)|\leq C h(x,t).
\end{equation}

We will work with a perturbed Rayleigh-Taylor function $\sigma_1(x,t)$, defined and analytic on $\{|\Im x|<h(\Re x,t)\}$ for each fixed $t\in[0,\tau]$; and smooth on the closure of that region.

We assume that
\begin{equation}\label{sietecatorce}
\left|\partial_x^j(\sigma_1(x,t)-\sigma_1^0(x,t))\right|\leq
\tau^{10}\end{equation} for $0\leq j\leq 2$, and for $|\Im
x|\leq h(\Re x,t)$, $0\leq t\leq \tau$.

We have not yet picked the function $h(x,t)$. This we do in the
next section.

\subsection{Two Height Functions}\label{seccion8}

In this section, we suppose we are given a Rayleigh-Taylor
function $\sigma_1^0$ as in the section (\ref{seccion7}). We pick
$A>>1$  and $0<\tau$, $\kappa$, $\lambda<<1$ real numbers. We first pick $A$ large
enough, then pick $\tau$ small enough (depending on $A$). Having picked $A$ and $\tau$, we then take $0<\lambda<\tau^{100}$ and $0<\kappa<\tau^{100}$.  We fix such $A$, $\tau$, $\lambda$ until the end of the paper.
The parameter $\kappa$ will be held fixed until section \ref{ultima} below, in which will let $\kappa$ tend to zero. From now on, we suppose $A$, $\tau$, $\lambda$ and $\kappa$ have been picked as above.

\begin{figure}
\centering
\includegraphics[width=.48\textwidth]{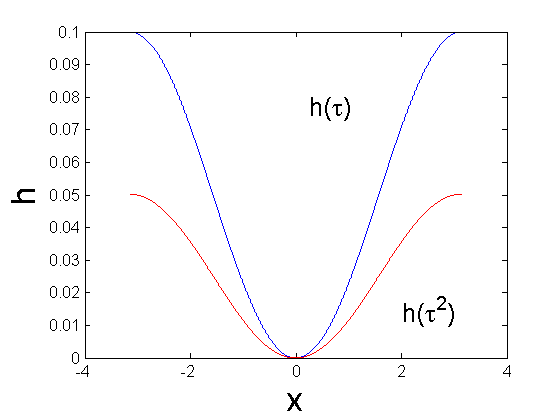}\caption{Red: Function $h(x,t)$ at time $t=\tau^2$, with $\kappa=0$. Blue: Function $h(x,t)$ at time $t=\tau$, with $\kappa=0$.}
\label{htau1}
\end{figure}

\begin{figure}
\centering
\includegraphics[width=.48\textwidth]{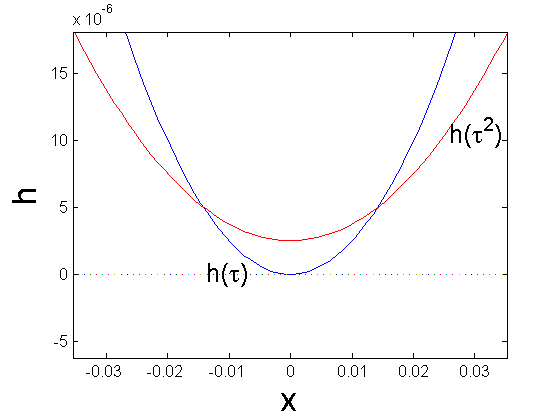}\caption{Zoom of Figure 1 around the point $(0,0)$.}
\label{htau2}
\end{figure}

\begin{figure}
\centering
\includegraphics[width=.5\textwidth]{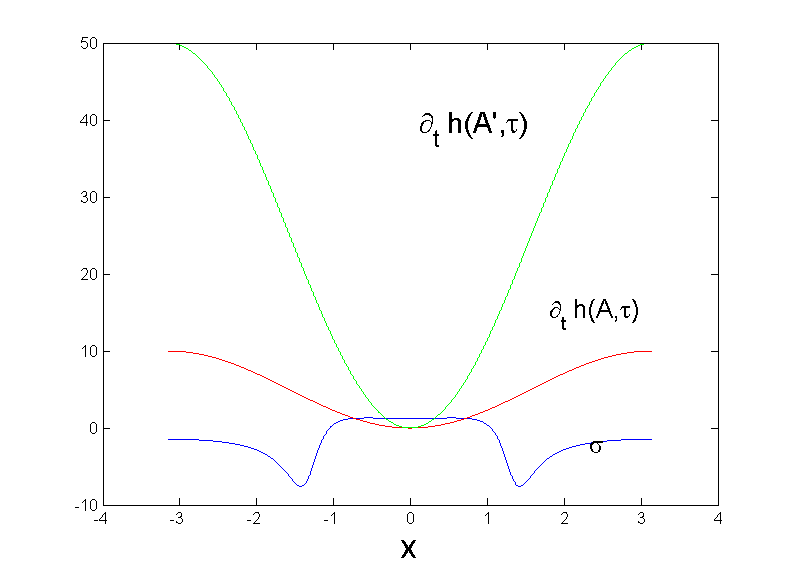}\caption{Blue: R-T function at time $t=\tau$. Red: Function $\pa_t h(x,t)$ at time $t=\tau$. Green: Function $\pa_t h(x,t)$ at time $t=\tau$ with a larger value of A.}
\label{sigma1}
\end{figure}

\begin{figure}
\centering
\includegraphics[width=.5\textwidth]{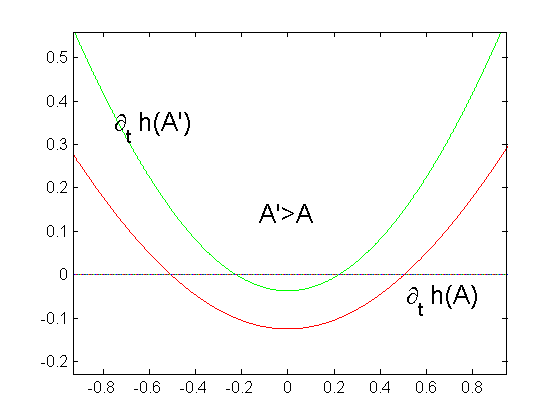}
\caption{Zoom of  Figure 3 around the point $(0,0)$.}
\label{sigma2}
\end{figure}

For $x\in\T$ and $\tau^2\leq t\leq \tau$, we define
\begin{equation}\label{ochouno}
h(x,t)=A^{-1}(\tau^2-t^2)+(A^{-1}-A(\tau-t))\sin^2\left(\frac{x}{2}\right)+\kappa.
\end{equation}
See Figures (\ref{htau1}) and (\ref{htau2}). Thus,
\begin{equation}\label{ochodos}
\partial_t h(x,t)=-2A^{-1}t+A\sin^2\left(\frac{x}{2}\right).
\end{equation}
See Figures (\ref{sigma1}) and (\ref{sigma2}).
If $||x||\geq 10 A^{-1}t^{\frac{1}{2}}$, then $\sin^2(x/2)\geq
10^{-2}||x||^2\geq A^{-2}t$, hence (\ref{ochodos}) yields
\begin{equation}\label{ochotres}
\left|\partial_t h(x,t)\right|\leq
2A^{-1}t+A\sin^2\left(\frac{x}{2}\right)\quad \text{for all $x\in
\T$, $t\in[\tau^2,\tau]$,}\end{equation} and
\begin{equation}\label{ochocuatro}
\left|\partial_t h(x,t)\right|\leq
3A\sin^2\left(\frac{x}{2}\right)\quad \text{if $||x||\geq
10A^{-1}t^{\frac{1}{2}}$, $t\in[\tau^2,\tau]$}.\end{equation} On
the other hand, (\ref{ochouno}) yields $h(x,t)\geq
\frac{1}{2}A^{-1}\sin^2(x/2)$ for all $x\in\T$, $t\in [\tau^2,\tau)$. Hence,
(\ref{ochocuatro}) implies
\begin{equation}\label{ochocinco}
\left|\partial_t h(x,t)\right|\leq6A^2h(x,t)\quad \text{for
$||x||\geq 10A^{-1}t^{\frac{1}{2}}$,
$t\in[\tau^2,\tau]$}.\end{equation} Next, let $||x||\leq
\frac{\overline{c}}{2}$ with $\overline{c}$ as in inequality
(\ref{sieteocho}). Applying that estimate we learn that
$$\sigma_1^0(x,t)+\pa_t h(x,t)-A^{\frac{1}{2}}h(x,t)$$
$$\geq
\left(\tilde{c}_1t-\tilde{C}_1\sin^2\left(\frac{x}{2}\right)\right)+\left(-2A^{-1}t+A\sin^2\left(\frac{x}{2}\right)\right)
-A^{\frac{1}{2}}\left[A^{-1}(\tau^2-t^2)+(A^{-1}-A(\tau-t))\sin^2\left(\frac{x}{2}\right)+\kappa\right]$$
$$=\left[\tilde{c}_1t-2A^{-1}t-A^{-\frac{1}{2}}(\tau^2-t^2)-A^\frac{1}{2}\kappa\right]+\left[-\tilde{C}_1+A-(A^{-\frac{1}{2}}
-A^{\frac{3}{2}}(\tau-t))\right]\sin^2\left(\frac{x}{2}\right).$$
The second term in square brackets is positive (since $A$ is large
but $A^{3/2}\tau$ is small), and the first term in
square brackets is greater than $c\tau^2$.

Hence,
$$\sigma_1^0(x,t)+\partial_th(x,t)-A^{\frac{1}{2}}h(x,t)\geq
c\tau^2\quad\text{if $||x||\leq\frac{1}{2}\overline{c}$,
$t\in[\tau^2,\tau]$}.$$
\begin{rem}The reason we restrict to $t\geq \tau^2$ is that, if we
allowed $t<<\tau^2$, then the first term in square brackets above
might be negative. \end{rem}
 On the other hand, suppose that
$||x||>\frac{1}{2}\overline{c}$, $t\in[\tau^2,\tau]$. Then
$|\sigma_1^0(x,t)|\leq C$, and (\ref{ochodos}) yields $\pa_t
h(x,t)\geq c A$, while (\ref{ochouno}) implies that $0\leq
h(x,t)\leq C'A^{-1}$. Hence, in this case,
$$\sigma_1^0(x,t)+\pa_t h(x,t)-A^{\frac{1}{2}}h(x,t)\geq
c\tau^2.$$ Thus, in all cases, we have
\begin{equation}\label{ochoseis}
\sigma_1^0(x,t)+\pa_th(x,t)-A^{\frac{1}{2}}h(x,t)\geq
c\tau^2\quad\text{for $x\in\T$, $t\in[\tau^2,\tau]$}\end{equation}
Next, suppose $||x||\leq 20A^{-1}t^{\frac{1}{2}}$,
$t\in[\tau^2,\tau]$. Then inequality (\ref{sieteocho}) gives
$\sigma_1^0(x,t)\geq ct$, while (\ref{ochouno}) gives $0\leq
h(x,t)\leq CA^{-1}t$ and \eqref{ochodos} gives $|\pa_t h(x,t)|\leq CA^{-1}t$. Hence,
\begin{equation}\label{ochosiete}
\sigma_1^0(x,t)-\left|\partial_t
h(x,t)\right|-A^\frac{1}{2}h(x,t)\geq c't\geq c'\tau^2\quad
\text{for $||x||\leq 20A^{-1}t^{\frac{1}{2}}$,
$t\in[\tau^2,\tau]$}.\end{equation} Next, we investigate the
$x-$derivatives of $h(x,t)$, $h_t(x,t)$ and $\sigma_1^0(x,t)$.
From (\ref{ochouno}) and (\ref{ochodos}), we see that
\begin{equation}\label{ochoocho}
\left|\partial_x^j h(x,t)\right|\leq CA^{-1}\quad \text{for $x\in
\T$, $t\in[\tau^2,\tau]$, $0\leq j\leq 3$};\end{equation} and
\begin{equation}\label{ochonueve}
\left|\pa_x^j h_t(x,t)\right|\leq CA\quad \text{for $x\in \T$,
$t\in[\tau^2,\tau]$, $0\leq j\leq 3$}.\end{equation} Also,
\begin{equation}\label{ochodiez}
\left|\partial_x^j\sigma_1^0(x,t)\right|\leq C\quad \text{for
$x\in \T$, $t\in[\tau^2,\tau]$, $0\leq j\leq 2$}.\end{equation}
(See section (\ref{seccion7})).

Now suppose that $||x||\leq 20A^{-1}t^\frac{1}{2}$. According to
estimate (\ref{sieteocho}), we have
\begin{equation}\label{ochoonce}
ct\leq \sigma_1^0(x,t)\leq Ct\quad \text{for $||x||\leq
20A^{-1}t^{\frac{1}{2}}$, $t\in[\tau^2,\tau]$}.
\end{equation}
Since (\ref{ochoonce}) holds for all
$x\in[-20A^{-1}t^\frac{1}{2},20A^{-1}t^\frac{1}{2}]$, the
mean-value theorem implies that $|\pa_x\sigma_1^0(x,t)|\leq
CAt^\frac{1}{2}$ somewhere in
$[-20A^{-1}t^\frac{1}{2},20A^{-1}t^\frac{1}{2}]$. Hence, by
(\ref{ochodiez}), we have
\begin{equation}\label{ochodoce}
\left|\pa_x\sigma_1^0(x,t)\right|\leq
CAt^\frac{1}{2},\quad\text{for all $||x||\leq
20A^{-1}t^\frac{1}{2}$, $t\in[\tau^2,\tau]$}.\end{equation} From
(\ref{ochodiez}), (\ref{ochoonce}), (\ref{ochodoce}), we see that
\begin{equation}\label{ochotrece}
\left|\pa_x^j\sigma_1^0(x,t)\right|\leq
CA\left(t^{\frac{1}{2}}\right)^{2-j}\quad\text{for all $||x||\leq
20A^{-1}t^\frac{1}{2}$, $t\in[\tau^2,\tau]$, $0\leq j\leq
2$}.\end{equation} Next, return to $h$. From
(\ref{ochouno}), we have
\begin{equation}\label{ochocatorce}
0\leq h(x,t)\leq C A^{-1}t\quad\text{for all $||x||\leq
20A^{-1}t^\frac{1}{2}$, $t\in[\tau^2,\tau]$}\end{equation} and
\begin{equation}\label{ochoquince}
\left|\pa_xh(x,t)\right|\leq C A^{-2}t^\frac{1}{2}\quad\text{for all
$||x||\leq 20A^{-1}t^\frac{1}{2}$,
$t\in[\tau^2,\tau]$}.\end{equation} From
(\ref{ochoocho}),(\ref{ochocatorce}),(\ref{ochoquince}), we see
that
\begin{equation}\label{ochodieciseis}
\left|\pa_x^jh(x,t)\right|\leq
C\left(t^\frac{1}{2}\right)^{2-j}\text{for all $||x||\leq
20A^{-1}t^\frac{1}{2}$, $t\in[\tau^2,\tau]$, $0\leq j\leq
2$}.\end{equation} Similarly, regarding $\partial_t h(x,t)$, we
learn from (\ref{ochodos}) that
\begin{equation}\label{ochodiecisiete}
\left|\pa_th(x,t)\right|\leq CA^{-1}t\quad\text{for all $||x||\leq
20A^{-1}t^\frac{1}{2}$, $t\in[\tau^2,\tau]$};\end{equation} and
\begin{equation}\label{ochodieciocho}
\left|\pa_x\pa_th(x,t)\right|\leq Ct^\frac{1}{2}\quad\text{for all $||x||\leq
20A^{-1}t^\frac{1}{2}$, $t\in[\tau^2,\tau]$}.\end{equation} From
(\ref{ochonueve}), (\ref{ochodiecisiete}) and
(\ref{ochodieciocho}) we obtain the estimate
\begin{equation}\label{ochodiecinueve}
\left|\pa_x^j\pa_th(x,t)\right|\leq
CA\left(t^{\frac{1}{2}}\right)^{2-j}\quad\text{for all $||x||\leq
20A^{-1}t^\frac{1}{2}$, $t\in[\tau^2,\tau]$, $0\leq j\leq
2$}.\end{equation}

It will also be useful to estimate $\pa_x^j\sigma_1^0(z,t)$ for
$z=x+iy$ with $||x||\leq 20A^{-1}t^\frac{1}{2}$ and $|y|\leq
h(x)$.

From the assumption 1 in section (\ref{seccion7}) and since
$0\leq h\leq CA^{-1}t$ we have $|\pa_y^{j+1}\sigma_1^0(x+iy,t)|\leq C$ in
this range, hence
$$\left|\pa_x^j\sigma_1^0(x+iy,t)\right|\leq
\left|\pa_x^j\sigma_1^0(x,t)\right|+C|y|\leq
CA\left(t^\frac{1}{2}\right)^{2-j}+CA^{-1}t\leq
C'A\left(t^\frac{1}{2}\right)^{2-j}.$$ Thus,
$$
\left|\pa_x^j\sigma_1^0(x+iy,t)\right|\leq
CA\left(t^\frac{1}{2}\right)^{2-j}$$\begin{equation}\label{ochoveinte}\text{for
$z=x+iy$, $||x||\leq 20A^{-1}t^\frac{1}{2}$, $t\in[\tau^2,\tau]$,
$|y|\leq h(x)$, $0\leq j\leq 2$}.\end{equation} We look at
$h(x,\tau^2)$. From (\ref{ochouno}), we have
\begin{equation}\label{ochoveintiuno}
h(x,\tau^2)\geq
\frac{1}{2}A^{-1}\tau^2+\frac{1}{2}A^{-1}\sin^2\left(\frac{x}{2}\right)\quad\text{for
$x\in\T$}.\end{equation} Let
\begin{equation}\label{ochoveintidos}
\hbar(x,t)=\frac{1}{4}\left(A^{-1}\tau^2+A^{-1}\sin^2\left(\frac{x}{2}\right)\right)+A^{-2}\tau
t+At\sin^2\left(\frac{x}{2}\right)\quad\text{for $x\in\T$,
$t\in[-\tau^2,\tau^2].$}\end{equation} Note that
\begin{equation}\label{ochoveintitres}
0\leq \hbar(x,\tau^2)\leq h(x,\tau^2)\quad \text{for all
$x\in\T$}.
\end{equation}
From (\ref{ochoveintidos}), we have
\begin{equation}\label{ochoveinticuatro}
\pa_t\hbar(x,t)=A^{-2}\tau+A\sin^2\left(\frac{x}{2}\right)\quad\text{for
$x\in\T$, $t\in[-\tau^2,\tau^2]$.}\end{equation} On the other hand,
(\ref{ochoveintidos}) also gives
\begin{equation}\label{ochoveinticinco}
\hbar(x,t)\geq
\frac{1}{8}\left(A^{-1}\tau^2+A^{-1}\sin^2\left(\frac{x}{2}\right)\right)\quad\text{for
$x\in\T$, $t\in[-\tau^2,\tau^2]$.}\end{equation} From
(\ref{ochoveinticuatro}), (\ref{ochoveinticinco}) we obtain the
estimate
\begin{equation}\label{ochoveintiseis}
\left|\pa_t \hbar(x,t)\right|\leq
C\tau^{-1}\hbar(x,t)\quad\text{for all $x\in\T$,
$t\in[-\tau^2,\tau^2]$.}\end{equation} Next, suppose $||x|| < \overline{c}$; from (\ref{sieteocho}) and \eqref{sieteocho1/2} we see that
\begin{align}
\sigma_1^0(x,t)\geq & ct-C\sin^2\left(\frac{x}{2}\right)\quad \text{for $t\in[0,\,\tau^2]$},\\
\sigma_1^0(x,t)\geq & c't-C'\sin^2\left(\frac{x}{2}\right)\quad \text{for $t\in[-\tau^2,\,0]$}
\end{align}

Therefore
$$\sigma_1^0(x,t)+\pa_t\hbar(x,t)-A^\frac{1}{2}\hbar(x,t)\geq \left(ct-C\sin^2\left(\frac{x}{2}\right)\right)
+\left(A^{-2}\tau+A\sin^2\left(\frac{x}{2}\right)\right)$$
$$-A^\frac{1}{2}\left[\frac{1}{4}\left(A^{-1}\tau^2+A^{-1}\sin^2\left(\frac{x}{2}\right)\right)+A^{-2}\tau
t+At\sin^2\left(\frac{x}{2}\right)\right]$$
$$=\left\{ct+A^{-2}\tau-\frac{1}{4}A^{-\frac{1}{2}}\tau^2-A^{-\frac{3}{2}}\tau
t\right\}+\left\{-C+A-\frac{1}{4}A^{-\frac{1}{2}}-A^{\frac{3}{2}}t\right\}\sin^2\left(\frac{x}{2}\right),$$
for $t\in [0,\tau^2]$. The first expression in curly brackets is greater than
$\frac{1}{2}A^{-2}\tau$, and the second expression in curly
brackets is positive.

Therefore,
$$\sigma_1^0(x,t)+\pa_t\hbar(x,t)-A^\frac{1}{2}\hbar(x,t)\geq
\frac{1}{2}A^{-2}\tau\quad \text{if $||x||< \overline{c}$,
$t\in[0,\tau^2]$}.$$ On the other hand, suppose $||x||\geq
\overline{c}$. Then (\ref{ochoveinticuatro}) gives
$\pa_t\hbar(x,t)\geq cA$, whereas (\ref{ochoveintidos}) gives
$0\leq \hbar(x,t)\leq CA^{-1}$, and we know that
$|\sigma_1^0(x,t)|\leq C$. Hence,
$$\sigma_1^0(x,t)+\pa_t\hbar(x,t)-A^\frac{1}{2}\hbar(x,t)\geq
cA-C-CA^{-\frac{1}{2}}>c'A>\frac{1}{2}A^{-2}\tau\quad\text{if $||x||\geq
\overline{c}$, $t\in [0,\tau^2]$.}$$ Thus, in all cases,
\begin{equation}\label{ochoveintisietea}
\sigma_1^0(x,t)+\pa_t\hbar(x,t)-A^\frac{1}{2}\hbar(x,t)\geq
\frac{1}{2}A^{-2}\tau\quad \text{for $x\in\T$,
$t\in[0,\tau^2]$}.\end{equation}
Similar computations yield
\begin{equation}\label{ochoveintisiete}
\sigma_1^0(x,t)+\pa_t\hbar(x,t)-A^\frac{1}{2}\hbar(x,t)\geq
\frac{1}{2}A^{-2}\tau\quad \text{for $x\in\T$,
$t\in[-\tau^2,\tau^2]$}.\end{equation}

Regarding the smoothness of
$\sigma_1^0$, $\hbar$, $\pa_t \hbar$, we see from
(\ref{ochoveintidos}) and (\ref{ochoveinticuatro}) that
\begin{equation}\label{ochoveintiocho}
\left|\pa_x^j \hbar (x,t)\right|\leq CA^{-1}\quad \text{for
$x\in\T$, $t\in[-\tau^2,\tau^2]$, $0\leq j\leq 3$};\end{equation} and
\begin{equation}\label{ochoveintinueve}
\left|\pa_x^j \pa_t\hbar (x,t)\right|\leq CA\quad\text{for
$x\in\T$, $t\in[-\tau^2,\tau^2]$, $0\leq j\leq 3$}.\end{equation}
Also, from our assumption on $\sigma_1^0$ in section
(\ref{seccion7}), we have
\begin{equation}\label{ochotreinta}
\left|\pa_z^j\sigma_1^0(z,t)\right|\leq C\quad \text{for $\Re z
\in \T$, $|\Im z|\leq c$, $t\in[-\tau^2,\tau]$, $0\leq j\leq
2$.}\end{equation} Note that (\ref{ochotreinta}) holds whenever
$t\in[-\tau^2,\tau^2]$ and $|\Im z|\leq \hbar(\Re z,t)$, since
(\ref{ochoveintidos}) gives $0\leq \hbar(x,t)\leq CA^{-1}$ for all
$x\in\T$, $t\in[-\tau^2,\tau^2]$. Similarly, (\ref{ochotreinta}) holds
whenever $t\in[\tau^2,\tau]$ and $|\Im z|\leq h(\Re z,t)$, since
(\ref{ochouno}) then gives $0\leq h(\Re z,t)\leq CA^{-1}$.

\subsection{Rayleigh Taylor in the Complex Domain}\label{seccion9}
Let $\sigma_1(x,t)$, $\sigma_1^0$, $h(x,t)$ $\hbar(x,t)$ be as in
the sections (\ref{seccion7}) and (\ref{seccion8}).

In the time interval $[\tau^2,\tau]$ we define the functions
\begin{eqnarray}\label{nueveuno}
\sigma_+^0(x,t)&=&\sigma_1^0(x+ih(x,t),t)\quad\text{for $x\in \T$,
$t\in [\tau^2,\tau]$}\\ \label{nuevedos}
\sigma_-^0(x,t)&=&\sigma_1^0(x-ih(x,t),t)\quad\text{for $x\in \T$,
$t\in [\tau^2,\tau]$}\\ \label{nuevetres}
\sigma_+(x,t)&=&\sigma_1(x+ih(x,t),t)\quad\text{for $x\in \T$,
$t\in [\tau^2,\tau]$}\\ \label{nuevecuatro}
\sigma_-(x,t)&=&\sigma_1(x-ih(x,t),t)\quad\text{for $x\in \T$,
$t\in [\tau^2,\tau]$}.\end{eqnarray} When $t\in[-\tau^2,\tau^2]$ we use
instead the following definitions:
\begin{eqnarray}\label{nuevecinco}
\sigma_+^0(x,t)&=&\sigma_1^0(x+i\hbar(x,t),t)\quad\text{for $x\in
\T$, $t\in [-\tau^2,\tau^2]$}\\ \label{nueveseis}
\sigma_-^0(x,t)&=&\sigma_1^0(x-i\hbar(x,t),t)\quad\text{for $x\in
\T$, $t\in [-\tau^2,\tau^2]$}\\ \label{nuevesiete}
\sigma_+(x,t)&=&\sigma_1(x+i\hbar(x,t),t)\quad\text{for $x\in \T$,
$t\in [-\tau^2,\tau^2]$}\\ \label{nueveocho}
\sigma_-(x,t)&=&\sigma_1(x-i\hbar(x,t),t)\quad\text{for $x\in \T$,
$t\in [-\tau^2,\tau^2]$}.\end{eqnarray} Note that at $t=\tau^2$, we are
giving two conflicting definitions of the Rayleigh Taylor
functions. That is a defect in the notation, but it will lead to
no trouble.

From (\ref{ochouno}) and (\ref{ochoveintidos}), we see that $0\leq
h \leq C A^{-1}$ and $0\leq \hbar\leq CA^{-1}$ whenever $h$,
$\hbar$ are defined. Hence, property 1 in section
(\ref{seccion7}), tells us that
\begin{equation}\label{nuevenueve}
\left|\sigma_{\pm}^0(x,t)-\sigma_1^0(x,t)\right|\leq C
h(x,t)\quad\text{for $x\in\T$,
$t\in[\tau^2,\tau]$;}\end{equation}and
\begin{equation}\label{nuevediez}
\left|\sigma_{\pm}^0(x,t)-\sigma_1^0(x,t)\right|\leq C
\hbar(x,t)\quad\text{for $x\in\T$,
$t\in[-\tau^2,\tau^2]$;}\end{equation} moreover, (\ref{sietecatorce})
(modified: use h when $t\in[\tau^2,\tau]$, use $\hbar$ when
$t\in[-\tau^2,\tau^2]$) together with (\ref{ochoocho}) and
(\ref{ochoveintiocho}), shows that
\begin{equation}\label{nueveonce}
\left|\pa_x^j\left(\sigma_{\pm}-\sigma^0_{\pm}\right)(x,t)\right|\leq
C\tau^{10}\quad\text{for $x\in\T$, $t\in[-\tau^2,\tau]$, $0\leq j\leq
2$.}\end{equation} Again using (\ref{ochoocho}) and
(\ref{ochoveintiocho}), and recalling the property 1 in the
section (\ref{seccion7}), we see that
\begin{equation}\label{nuevedoce}
\left|\pa_x^j\sigma_{\pm}^0(x,t)\right|\leq C \quad\text{for
$x\in\T$, $t\in[-\tau^2,\tau]$, $0\leq j\leq 2$.}\end{equation} From
(\ref{nueveonce}) and (\ref{nuevedoce}) we obtain
\begin{equation}\label{nuevetrece}
\left|\pa_x^j\sigma_{\pm}(x,t)\right|\leq C \quad\text{for
$x\in\T$, $t\in[-\tau^2,\tau]$, $0\leq j\leq 2$.}\end{equation} Next,
we study $\sigma_\pm(x,t)$ for $||x||\leq
20A^{-1}t^{\frac{1}{2}}$, $t\in[\tau^2,\tau]$. For such $(x,t)$ we
argue as follows.

From (\ref{nueveuno}), together with the assumption 1 in section
(\ref{seccion7}), we see that
$$\left|\sigma_+^0(x,t)-\sigma_1^0\right(x,t)|=\left|\sigma_1^0(x+ih(x,t),t)-\sigma_1^0(x,t)\right|\leq
Ch(x,t);$$ and
$$\left|\pa_z\sigma_1^0(z,t)|_{z=x+ih(x,t)}-\pa_x\sigma_1^0(x,t)\right|\leq
Ch(x,t).$$ Hence
\begin{equation}\label{nuevecatorce}
\left|\sigma_+^0(x,t)\right|\leq
\left|\sigma_1^0(x,t)\right|+Ch(x,t)\leq Ct,\end{equation} thanks to (\ref{ochoonce}) and (\ref{ochocatorce}); and also
\begin{equation}\label{nuevequince}
\left|\pa_z\sigma_1^0(z,t)|_{z=x+ih(x,t)}\right|\leq\left|\pa_x\sigma_1^0(x,t)\right|+Ch(x,t)\leq
CAt^\frac{1}{2},\end{equation} thank to (\ref{ochodoce}) and
(\ref{ochocatorce}).

Differentiating (\ref{nueveuno}), we find that
\begin{equation}\label{nuevedieciseis}
\pa_x\sigma_+^0(x,t)=\left[\pa_z
\sigma_1^0(z,t)|_{z=x+ih(x,t)}\right](1+i\pa_xh(x,t)).\end{equation}
Since $|\pa_x h(x,t)|\leq A^{-2}t^\frac{1}{2}<1$ by
(\ref{ochoquince}), it follows from (\ref{nuevequince}) and
(\ref{nuevedieciseis}) that
\begin{equation}\label{nuevediecisiete}
\left|\pa_x\sigma^+_0(x,t)\right|\leq CAt^\frac{1}{2}.
\end{equation}
From (\ref{nuevedoce}), (\ref{nuevecatorce}),
(\ref{nuevediecisiete}), we obtain the estimate
\begin{equation}\label{nuevedieciocho}
\left|\pa_x^j\sigma_+^0(x,t)\right|\leq
CA\left(t^\frac{1}{2}\right)^{2-j}\quad\text{for $||x||\leq
20A^{-1}t^{\frac{1}{2}}$, $t\in[\tau^2,\tau]$, $0\leq j\leq
2$.}\end{equation} Similarly
\begin{equation}\label{nuevediecinueve}
\left|\pa_x^j\sigma_-^0(x,t)\right|\leq
CA\left(t^\frac{1}{2}\right)^{2-j}\quad\text{for $||x||\leq
20A^{-1}t^{\frac{1}{2}}$, $t\in[\tau^2,\tau]$, $0\leq j\leq
2$.}\end{equation} From (\ref{nueveonce}),
(\ref{nuevedieciocho}), (\ref{nuevediecinueve}), we obtain the
estimate
\begin{equation}
\label{nueveveinte} \left|\pa_x^j\sigma_\pm(x,t)\right|\leq
CA\left(t^\frac{1}{2}\right)^{2-j}\quad\text{for $||x||\leq
20A^{-1}t^{\frac{1}{2}}$, $t\in[\tau^2,\tau]$, $0\leq j\leq
2$.}\end{equation} The main results of this sections are
(\ref{nuevetrece}), (\ref{nueveveinte}) and the following
consequences of (\ref{nuevenueve}), (\ref{nuevediez}),
(\ref{nueveonce}):
\begin{equation}\label{nueveveintiuno}
\Re \sigma_{\pm}(x,t)\geq
\sigma_1^0(x,t)-Ch(x,t)-C\tau^{10}\quad\text{for $x\in\T$,
$t\in[\tau^2,\tau]$;}\end{equation} and
\begin{equation}\label{nueveveintidos}
\Re \sigma_{\pm}(x,t)\geq
\sigma_1^0(x,t)-C\hbar(x,t)-C\tau^{10}\quad\text{for $x\in\T$,
$t\in[-\tau^2,\tau^2]$.}\end{equation} We also need to estimate $\Im
\sigma_\pm(x,t)$. From (\ref{nuevenueve}),(\ref{nuevediez}),
(\ref{nueveonce}) and the fact that $\sigma_1^0(x,t)$ is real for
$x\in\T$, we have
\begin{equation}\label{nueveveintitres}
\left|\Im \sigma_{\pm}(x,t)\right|\leq
Ch(x,t)+C\tau^{10}\quad\text{for $x\in\T$,
$t\in[\tau^2,\tau]$;}\end{equation} and
\begin{equation}\label{nueveveinticuatro}
\left|\Im \sigma_{\pm}(x,t)\right|\leq
C\hbar(x,t)+C\tau^{10}\quad\text{for $x\in\T$,
$t\in[-\tau^2,\tau^2]$}\end{equation}

\subsection{Sobolev Embedding in $\Omega(t)$}\label{seccionsobo}

Let $h(x)=h(x,t)$ for some fixed $t\in[\tau^2,\,\tau]$, or else let $h(x)=\hbar(x,t)$ for some fixed $t\in[0,\,\tau^2]$;
here, $h(x,t)$, $\hbar(x,t)$ are as in section \ref{seccion8}.

The goal of this section is to check that a simple Sobolev-type embedding theorem holds on the domain
\begin{equation*}
\Omega=\{\zeta\in\C\,:\,|\Im \zeta|<h(\Re\zeta)\}
\end{equation*}
with constants independent of the parameters $A$, $\tau$, $\kappa$ used to define $h(x,t)$ and $\hbar(x,t)$.

We write
\begin{equation*}
\Gamma_\pm=\{\zeta\in\C\,:\,\Im \zeta=\pm h(\Re\zeta)\}.
\end{equation*}
Note that $|h'(x)|\leq C$ for all $x\in\R$, with $C$ independent of $A$, $\tau$, $\kappa$ (See \eqref{ochouno}, \eqref{ochoveintidos}). Let $F(\zeta)$ be continuous on $\Omega^{\text{closure}}$ and holomorphic on $\Omega$. Assume that $F(\zeta)$ is $2\pi-$periodic and real for real $\zeta$. We will prove the following assertions:
\begin{description}
\item $X_1$. Suppose $F|_{\Gamma_+}$ belongs to $\text{Lip}(\alpha)$, with $0<\alpha<1$. Then $F$ belongs to Lip$(\alpha)$ as a function on $\Omega^\text{closure}$, and the $\text{Lip}(\alpha)$-norm of $F$ is at most $C$ times that of $F|_{\Gamma_+}$, where $C$ depends only on $\alpha$.
\item $X_2$. Suppose $F$ belongs to the Sobolev space $H^4(\Gamma_+)$. Then $F$ belongs to $C^{3.5}(\Omega)$ and  $$||F||_{C^{3.5}(\Omega)}\leq C||F||_{H^4(\Gamma_+)}$$
    for a universal constant $C$.
\end{description}
Proof of $X_1$: It is enough to show that $F\in \text{Lip}(\alpha)$ on $\Omega$, and to bound its $\text{Lip}(\alpha)$, and to bound its $\text{Lip}(\alpha)$ norm there, assuming that $F|_{\Gamma_+}$ has $\text{Lip}(\alpha)$ norm less or equal to one.

Let $z$, $z'\in \Omega$. By the Cauchy integral formulas we have
\begin{equation*}
F(z)-F(z')=\frac{1}{2\pi i}\int_{\Gamma_+}F(\zeta)\left[\frac{1}{\zeta-z}-\frac{1}{\zeta-z'}\right]d\zeta-\frac{1}{2\pi i}\int_{\Gamma_-}F(\zeta)\left[\frac{1}{\zeta-z}-\frac{1}{\zeta-z'}\right]d\zeta.
\end{equation*}
We show that both the integrals on the right-hand side are dominated by $$C(\alpha)|z-z'|^\alpha.$$ This will complete the proof of $X_1$. We examine only the first integral here; the second is completely analogous, and we have $||F|_{\Gamma_-}||_{\text{Lip}(\alpha)}=||F|_{\Gamma_+}||_{\text{Lip}(\alpha)}=1$ since $F(\zeta)$ is real for real $\zeta$.

Fix points $\zeta_0$, $\zeta_0'\in \Gamma_+$, with $\Re \zeta_0=\Re z$ and $\Re \zeta_0'=\Re z'$.

Note that $|\zeta-\zeta'|\leq C_1|\Re \zeta-\Re \zeta'|$ for $\zeta$, $\zeta'\in \Gamma_+$. In particular,
$|\zeta_0-\zeta_0'|\leq C_1|z-z'|$.

We subdivide the contour $\Gamma_+$ into

\begin{equation*}
\Gamma_+^\text{near}=\{\zeta\in \Gamma_+\,:\,|\Re \zeta-\Re z|\leq 100(C_1+1)|z-z'|\}
\end{equation*}
and
\begin{equation*}
\Gamma_+^\text{far}=\Gamma_+\setminus \Gamma_+^\text{near}.
\end{equation*}
Since
\begin{equation*}
\int_{\Gamma_+}\left[\frac{1}{\zeta-z}-\frac{1}{\zeta-z'}\right]d\zeta=0
\end{equation*}
by routine contour integration, we have
\begin{align}
&\int_{\Gamma_+} F(\zeta)\left[\frac{1}{\zeta-z}-\frac{1}{\zeta-z'}\right]d\zeta\nonumber\\
&=\int_{\Gamma_+}(F(\zeta)-F(\zeta_0))\left[\frac{1}{\zeta-z}-\frac{1}{\zeta-z'}\right]d\zeta\nonumber\\
&=\int_{\Gamma_+^\text{near}}(F(\zeta)-F(\zeta_0))\frac{d\zeta}{\zeta-z}-\int_{\Gamma_+^\text{near}}
(F(\zeta)-F(\zeta_0'))\frac{d\zeta}{\zeta-z'}\nonumber\\
&+\int_{\Gamma_+^\text{near}}
(F(\zeta_0)-F(\zeta_0'))\frac{d\zeta}{\zeta-z'}+\int_{\Gamma_+^\text{far}}(F(\zeta)-F(\zeta_0))
\left[\frac{1}{\zeta-z}-\frac{1}{\zeta-z'}\right]d\zeta\label{sobo1}.
\end{align}
Since
\begin{align}
|F(\zeta)-F(\zeta_0)|&\leq  C|\zeta-\zeta_0|^\alpha\nonumber\\
\leq & C' |\Re \zeta-\Re \zeta_0|^\alpha\nonumber\\
= & C'|\Re \zeta-\Re z|^\alpha\nonumber\\
\leq & C'|\zeta-z|^\alpha,\nonumber\end{align}
it follows that
\begin{equation}\label{sobo2}
\left|\int_{\Gamma_+^\text{near}}(F(\zeta)-F(\zeta_0))\frac{d\zeta}{\zeta-z}\right|\leq C''|z-z'|^\alpha.
\end{equation}
Similarly,
\begin{equation}\label{sobo3}
\left|\int_{\Gamma_+^\text{near}}(F(\zeta)-F(\zeta_0'))\frac{d\zeta}{\zeta-z'}\right|\leq C''|z-z'|^\alpha.
\end{equation}
Also,
\begin{equation*}
\left|\int_{\Gamma_+^\text{near}}\frac{d\zeta}{\zeta-z'}\right|\leq C,
\end{equation*}
since we may deform $\Gamma_+^\text{near}$ to a contour of length less or equal than $C|z-z'|$, and distance larger or equal than $|z-z'|$ from $z$.

Therefore,
\begin{align}
&\left|\int_{\Gamma_+^\text{near}}(F(\zeta_0)-F(\zeta_0'))\frac{d\zeta}{\zeta-z'}\right|\leq C|F(\zeta_0)-F(\zeta_0')|\nonumber\\
& \leq C'|\zeta_0-\zeta_0'|^\alpha\leq C''|\Re \zeta_0-\Re \zeta_0'|^\alpha=C''|\Re z-\Re z'|^\alpha\leq C'' | z- z'|^\alpha\label{sobo4}.
\end{align}
For $\zeta\in \Gamma_+^\text{far}$, we have

\begin{equation*}
|F(\zeta)-F(\zeta_0)|\cdot \left|\frac{1}{\zeta-z}-\frac{1}{\zeta-z'}\right|\leq C
|\zeta-\zeta_0|^\alpha \frac{|z-z'|}{|\zeta-z|^2}.
\end{equation*}
Since also $|\zeta-\zeta_0|\leq C|\Re \zeta-\Re \zeta_0|=C|\Re \zeta -\Re z|\leq C|\zeta-z|$, it follows that
\begin{equation}
|F(\zeta)-F(\zeta_0)|\cdot \left|\frac{1}{\zeta-z}-\frac{1}{\zeta-z'}\right|\leq C |z-z'||\zeta-z|^{\alpha-2}.
\end{equation}
Consequently,
\begin{equation}\label{sobo5}
\left|\int_{\Gamma_+^\text{far}}(F(\zeta)-F(\zeta_0))\left[\frac{1}{\zeta-z}-\frac{1}{\zeta-z'}\right]d\zeta\right|\leq C|z-z'|^\alpha.
\end{equation}
Putting \eqref{sobo2},...,\eqref{sobo5} into \eqref{sobo1}, we learn that
\begin{equation*}
\left|\int_{\Gamma_+} F(\zeta)\left[\frac{1}{\zeta-z}-\frac{1}{\zeta-z'}\right]d\zeta\right|\leq C|z-z'|^\alpha.
\end{equation*}
This completes the proof of $X_1$.

Proof of $X_2$: Suppose $F|_{\Gamma_+}\in H^4$. Then $F$ and its derivatives up to order 3 are bounded on $\Gamma_+$, hence also on $\Gamma_-$ (since $F(\overline{\zeta})=\overline{F(\zeta)}$). Therefore, by the maximum principle, $F$ and its derivatives up to $3^{\underline{\text{rd}}}$ order are bounded on $\Omega$. Also, since $F|_{\Gamma_+}\in H^4$, we know that $F'''(\zeta)$ belongs to $\text{Lip}(\frac{1}{2})$ on $\Gamma_+$. Applying $X_1$, we conclude that $F'''$ belongs to Lip$(\frac{1}{2})$ on $\Omega$. Thus $F\in C^{3.5}(\Omega)$, completing the proof of $X_2$.

When we apply $X_2$ in later sections, we will typically take $F(\zeta)=z_\mu(\zeta,t)-\un_\mu(\zeta,t)$, where $(z_\mu(\zeta,t))_{\mu=1,\,2}$ and $(\un_\mu(\zeta,t))_{\mu=1,\,2}$ are two solutions of the Muskat equation.


\section{Perturbing the Muskat solution}\label{seccionM}

Here we expose the necessary a priori estimates of the difference of two solutions of the Muskat equation in order to obtain the theorem \ref{mean} of subsection \ref{subseccion1}.

\subsection{Differentiating The Muskat Equation}\label{seccion12}
In this section, we suppose we are given a pair of functions
$z_\mu(x,t)$ ($\mu=1,2$),  $z_1(x,t)-x$ and $z_2(x,t)$ periodic of
period $2\pi$ in $x$, analytic in the complex variable $x$ in the
region \begin{equation*}
\Omega(t)=\{|\Im x|<h(\Re x,t)\},\end{equation*} and smooth in
$(x,t)$ on the closure of this region. We write
$$\Gamma_\pm(t)=\{x\pm ih(x,t)\,:\, x\in \T\}$$ for the two contours
bounding the above region.

Later, we will take $h(x,t)$ (of this section) to be either our
previous $h(x,t)$ or our previous $\hbar(x,t)$, depending on
whether we work with $t\in[\tau^2,\tau]$ or with $t\in[-\tau^2,\tau^2]$.
For now, we simply suppose that $h(x,t)$ is a positive, smooth
function of $x$ for each fixed $t$.

We suppose that our $z_\mu$ satisfy the Muskat equation
\begin{equation}\label{doceuno}
\pa_tz_\mu(x,t)=\int_{y\in\T}\frac{\sin\left(z_1(x,t)-z_1(x+y)\right)
\left(\pa_xz_\mu(x,t)-\pa_xz_\mu(x+y,t)\right)}
{\cosh\left(z_2(x,t)-z_2(x+y)\right)-\cos\left(z_1(x,t)-z_1(x+y)\right)}dy.
\end{equation}
Note the sign in (\ref{doceuno}) and note that there is no factor
$\frac{1}{2\pi}$ in front of the integral.

In (\ref{doceuno}), $x$ and $y \in \T$ are real. Later we
will use the analyticity of $z_\mu(x,t)$ in $x$ to move into the
complex plane.

We will assume the following COMPLEX CHORD-ARC CONDITION:

Let $x$, $y\in \Omega(t)$. Then
\begin{equation}\label{docedos}
|\cosh\left(z_2(x,t)-z_2(y,t)\right)-\cos\left(z_1(x,t)-z_1(y,t)\right)|\geq
c_{CA}\left[||\Re (x-y)||+|\Im (x-y)|\right]^2.
\end{equation}
Let $\underline{k}$ be a large enough positive integer constant (here, $\underline{k}=4$ is enough).

We apply the $\underline{k}-$th derivative $\pa_x^{\underline{k}}$
to both sides of (\ref{doceuno}). The goal of this section is to
see what results.

We obtain
$$
\pa_t \left[\pa_x^{\underline{k}}z_\mu(x,t)\right]$$$$
=\sum_{k+k'=\underline{k}}c(k,\underline{k})\int_{y\in\T}\pa_x^k\left\{\frac{\sin\left(z_1(x,t)-z_1(x+y)\right)}
{\cosh\left(z_2(x,t)-z_2(x+y)\right)-\cos\left(z_1(x,t)-z_1(x+y)\right)}
\right\}$$
\begin{equation}\label{docetres}
\times
\pa_x^{k'}\left\{\pa_xz_\mu(x,t)-\pa_xz_\mu(x+y,t)\right\}dy,
\end{equation}
where $c(k,\underline{k})$ are harmless coefficients and
$c(0,\underline{k})=1$.

We fix $t$, $y$ and write
$\hat{z}_\mu(x)=z_\mu(x,t)-z_\mu(x+y,t)$. The integrand in
(\ref{docetres}) is then (with $k+k'=\underline{k}$):
$$\left[\pa_x^k\left\{\frac{\sin(\hat{z}_1(x))}{\cosh(\hat{z}_2(x))-\cos(\hat{z}_1(x))}\right\}\right]
\left[\pa_x^{k'+1}\hat{z}_\mu(x)\right]$$
$$=\sum_{k^\sharp+k^b=k}\tilde{c}(k^\sharp,k^b)\left[\pa_x^{k^\sharp}\{\sin(\hat{z}_1(x))\}\right]
\left[\pa_x^{k^b}\left\{\frac{1}{\cosh(\hat{z}_2(x))-\cos(\hat{z}_1(x))}\right\}\right]
\left[\pa_x^{k'+1}\hat{z}_\mu(x)\right],$$ where
$\tilde{c}(k^\sharp,k^b)$ are harmless coefficients.

Now $\pa_x^{k^\sharp}\{\sin(\hat{z}_1(x))\}$ is a linear
combination (with harmless coefficients) of terms
$$\left[\sin^{(m^\sharp)}(\hat{z}_1(x))\right]\prod_{j=1}^{m^\sharp}\left(\pa_x^{k_j^\sharp}\hat{z}_1(x)\right),$$
where $k_j^\sharp\geq 1$ for each $j$,
$k_1^\sharp+\cdot\cdot\cdot+k_{m^\sharp}^\sharp=k^\sharp$ and
$\sin^{(m^\sharp)}$ denotes the $m^{\sharp\,\,\text{th}}$
derivative of the sine function.

Also,
$$\pa_x^{k^b}\left\{\frac{1}{\cosh(\hat{z}_2(x))-\cos(\hat{z}_1(x))}\right\}$$
is a linear combination, with harmless coefficients, of terms
$$\left[\cosh(\hat{z}_2(x))-\cos(\hat{z}_1(x))\right]^{-(m+1)}
\prod_{j=1}^m\pa_x^{k_j}\{\cosh(\hat{z}_2(x))-\cos(\hat{z}_1(x))\},$$
with $k_j\geq1$ for each $j$, and with
$k_1+\cdot\cdot\cdot+k_m=k^b$.

Therefore, the integrand in (\ref{docetres}) is a linear
combination, with harmless coefficients, of terms
$$\left[z_\mu^{(k'+1)}(x)\right]\left[\sin^{(m^\sharp)}(\hat{z}_1(x))\right]
\prod_{j=1}^{m^\sharp}\left(\hat{z}_1^{(k_j^\sharp)}(x)\right)$$$$\times\left[\cosh(\hat{z}_2(x))-\cos(\hat{z}_1(x))\right]^{-(m+1)}
\prod_{j=1}^m\pa_x^{k_j}\{\cosh(\hat{z}_2(x))-\cos(\hat{z}_1(x))\},$$
with $k_j$, $k_j^\sharp \geq 1$ for each $j$, and with
$k^\sharp_1+\cdot\cdot\cdot+k^\sharp_{m^\sharp}+k_1+\cdot\cdot\cdot+k_m+k'=\underline{k}$.

Furthermore, for $k_j\geq 1$,
$$\pa_x^{k_j}\{\cosh(\hat{z}_2(x))\}$$
is a sum with harmless coefficients of terms
$$\cosh^{(m_j)}(\hat{z}_2(x))\prod_{i=1}^{m_j}\left(
\hat{z}_2^{k_{ji}}(x)\right),$$ with $m_j\geq 1$, with $k_{ji}\geq
1$, and with $\sum_{i=1}^{m_j}k_{ji}=k_j$.

Similarly, $$\pa_x^{k_j}\{\cos(\hat{z}_1(x))\}$$ is a  sum with
harmless coefficients of terms
$$\cos^{(m_j)}(\hat{z}_1(x))\prod_{i=1}^{m_j}\left(\hat{z}_1^{(k_{ji})}(x)\right),$$
with $m_j\geq 1$, with each $k_{ji}\geq1$, and with
$\sum_{i=1}^{m_j}k_{ji}=k_j$.

Consequently, the integrand in (\ref{docetres}) is a linear
combination with harmless coefficients, of terms
$$\left[z_\mu^{(k'+1)}(x)\right]\left\{\left[\sin^{(m^\sharp)}(\hat{z}_1(x))\right]
\prod_{j=1}^{m^\sharp}\left[\hat{z}_1^{(k_j^\sharp)}(x)\right]\right\}$$$$\times\left[\cosh(\hat{z}_2(x))-\cos(\hat{z}_1(x))\right]^{-(m+1)}
\prod_{j=1}^{m'}\left\{\left[\cosh^{(m_j')}(\hat{z}_2(x))\right]\prod_{i=1}^{m_j'}\left[\hat{z}_2^{(k'_{ji})}(x)\right]\right\}$$
\begin{equation}\label{docecuatro}
\times\prod_{j=1}^{m''}\left\{\left[\cos^{(m_j'')}(\hat{z}_1(x))\right]\prod_{i=1}^{m_j''}\left[\hat{z}_1^{(k''_{ji})}(x)\right]\right\},\end{equation}
where all the
\begin{equation}\label{docecincoa}
\text{$k^\sharp_j$, $k'_{ji}$, $k''_{ji}\geq 1$, all the
$m'_j$, $m''_j\geq 1$,}
 \end{equation}
 and
\begin{equation}\label{docecincob}
k'+\sum_{j=1}^{m^\sharp}k^\sharp_j+\sum_{j=1}^{m'}\sum_{i=1}^{m_j'}k'_{ji}
+\sum_{j=1}^{m''}\sum_{i=1}^{m_j''}k''_{ji}=\underline{k},\quad
m'+m''=m.\end{equation} Let us call each expression in square
brackets in (\ref{docecuatro}) a BASIC FACTOR. Thus, the integrand
in (\ref{docetres}) is a linear combination, with harmless
coefficients, of products of basic factors. We assign to each
basic factor a weight, as follows.
\begin{itemize}
\item $[\hat{z}_1^{(k)}(x)]$ has weight 1. \item
$[\hat{z}_2^{(k)}(x)]$ has weight 1. \item
$[\sin^{(m^\sharp)}(\hat{z}_1(x))]$ has weight 1 if $m^\sharp$ is
even, 0 if $m^\sharp$ is odd. \item $[\cos^{(m''_j)}(\hat{z}_1(x))]$
has weight 1 if $m''_j$ is odd, 0 if $m''_j$ is even.  \item
$[\cosh^{(m'_j)}(\hat{z}_2(x))]$ has weight 1 if $m'_j$ is odd, 0 if
$m'_j$ is even. \item
$[\cosh(\hat{z}_2(x)-\cos(\hat{z}_1(x))]^{-(m+1)}$ has weight
$-2(m+1)$.
\end{itemize}
The significance of this notion is that if $\phi(x)$ is a basic
factor with weight $w$, then $|\phi(x)|=O(||y||^w)$ when $||y||$
is small, uniformly in $x$. That is because
$\hat{z}_\mu^{(k)}=\pa_x^k z_\mu(x,t)-\pa_x^k z_\mu(x+y,t)$ and
because of the COMPLEX CHORD-ARC CONDITION.

We define the TOTAL WEIGHT of a product of basic factors to be the
sum of the weights of those basic factors.

 We now check that
(\ref{docecincoa}) and \eqref{docecincob} implies that the sum of the weight of the basic
factors in (\ref{docecuatro}) is greater than or equal to zero.

To see this, we note that the total weight of
$$\left\{\cosh^{(m'_j)}(\hat{z}_2(x))\prod_{i=1}^{m'_j}\left(
\hat{z}_2^{k'_{ji}}(x)\right)\right\}$$ is
$1_{m_j'\,\text{odd}}+m'_j\geq 2$, since $m'_j\geq1$ by
(\ref{docecincoa}).

Similarly, the total weight of
$$\left\{\cos^{(m''_j)}(\hat{z}_1(x))\prod_{i=1}^{m''_j}\left(\hat{z}_1^{(k''_{ji})}(x)\right)\right\}$$
is $1_{m_j''\,\,\text{odd}}+m_j''\geq 2$, since $m_j''\geq 1$ by
(\ref{docecincob}).

Also, the total weight of
$$\left\{\left[\sin^{(m^\sharp)}(\hat{z}_1(x))\right]\prod_{j=1}^{m^\sharp}\left[\hat{z}_1^{(k_j^\sharp)}(x)\right]\right\},$$
is $1_{m^\sharp\,\,\text{even}}+m^\sharp\geq 1$, since $m^\sharp\geq0$.

Consequently, the total weight of the terms (\ref{docecuatro}) is
greater than or equal to
$$1+1-2(m+1)+2m'+2m''=0,$$
since $m'+m''=m$ by (\ref{docecincob}).

This completes the verification that the total weight of
(\ref{docecuatro}) is bigger than or equal to  zero.

We classify the terms (\ref{docecuatro}) as
follows.
\begin{itemize}
\item Suppose $[\hat{z}_1^{(\underline{k}+1)}(x)]$ or
$[\hat{z}_2^{(\underline{k}+1)}(x)]$ appears in the product
(\ref{docecuatro}). Then we say that the term (\ref{docecuatro})
is dangerous. \item If $[\hat{z}_1^{(\underline{k})}(x)]$ or
$[\hat{z}_2^{(\underline{k})}(x)]$ appears in a term
(\ref{docecuatro}), then we say that the term (\ref{docecuatro})
is safe. \item If the only $[\hat{z}_1^{(k)}(x)]$,
$[\hat{z}_2^{(k)}(x)]$ appearing in a term (\ref{docecuatro}) have
$k\leq \underline{k}-1$, then we call term (\ref{docecuatro})
easy.
\end{itemize}
Thus, any term (\ref{docecuatro}) (for which (\ref{docecincoa}) and \eqref{docecincob}
holds) is dangerous, safe, or easy. These conditions are mutually
exclusive thanks to (\ref{docecincoa}) and \eqref{docecincob} and because we will take  $\underline{k}\geq 4$.

Let us identify the dangerous and the safe terms of the form
(\ref{docecuatro}). Since
$$k'+\sum_{j=1}^{m^\sharp}k^\sharp_j+\sum_{j=1}^{m'}\sum_{i=1}^{m_j'}k'_{ji}
+\sum_{j=1}^{m''}\sum_{i=1}^{m_j''}k''_{ji}=\underline{k},$$ we
see that the only possible dangerous term (\ref{docecuatro}) is
the term arising from $k'=\underline{k}$, all $k^\sharp_j$,
$k'_{ji}$, $k''_{ji}=0$. Since all the $k'_{ji}$, $k''_{ji}\geq 1$
by (\ref{docecincoa}) and \eqref{docecincob}, we must have $k'=\underline{k}$,
$m^\sharp=0$, $m'=0$, $m''=0$. That is, the only dangerous term of
the form (\ref{docecuatro}) is
$$\left[\hat{z}_\mu^{(\underline{k}+1)}(x)\right]\left[\sin(\hat{z}_1(x))\right]\left[\frac{1}{\cosh(\hat{z}_2(x))-\cos(\hat{z}_1(x))}\right]$$
This term appears in the integrand of (\ref{docetres}) with a
coefficient 1.

Next, we look for all the safe terms (\ref{docecuatro}).

We distinguish several cases.

\emph{Case 1}: The basic factor $[\hat{z}_1^{(k'+1)}(x)]$ may have
$k'=\underline{k}-1$. In this case (\ref{docecincoa}) and \eqref{docecincob} imply that
there is exactly one nonzero number among the $k^\sharp_j$,
$k'_{ji}$, $k''_{ji}$; and that number is 1. Since all the
$k^\sharp_j$, $k'_{ji}$, $k''_{ji}$ are at least 1 by
(\ref{docecincoa}) and \eqref{docecincob}, there remain only the following possibilities:
\begin{itemize}
\item $m^\sharp=1$, $k_1^\sharp=1$, $m'=0$, $m''=0$ \item
$m^\sharp=0$, $m'=1$, $m_1'=1$, $k_{11}'=1$, $m''=0$\item
$m^\sharp=0$, $m'=0$, $m''=1$, $m_1''=1$, $k_{11}''=1$.
\end{itemize}
These possibilities give rise to the following term of the form
(\ref{docecuatro}):
\begin{itemize}
\item
$$\left[\hat{z}_\mu^{(\underline{k})}(x)\right][\cos(\hat{z}_1(x))]
\left[\frac{1}{\cosh(\hat{z}_2(x))-\cos(\hat{z}_1(x))}\right][\hat{z}_1'(x)]$$
\item
$$\left[\hat{z}_\mu^{(\underline{k})}(x)\right][\sin(\hat{z}_1(x))]
\left[\frac{1}{\cosh(\hat{z}_2(x))-\cos(\hat{z}_1(x))}\right]^2
[\sinh(\hat{z}_2(x))][\hat{z}_2'(x)]$$ \item
$$\left[\hat{z}_\mu^{(\underline{k})}(x)\right][\sin(\hat{z}_1(x))]
\left[\frac{1}{\cosh(\hat{z}_2(x))-\cos(\hat{z}_1(x))}\right]^2
[\sin(\hat{z}_1(x))][\hat{z}_1'(x)]$$
\end{itemize}
\emph{Case 2}: A basic factor $[\hat{z}_1^{(k_j^\sharp)}(x)]$ in
(\ref{docecuatro}) has $k^\sharp_j=\underline{k}$. Then
(\ref{docecincoa}) and \eqref{docecincob} show that $k'=0$, $m^\sharp=1$, $m'=0$,
$m''=0$, $m=0$. This gives rise to a term (\ref{docecuatro}) of
the form
\begin{itemize}
\item
$$\left[\hat{z}_\mu'(x)\right][\cos(\hat{z}_1(x))]
\left[\hat{z}_1^{(\underline{k})}(x)\right]
\left[\frac{1}{\cosh(\hat{z}_2(x))-\cos(\hat{z}_1(x))}\right]$$
\end{itemize}
\emph{Case 3}: A basic factor $[\hat{z}_2^{(k_{ji}')}(x)]$ in
(\ref{docecuatro}) has $k'_{ji}=\underline{k}$. Then
(\ref{docecincoa}) and \eqref{docecincob} show that $k'=0$, $m^\sharp=0$, $m'=1$,
$m_1'=1$, $k_{11}'=\underline{k}$, $m''=0$, and $m=1$. This gives
rise to a term (\ref{docecuatro}) of the form
\begin{itemize}
\item
$$\left[\hat{z}_\mu'(x)\right][\sin(\hat{z}_1(x))]
\left[\frac{1}{\cosh(\hat{z}_2(x))-\cos(\hat{z}_1(x))}\right]^2
[\sinh(\hat{z}_2(x))][\hat{z}_2^{(\underline{k})}(x)].$$
\end{itemize}
\emph{Case 4}: A basic factor $[\hat{z}_1^{(k_{ji}'')}(x)]$ in
(\ref{docecuatro}) has $k''_{ji}=\underline{k}$. Then
(\ref{docecincoa}) and \eqref{docecincob} show that $k'=0$, $m^\sharp=0$, $m'=0$,
$m''=1$, $m''_1=1$, $k_{11}''=\underline{k}$, $m=1$. This give
rise to a term (\ref{docecuatro}) of the form
\begin{itemize}
\item$$\left[\hat{z}_\mu'(x)\right][\sin(\hat{z}_1(x))]
\left[\frac{1}{\cosh(\hat{z}_2(x))-\cos(\hat{z}_1(x))}\right]^2
[\sin(\hat{z}_1(x))][\hat{z}_1^{(\underline{k})}(x)].$$
\end{itemize}
These are all the safe terms.

We have now identified all the dangerous and the safe terms.
All the remaining terms are easy.

We now substitute into (\ref{docetres}) the results of the above
discussion of the integrand of (\ref{docetres}). Recalling our
definition of $\hat{z}_1(x)$, $\hat{z}_1(x)$, and making the
change of variable $u=x+y$ in each of our integrals, we obtain the
formula:
$$
\pa_t\left[\pa_x^{\underline{k}}z_\mu(x,t)\right]$$
$$=\int_{u\in\T} \frac{\sin(z_1(x,t)-z_1(u,t))}
{\cosh(z_2(x,t)-z_2(u,t))-\cos(z_1(x,t)-z_1(u,t))}
\left[\pa_x^{\underline{k}+1}z_\mu(x,t)-\pa_u^{\underline{k}+1}z_\mu(u,t)\right]du
$$ $$+\text{(coeff)}\int_{u\in\T} \frac{[\pa_x
z_1(x,t)-\pa_uz_1(u,t)]\cos(z_1(x,t)-z_1(u,t))}
{\cosh(z_2(x,t)-z_2(u,t))-\cos(z_1(x,t)-z_1(u,t))}
\left[\pa_x^{\underline{k}}z_\mu(x,t)-\pa_u^{\underline{k}}z_\mu(u,t)\right]du
$$ $$+\text{(coeff)}\int_{u\in\T} \frac{[\pa_x
z_2(x,t)-\pa_uz_2(u,t)]\sinh(z_2(x,t)-z_2(u,t))\sin(z_1(x,t)-z_1(u,t))}
{(\cosh(z_2(x,t)-z_2(u,t))-\cos(z_1(x,t)-z_1(u,t)))^2} $$
$$\qquad\qquad\qquad\qquad\qquad\qquad\qquad\qquad\qquad\qquad\qquad\qquad\qquad
\times\left[\pa_x^{\underline{k}}z_\mu(x,t)-\pa_u^{\underline{k}}z_\mu(u,t)\right]du
$$
$$+\text{(coeff)}\int_{u\in\T} \frac{[\pa_x
z_1(x,t)-\pa_uz_1(u,t)][\sin(z_1(x,t)-z_1(u,t))]^2}
{(\cosh(z_2(x,t)-z_2(u,t))-\cos(z_1(x,t)-z_1(u,t)))^2}
\left[\pa_x^{\underline{k}}z_\mu(x,t)-\pa_u^{\underline{k}}z_\mu(u,t)\right]du$$
$$+\text{(coeff)}\int_{u\in\T} \frac{[\pa_x
z_\mu(x,t)-\pa_uz_\mu(u,t)]\cos(z_1(x,t)-z_1(u,t))}
{\cosh(z_2(x,t)-z_2(u,t))-\cos(z_1(x,t)-z_1(u,t))}
\left[\pa_x^{\underline{k}}z_1(x,t)-\pa_u^{\underline{k}}z_1(u,t)\right]du
$$ $$+\text{(coeff)}\int_{u\in\T} \frac{[\pa_x
z_\mu(x,t)-\pa_uz_\mu(u,t)]\sinh(z_2(x,t)-z_2(u,t))\sin(z_1(x,t)-z_1(u,t))}
{(\cosh(z_2(x,t)-z_2(u,t))-\cos(z_1(x,t)-z_1(u,t)))^2 }$$
$$\qquad\qquad\qquad\qquad\qquad\qquad\qquad\qquad\qquad\qquad\qquad\qquad\qquad
\times\left[\pa_x^{\underline{k}}z_2(x,t)-\pa_u^{\underline{k}}z_2(u,t)\right]du
$$ $$+\text{(coeff)}\int_{u\in\T} \frac{[\pa_x
z_\mu(x,t)-\pa_uz_\mu(u,t)][\sin(z_1(x,t)-z_1(u,t))]^2}
{(\cosh(z_2(x,t)-z_2(u,t))-\cos(z_1(x,t)-z_1(u,t)))^2}
\left[\pa_x^{\underline{k}}z_1(x,t)-\pa_u^{\underline{k}}z_1(u,t)\right]du$$
\begin{equation}\label{doceseis}
+\sum_{\nu=1}^{\nu_{\text{max}}}c^\mu_{\nu} \text{Easy Term}^\mu_\nu
(x,t),\end{equation} where $$\text{Easy Term}^\mu_\nu(x,t)$$
$$=\int_{u\in\T}\prod_{i=1}^{m_1}\left[\pa_x^{k_i}z_1(x,t)-\pa_u^{k_i}z_1(u,t)\right]
\prod_{i=1}^{m_2}\left[\pa_x^{k'_i}z_2(x,t)-\pa_u^{k'_i}z_2(u,t)\right]$$
$$\times
[\sin(z_1(x,t)-z_1(u,t))]^{m_3}[\cos(z_1(x,t)-z_1(u,t))]^{m_4}$$
$$\times
[\sinh(z_2(x,t)-z_2(u,t))]^{m_5}[\cosh(z_2(x,t)-z_2(u,t))]^{m_6}$$
\begin{equation}\label{docesiete}\times
\left(\cosh(z_2(x,t)-z_2(u,t))-\cos(z_1(x,t)-z_1(u,t))\right)^{-m_7}du,\end{equation}
where the $m_i$, $k_i$, $k'_i$ may depend on $\mu$, $\nu$, although we do
not indicate that in the notation, and where
\begin{equation}\label{doceocho}
\text{all the $k_i$ and $k'_i$ are $\geq 1$, and $m_7\geq
1$,}\end{equation} and
\begin{equation}\label{docenueve}
m_1+m_2+m_3+m_5-2m_7\geq 0 ;\end{equation} and
\begin{equation}\label{docediez}
\text{all $k_i$ and $k'_i$ are $\leq
\underline{k}-1$.}\end{equation} Now we recall that we assumed
$z_\mu(x,t)$ are analytic in $\Omega(t)$, and that
$(z_1(x,t),z_2(x,t))$ are assumed to satisfy the COMPLEX CHORD-ARC
CONDITION.

The integrands (\ref{doceseis}), (\ref{docesiete}) have no
singularities; they are holomorphic in $(x,u)\in \Omega(t)\times
\Omega(t)$. Hence, we may deform the contour in (\ref{doceseis}),
(\ref{docesiete}) from $\T$ to $\Gamma_\pm(t)$, and the formulas
holds, not merely for $x$ real, but for all $x$ in the closure of
$\Omega(t)$.

We record this result as a lemma.
\begin{lemma}\label{lemma6}
Let $\underline{k}\geq 4$.
Let $h(x,t)$ be a positive smooth function of $x\in \T$ for each
fixed $t$. Set $\Omega(t)=\{x\in \C/2\pi\Z\,:\,|\Im x|\leq h(\Re
x,t)\}$, $\Gamma_{\pm}=\{x\pm ih(x,t)\,:\,x\in\T\}$. Let
$z_1(x,t)$, $z_2(x,t)$ be smooth functions of $(x,t)$ on $\{
x\in\Omega(t)^{\text{closure}}$, $t\in$ Time interval$\}$.

For each fixed t, suppose $z_1(x,t)$, $z_2(x,t)$ are analytic in
$\Omega(t)$.

Assume that $z_1(x,t)$, $z_2(x,t)$ satisfy the Muskat equation
$$(ME):\,\,\pa_tz_\mu(x,t)=\int_{y\in\T}\frac{\sin\left(z_1(x,t)-z_1(x+y)\right)
\left(\pa_xz_\mu(x,t)-\pa_xz_\mu(x+y,t)\right)}
{\cosh\left(z_2(x,t)-z_2(x+y)\right)-\cos\left(z_1(x,t)-z_1(x+y)\right)}dy,$$
for $\mu=1,2$.

Assume also that $z_1(x,t)$, $z_2(x,t)$ satisfy the COMPLEX
CHORD-ARC CONDITION
$$(CCA):\,\,|\cosh\left(z_2(x,t)-z_2(y,t)\right)-\cos\left(z_1(x,t)-z_1(y,t)\right)|$$$$\geq
c_{CA}\left[||\Re (x-y)||+|\Im(x-y)|\right]^2.$$ for $x$,
$y\in\Omega(t)$.

Then, for $\zeta\in \Gamma_+(t)$, we have that
$$(I)\,\,\pa_t\left[\pa_\zeta^{\underline{k}}
z_\mu(\zeta,t)\right]=
\text{Dangerous}_\mu(\zeta,t)+\sum_{j=1}^6\text{Safe}^\mu_j(\zeta,t)+\sum_{\nu=1}^{\nu_{\text{max}}}c^\mu_\nu
\text{Easy}^\mu_\nu(\zeta,t),$$ where
\begin{align*}
(II)\,\,& \text{Dangerous}_\mu(\zeta,t)\\
&=\int\limits_{w\in\Gamma_+(t)}\frac{\sin(z_1(\zeta,t)-z_2(w,t))}
{\cosh(z_2(\zeta,t)-z_2(w,t))-\cos(z_1(\zeta,t)-z_1(w,t))}\\
&\qquad\qquad\qquad\qquad\qquad\qquad\qquad\qquad\qquad\qquad\qquad\times
(\pa_\zeta^{\underline{k}+1}z_\mu(\zeta,t)-\pa_w^{\underline{k}+1}z_\mu(w,t))dw
\end{align*}\begin{align*}
(III)\,\,& \text{Safe}^\mu_1(\zeta,t)\\
& =(\text{coeff})\int\limits_{w\in\Gamma_+(t)}\frac{(\pa_\zeta
z_1(\zeta,t)-\pa_w z_1(w,t))\cos(z_1(\zeta,t)-z_1(w,t))}
{\cosh(z_2(\zeta,t)-z_2(w,t))-\cos(z_1(\zeta,t)-z_1(w,t))}\\
&\qquad\qquad\qquad\qquad\qquad\qquad\qquad\qquad\qquad\qquad\qquad\qquad\times
(\pa_\zeta^{\underline{k}}z_\mu(\zeta,t)-\pa_w^{\underline{k}}z_\mu(w,t))dw
\end{align*}\begin{align*}(IV)\,\,&
\text{Safe}^\mu_2(\zeta,t)\\
&=(coeff)\int\limits_{w\in\Gamma_+(t)}\frac{(\pa_\zeta z_2(\zeta,t)-\pa_w
z_2(w,t))\sinh(z_2(\zeta,t)-z_2(w,t))\sin(z_1(\zeta,t)-z_1(w,t))}
{\left(\cosh(z_2(\zeta,t)-z_2(w,t))-\cos(z_1(\zeta,t)-z_1(w,t))\right)^2}\\
&\qquad\qquad\qquad\qquad\qquad\qquad\qquad\qquad\qquad\qquad\qquad\qquad\times
(\pa_\zeta^{\underline{k}}z_\mu(\zeta,t)-\pa_w^{\underline{k}}z_\mu(w,t))dw\end{align*}\begin{align*}
(V)\,\,& \text{Safe}^\mu_3(\zeta,t)\\
& = \text{(coeff)}\int\limits_{w\in\Gamma_+(t)}\frac{(\pa_\zeta
z_1(\zeta,t)-\pa_w z_1(w,t))\sin^2(z_1(\zeta,t)-z_1(w,t))}
{\left(\cosh(z_2(\zeta,t)-z_2(w,t))-\cos(z_1(\zeta,t)-z_1(w,t))\right)^2}\\
&\qquad\qquad\qquad\qquad\qquad\qquad\qquad\qquad\qquad\qquad\qquad\qquad\times
(\pa_\zeta^{\underline{k}}z_\mu(\zeta,t)-\pa_w^{\underline{k}}z_\mu(w,t))dw\end{align*}\begin{align*}
(VI)\,\,&\text{Safe}^\mu_4(\zeta,t)\\
& = \text{(coeff)}\int\limits_{w\in\Gamma_+(t)}\frac{(\pa_\zeta
z_\mu(\zeta,t)-\pa_wz_\mu(w,t))\cos(z_1(\zeta,t)-z_1(w,t))}
{\left(\cosh(z_2(\zeta,t)-z_2(w,t))-\cos(z_1(\zeta,t)-z_1(w,t))\right)}\\
&\qquad\qquad\qquad\qquad\qquad\qquad\qquad\qquad\qquad\qquad\qquad\qquad\times
(\pa_\zeta^{\underline{k}}z_1(\zeta,t)-\pa_w^{\underline{k}}z_1(w,t))dw\end{align*}\begin{align*}
(VII)\,\,& \text{Safe}^\mu_5\\
& = \text{(coeff)}\int\limits_{w\in\Gamma_+(t)}\frac{(\pa_\zeta
z_\mu(\zeta,t)-\pa_wz_\mu(w,t))\sin(z_1(\zeta,t)-z_1(w,t))\sinh(z_2(\zeta,t)-z_2(w,t))}
{\left(\cosh(z_2(\zeta,t)-z_2(w,t))-\cos(z_1(\zeta,t)-z_1(w,t))\right)^2}\\
&\qquad\qquad\qquad\qquad\qquad\qquad\qquad\qquad\qquad\qquad\qquad\qquad\times
(\pa_\zeta^{\underline{k}}z_2(\zeta,t)-\pa_w^{\underline{k}}z_2(w,t))dw\end{align*}\begin{align*}
(VIII)\,\,& \text{Safe}^\mu_6(\zeta,t)\\
& = \text{(coeff)}\int\limits_{w\in\Gamma_+(t)}\frac{(\pa_\zeta
z_\mu(\zeta,t)-\pa_w z_\mu(w,t))\sin^2(z_1(\zeta,t)-z_1(w,t))}
{\left(\cosh(z_2(\zeta,t)-z_2(w,t))-\cos(z_1(\zeta,t)-z_1(w,t))\right)^2}\\
&\qquad\qquad\qquad\qquad\qquad\qquad\qquad\qquad\qquad\qquad\qquad\qquad\times
(\pa_\zeta^{\underline{k}}z_1(\zeta,t)-\pa_w^{\underline{k}}z_1(w,t))dw
\end{align*}

and each term Easy$^\mu_\nu(\zeta,t)$ has the form
\begin{align*}
(IX)\,\, & \text{Easy}^\mu_\nu(\zeta,t)\\
& = \int\limits_{w\in \Gamma_+(t)} \prod_{i=1}^{m_1}(\partial_\zeta^{k_i}z_1(\zeta,t)-\pa_w^{k_i}z_1(w,t))
\prod_{i=1}^{m_2}(\partial_\zeta^{k'_i}z_2(\zeta,t)-\pa_w^{k'_i}z_2(w,t))\\
&\times [\sin(z_1(\zeta,t)-z_1(w,t))]^{m_3}[\cos(z_1(\zeta,t)-z_1(w,t))]^{m_4}\\
&\times [\sinh(z_2(\zeta,t)-z_2(w,t))]^{m_5}[\cosh(z_2(\zeta,t)-z_2(w,t))]^{m_6}\\
&\times [\left(\cosh(z_2(\zeta,t)-z_2(w,t))-\cos(z_1(\zeta,t)-z_1(w,t))\right)]^{-m_7}dw,
\end{align*}
where
\begin{align*}
(X)\,\, & m_7\geq 1,\quad 1\leq k_i\leq \underline{k}-1,\quad 1\leq k'_i\leq \underline{k}-1,
\end{align*}
and
\begin{align*}
(XI)\,\, & m_1+m_2+m_3+m_5-2m_7\geq 0;
\end{align*}
also
\begin{align*}
(XII)\,\, & \sum_ik_i+\sum_i k'_i=\underline{k}+1.
\end{align*}
Here, the $m_i$, $k_i$, $k_i'$ may depend on $\mu$ and $\nu$, though we have suppressed this dependence in the notation.
\end{lemma}

\subsection{Perturbing Easy and Safe terms}\label{seccion13}
Let $z=(z_\mu(x,t))_{\mu=1,2}$ and $\underline{z}=(\underline{z}_\mu(x,t))_{\mu=1,2}$ be two solutions of the Muskat
equation (ME) as in the lemma (\ref{lemma6}).

We assume that $z_\mu(x,t)$ and $\underline{z}_\mu(x,t)$ are defined on $\T\times T$, where the time interval $T$ contains some given $t_0\in[\tau^2,\tau]$.

We make the following assumptions on $z_\mu$ and $\underline{z}_\mu(x,t)$:
\begin{enumerate}
\item $z_\mu(x,t_0)$ and $\underline{z}_\mu(x,t_0)$ extend to functions analytic on $\Omega=\{|\Im x|< h(\Re x,t_0)\}$ and smooth on the closure of $\Omega$.
\item $$\int_{\zeta\in \Gamma_{\pm}}\left|\pa_\zeta^{14}\underline{z}_\mu(\zeta,t_0)\right|^2d\Re \zeta \leq C,$$
where $\Gamma_{\pm}=\{\Im \zeta=\pm h(\Re \zeta,t_0)\}$. Here, $h(x,t)$ is as in section \ref{seccion8}.
\item $$|\cosh\left(\underline{z}_2(\zeta,t_0)-\underline{z}_2(w,t_0)\right)-\cos\left(\underline{z}_1(\zeta,t_0)-\underline{z}_1(w,t_0)\right)|$$$$\geq
c_{CA}\left[||\Re (\zeta-w)||+|\Im(\zeta-w)|\right]^2,$$ for $\zeta$, $w\in \Omega$.
\item $$||z(\cdot,t)-\un(\cdot,t)||_{H^4(\Omega(t))}\leq \lambda,$$
with $\lambda$ as in section \ref{seccion8} .

Here we recall that
$$||f||^2_{H^4(\Omega(t))}=\sum_{\pm}\int_{\zeta\in\Gamma_{\pm}(t)}|f(\zeta)|^2d\Re\zeta+\sum_{\pm}\int_{\zeta\in\Gamma_{\pm}(t)}|\pa^4_\zeta f(\zeta)|^2d\Re\zeta.$$

\end{enumerate}

Note that our unperturbed solution $\underline{z}$ is assumed to have $4+10$ controlled derivatives in  hypothesis 2, whereas the perturbation $z-\underline{z}$ is assumed to have only $4$ controlled derivatives in the hypothesis 4. In this section, we write $c$, $C$, $C'$, etc. to denote constants determined by the constants $C$ and $c_{CA}$ in hypotheses 2 and 3. These symbols may denote different constants in different occurrences.

From the previous hypotheses, we conclude that
\begin{equation}\label{trececinco}
\int_{\zeta\in \Gamma_{\pm}}\left | \pa_\zeta^{4} z_\mu (\zeta,t_0)\right|^2 d\Re \zeta\leq C,
\end{equation}
\begin{equation}\label{treceseis}
\left|\pa_\zeta^j z_\mu(\zeta,t_0)\right|,\,\,\left|\pa_\zeta^j \underline{z}_\mu(\zeta,t_0)\right|\leq C\quad \text{for $0\leq j\leq 3$, $\zeta\in \Omega^{closure}$},\end{equation}
$$|\cosh\left(\underline{z}_2(\zeta,t_0)-\underline{z}_2(w,t_0)\right)-\cos\left(\underline{z}_1(\zeta,t_0)-\underline{z}_1(w,t_0)\right)|$$
\begin{equation}\label{trecesiete}
\geq
c_{CA}\left[||\Re (\zeta-w)||+|\Im(\zeta-w)|\right]^2, \quad \text{ for $\zeta$, $w\in \Omega^{closure}$}.\end{equation}

Let $Safe^\mu_j(\zeta,t_0)$ ($j=1,...,6)$) and $Easy^\mu_\nu(\zeta,t_0)$ be as in the lemma (\ref{lemma6}), arising from the Muskat solution $(z_1(\zeta,t),z_2(\zeta,t))$, taking $t=t_0$ in that lemma. Similarly,
let $\underline{Safe}^\mu_j(\zeta,t_0)$
and $\underline{Easy}^\mu_\nu(\zeta,t_0)$ arise from the Muskat solution $(\underline{z}_1(\zeta,t), \underline{z}_2(\zeta,t))$. Our goal in this section is to estimate
 $\left| Easy^\mu_\nu (\zeta, t_0)-\underline{Easy}^\mu_\nu (\zeta, t_0)\right|$ and $\left| Safe^\mu_j (\zeta, t_0)-\underline{Safe}^\mu_j (\zeta, t_0)\right|$. We begin with the $Easy$ terms.

From (1)...(5), we obtain the estimates
\begin{equation}\label{treceocho}
\left|\pa_\zeta^j (z_\mu-\underline{z}_\mu)(\zeta,t_0)\right|\leq C\lambda \quad \text{$0\leq j\leq 3$, $\zeta\in \Gamma_{\pm}$};\end{equation}
and
\begin{equation}\label{trecenueve}
\left|\pa_\zeta^{3} (z_\mu-\underline{z}_\mu)(\zeta,t_0)-\pa_w^{3} (z_\mu-\underline{z}_\mu)(w,t_0)\right|\leq C ||\Re (\zeta-w)||M\left[\pa_\zeta^{4}(z_\mu-\underline{z}_\mu)\right](\zeta)\end{equation}
for $\zeta$, $w\in \Gamma_+$, where $M$ is a Hardy-Littlewood maximal function.

We prepare to estimate
$$\left|\left[\cosh\left(z_2(\zeta,t_0)-z_2(w,t_0)\right)-\cos\left(z_1(\zeta,t_0)-z_1(w,t_0)\right)\right]\right.$$
$$-\left. \left [\cosh\left(\underline{z}_2(\zeta,t_0)-\underline{z}_2(w,t_0)\right)-\cos\left(\underline{z}_1(\zeta,t_0)-\underline{z}_1(w,t_0)\right)\right]\right|$$
for $\zeta$, $w\in\Gamma_+$.

To do so, let $\xi$, $\eta$, $\underline{\xi}$, $\underline{\eta}$ be complex numbers, with $|\xi|$, $|\eta|$, $|\underline{\xi}|$, $|\underline{\eta}| \leq C$. We set $\xi_s=s\xi+(1-s)\underline{\xi}$ and $\eta_s=s\eta+(1-s)\underline{\eta}$ for $0\leq s\leq 1$. We then have
$$\left|\frac{d}{ds}\left[\cosh(\xi_s)-\cos(\eta_s)\right]\right|=\left|\sinh(\xi_s)(\xi-\underline{\xi})+\sin(\eta_s)(\eta-\underline{\eta})\right|$$
$$\leq C' |\xi_s||\xi-\underline{\xi}|+C'|\eta_s||\eta-\underline{\eta}|\leq C' (|\xi|+|\underline{\xi}|)|\xi-\underline{\xi}|+C'(|\eta|+|\underline{\eta}|)|\eta-\underline{\eta}|. $$
We integrate in $s$.

We take
\begin{eqnarray*}
\xi &=&z_2(\zeta,t_0)-z_2(w,t_0),\\
\underline{\xi} &=&\underline{z}_2(\zeta,t_0)-\underline{z}_2(w,t_0),\\
\eta &=& z_1(\zeta,t_0)-z_1(w,t_0),\\
\underline{\eta}&=&\underline{z}_1(\zeta,t_0)-\underline{z}_1(w,t_0)
\end{eqnarray*}
 for $|\zeta-w|\leq c$, where $c$ is a small constant and thus obtain:
$$\left|\left[\cosh\left(z_2(\zeta,t_0)-z_2(w,t_0)\right)-\cos\left(z_1(\zeta,t)-z_1(w)\right)\right]\right.$$
$$-\left. \left[\cosh\left(\underline{z}_2(\zeta,t_0)-\underline{z}_2(w,t_0)\right)-\cos\left(\underline{z}_1(\zeta,t)-\underline{z}_1(w)\right)\right]\right|$$
$$\leq C' [|z_2(\zeta,t_0)-z_2(w,t_0)|+|\un_2(\zeta,t_0)-\un_2(w,t_0)|]\times\left|[z_2(\zeta,t_0)-z_2(w,t_0)]-[\un_2(\zeta,t_0)-\un_2(w,t_0)]\right|$$
\begin{equation}\label{trecediez}+ C' [|z_1(\zeta,t_0)-z_1(w,t_0)|+|\un_1(\zeta,t_0)-\un_1(w,t_0)|]\times\left|[z_1(\zeta,t_0)-z_1(w,t_0)]-[\un_1(\zeta,t_0)-\un_1(w,t_0)]\right|\end{equation}
for $\zeta$, $w\in \Gamma_+$, $|\zeta-w|\leq c$.

Since
\begin{eqnarray*}
&|z_1(\zeta,t_0)-z_1(w,t_0)|,\,
|\un_1(\zeta,t_0)-\un_1(w,t_0)|,\
|z_2(\zeta,t_0)-z_2(w,t_0)|,\
|\un_2(\zeta,t_0)-\un_2(w,t_0)|\\&\leq C||\Re(\zeta-w)||,
\end{eqnarray*}
for $\zeta$, $w\in \Gamma_+$; and since
$$|[z_1(\zeta,t_0)-z_1(w,t_0)]-[\un_1(\zeta,t_0)-\un_1(w,t_0)]|=|[z_1(\zeta,t_0)-\un_1(\zeta,t_0)]-[z_1(w,t_0)-\un_1(w,t_0)]|$$
$$\leq\max_{\xi\in \Omega}|\pa_\xi(z_1-\un_1)(\xi,t_0)|\cdot |\Re(\zeta-w)|\leq C' \lambda |\Re(\zeta-w)|$$
for $\zeta$, $w\in \Gamma_+$, $|\zeta-w|\leq c$, thanks to (\ref{treceocho}); and since, similarly,
$$|[z_2(\zeta,t_0)-z_2(w,t_0)]-[\un_2(\zeta,t_0)-\un_2(w,t_0)]|=|[z_2(\zeta,t_0)-\un_2(\zeta,t_0)]-[z_2(w,t_0)-\un_2(w,t_0)]|$$
$$\leq\max_{\xi\in \Omega}|\pa_\xi(z_2-\un_2)(\xi,t_0)|\cdot |\Re(\zeta-w)|\leq C' \lambda |\Re(\zeta-w)|$$
for $\zeta$, $w\in \Gamma_+$, $|\zeta-w|\leq c$, we conclude from (\ref{trecediez}) that
$$\left|\left[\cosh\left(z_2(\zeta,t_0)-z_2(w,t_0)\right)-\cos\left(z_1(\zeta,t_0)-z_1(w,t_0)\right)\right]\right.$$
$$-\left. \left[\cosh\left(\underline{z}_2(\zeta,t_0)-\underline{z}_2(w,t_0)\right)-\cos\left(\underline{z}_1(\zeta,t_0)-\underline{z}_1(w,t_0)\right)\right]\right|$$
\begin{equation}\leq C\lambda \left|\sin\left(\frac{\zeta-w}{2}\right)\right|^2\label{nombrada}\end{equation}
for $\zeta$, $w\in\Gamma_+$. We have proven this for $|\zeta-w|\leq c$. For $|\zeta-w|>c$, $\zeta$, $w\in \Gamma_+$ the estimate \eqref{nombrada} is trivial from \eqref{treceocho}.
Together with the chord-arc condition (\ref{trecesiete}), this yields the following estimate
$$\left|\left[\cosh\left(z_2(\zeta,t_0)-z_2(w,t_0)\right)-\cos\left(z_1(\zeta,t_0)-z_1(w,t_0)\right)\right]^{-1}\right.$$
$$-\left. \left[\cosh\left(\underline{z}_2(\zeta,t_0)-\underline{z}_2(w,t_0)\right)-\cos\left(\underline{z}_1(\zeta,t_0)-\underline{z}_1(w,t_0)\right)\right]^{-1}\right|$$
\begin{equation}\label{treceonce}\leq C\lambda \left|\sin\left(\frac{\zeta-w}{2}\right)\right|^{-2},\end{equation}
for $\zeta$, $w\in \Gamma_+$.

Next, from (\ref{treceocho}) we obtain easily the estimates for $\zeta$, $w\in \Gamma_+$:
\begin{eqnarray}
\left| \cos\left(z_1(\zeta,t_0)-z_1(w,t_0)\right)-\cos\left(\un_1(\zeta,t_0)-\un_1(w,t_0)\right)\right|&\leq& C\lambda\label{trecedoce}\\
\left| \cosh\left(z_2(\zeta,t_0)-z_2(w,t_0)\right)-\cosh\left(\un_2(\zeta,t_0)-\un_2(w,t_0)\right)\right|&\leq& C\lambda\label{trecetrece}.
\end{eqnarray}
Next, for $\zeta$, $w\in \Gamma_+$, we have
$$\left| \sinh\left(z_2(\zeta,t_0)-z_2(w,t_0)\right)-\sinh\left(\un_2(\zeta,t_0)-\un_2(w,t_0)\right)\right|$$
$$\leq C\left |[z_2-\un_2](\zeta,t_0)-[z_2-\un_2](w,t_0)\right|$$
$$\leq C' \lambda ||\Re(\zeta-w)||$$
thanks to (\ref{treceocho}).

Because of the periodicity we conclude that
$$\left| \sinh\left(z_2(\zeta,t_0)-z_2(w,t_0)\right)-\sinh\left(\un_2(\zeta,t_0)-\un_2(w,t_0)\right)\right|$$
\begin{equation}\label{trececatorce}\leq C\lambda\left|\sin\left(\frac{\zeta-w}{2}\right)\right|\end{equation}
for $\zeta$, $w\in\Gamma_+$.

Similarly,

$$\left| \sin\left(z_1(\zeta,t_0)-z_1(w,t_0)\right)-\sin\left(\un_1(\zeta,t_0)-\un_1(w,t_0)\right)\right|$$
\begin{equation}\label{trecequince}\leq C\lambda\left|\sin\left(\frac{\zeta-w}{2}\right)\right|\end{equation}
for $\zeta$, $w\in\Gamma_+$.

From (\ref{treceocho}), we see that (for $\zeta$, $w\in \Gamma_+$) we have
$$\left|\left[\pa_\zeta^k z_\mu(\zeta,t_0)-\pa_w^k z_\mu (w,t_0)\right]-\left[\pa_\zeta^k \un_\mu(\zeta,t_0)-\pa_w^k \un_\mu (w,t_0)\right]\right|$$
\begin{equation}\label{trecedieciseis}
\leq C\lambda\left|\sin\left(\frac{\zeta-w}{2}\right)\right|,\end{equation}
when $1\leq k\leq 2$.

From (\ref{trecenueve}), we have (for $\zeta$, $w\in \Gamma_+$) that
$$\left|\left[\pa_\zeta^{3} z_\mu(\zeta,t_0)-\pa_w^{3} z_\mu (w,t_0)\right]-\left[\pa_\zeta^{3} \un_\mu(\zeta,t_0)-\pa_w^{3} \un_\mu (w,t_0)\right]\right|$$
\begin{equation}\label{trecediecisiete}
\leq C\left|\sin\left(\frac{\zeta-w}{2}\right)\right|M\left[\pa_\zeta^{4}(z_\mu-\un_\mu)\right](\zeta,t_0).\end{equation}
We now recall the form (IX) of the term $Easy^\mu_\nu(\zeta,t)$ in  lemma (\ref{lemma6})with $\underline{k}=4$. Changing the factors in (IX) one at a time to their analogues for $\un$, and applying estimates (\ref{treceonce})...(\ref{trecediecisiete}), we learn that
$$\left|Easy^\mu_\nu(\zeta,t_0)-\underline{Easy}^\mu_\nu(\zeta,t_0)\right|$$
$$\leq \int_{w\in\Gamma_+}\left\{ C\lambda+\sum_j C M\left[\pa_\zeta^{4}(z_j(\cdot,t_0)-\un_j(\cdot,t_0))\right](\zeta)\right\}d \Re w$$
\begin{equation}\label{trecedieciocho}\leq C\lambda+\sum_\mu C M\left[\pa_\zeta^{4}(z_\mu(\cdot,t_0)-\un_\mu(\cdot,t_0))\right](\zeta),\end{equation}
for $\zeta\in \Gamma_+$; here we use also (XI) in  lemma (\ref{lemma6}) to see that we do not pick up negative powers of $|\sin(\frac{\zeta-w}{2})|$.

It follows from 4, (\ref{trecedieciocho}) and the Hardy-Littlewood maximal  theorem that
\begin{equation}\label{trecediecinueve}
\int_{\zeta\in\Gamma_+}\left| Easy^\mu_\nu(\zeta,t_0)-\underline{Easy^\mu}_\nu(\zeta,t_0)\right|^2 d \Re \zeta\leq C\lambda^2.\end{equation}
Next, we estimate
$$\left| Safe^\mu_j(\zeta,t_0)-\underline{Safe}^\mu_j(\zeta,t_0)\right|$$
for $j=1,...,6$; see lemma (\ref{lemma6}) for the definitions of these quantities.

For instance let us examine $Safe^\mu_2(\zeta,t_0)$. Recalling (IV) from  lemma (\ref{lemma6}), we have
\begin{equation}\label{treceveinte}
Safe^\mu_2(\zeta,t_0)=c\int_{w\Gamma_+}K(\zeta,w)\left [\pa_\zeta^{4} z_\mu(\zeta,t_0)-\pa_w^{4} z_\mu(w,t_0)\right]dw,\end{equation}
for $\zeta\in \Gamma_+$, where
\begin{eqnarray}
&K(\zeta,w)=\nonumber\\
&\left[\pa_\zeta z_2(\zeta,t_0)-\pa_wz_2(w,t_0)\right]\cdot\left[\sinh(z_2(\zeta,t_0)-z_2(w,t_0))\right]\cdot\left[\sin\left(z_1(\zeta,t_0)-z_1(w,t_0)\right)\right]\nonumber\\
&\left[\cosh\left(z_2(\zeta,t_0)-z_2(w,t_0)\right)-\cos\left(z_1(\zeta,t_0)-z_1(w,t_0)\right)\right]^{-2}\nonumber\\
&\equiv [Factor1]\cdot[Factor2]\cdot[Factor3]\cdot[Factor4]\label{treceveintiuno}.
\end{eqnarray}
Similarly, starting from the Muskat solution $\un$ in place of $z$, we obtain
\begin{equation}\label{treceveintidos}
\underline{Safe}^\mu_2(\zeta,t_0)=c\int_{w\in\Gamma_+}\underline{K}(\zeta,w)\left [\pa_\zeta^{4} \un_\mu(\zeta,t_0)-\pa_w^{4} \un_\mu(w,t_0)\right]dw,\end{equation}
for $\zeta\in \Gamma_+$, where
\begin{equation}\label{treceveintitres}
\underline{K}(\zeta,w)=[\underline{Factor1}]\cdot[\underline{Factor2}]\cdot[\underline{Factor3}]\cdot[\underline{Factor4}].
\end{equation}
We write
$$Safe^\mu_2(\zeta,t_0)-\underline{Safe^\mu}_2(\zeta,t_0)$$
$$=c\int_{w\in \Gamma_+}K(\zeta,w)\left [\pa_\zeta^{4}(z_\mu-\un_\mu)(\zeta,t_0)-\pa_w^{4}(z_\mu-\un_\mu)(w,t_0)\right]dw$$
$$+c\int_{w\in\Gamma_+}\left[K(\zeta,w)-\underline{K}(\zeta,w)\right]\left[\pa_\zeta^{4}\un_\mu(\zeta,t_0)-\pa_w^{4}\un_\mu(w,t_0)\right]dw$$
\begin{equation}\label{treceveinticuatro}
\equiv Term1(\zeta,t_0)+Term2(\zeta,t_0),\end{equation}
for $\zeta\in \Gamma_+$.

It is convenient that the $4-$th derivatives of $\un_\mu$ enter into $Term2$, rather than the $4-$th derivatives of $z_\mu$.

From our estimates (\ref{treceseis}) and (\ref{trecesiete}),  we see that, for $|\zeta-w|\leq c$,
\begin{eqnarray}
Factor1 &=& A_1(\zeta,w)\cdot (\zeta-w)\nonumber\\
Factor2 &=& A_2(\zeta,w)\cdot (\zeta-w)\nonumber\\
Factor3 &=& A_3(\zeta,w)\cdot (\zeta-w)\nonumber\\
Factor4 &=& A_4(\zeta,w)\cdot (\zeta-w)^{-4}\nonumber\\
\end{eqnarray}
with the derivatives of $A_1$,...,$A_4$ up to order 1 less or equal than $C$ in absolute value.

For $||\Re(\zeta-w)||>c$, $Factor1$, $Factor2$, $Factor3$ and $Factor4$ we have derivatives up to order 1 bounded by $C$.

Therefore,
$$K(\zeta,w)=A(\zeta,w)\cdot (\zeta-w)^{-1}$$
 for $|\zeta-w|\leq c$, with the $C^1-$ norm of $A(\zeta,w)$ on $\Gamma_+\times\Gamma_+$ bounded by $C$. On $||\Re(\zeta-w)||>c$ the $C^
  1-$norm of the kernel $K(\zeta,w)$ is at most C. Recalling how $K(\zeta,w)$ enters into $Term1$ in (\ref{treceveinticuatro}), we see that $Term1(\zeta,t_0)$ is a singular integral operator applied to
$$\pa^{4}_\zeta(z_\mu-\un_\mu)(\cdot,t_0).$$
Consequently,
\begin{equation}\label{treceveinticinco}
\left(\int_{\zeta\in \Gamma_+}\left|Term1(\zeta,t_0)\right|^2d\Re \zeta\right)^{\frac{1}{2}}\leq C\lambda.\end{equation}

To estimate $Term2$ in (\ref{treceveinticuatro}), we first estimate
$$\left|K(\zeta,w)-\underline{K}(\zeta,w)\right|$$
$$=\left|[Factor1]\cdots[Factor4]-[\underline{Factor1}]\cdots[\underline{Factor4}]\right|$$
$$=\left|\sum_{j=1}^4 \prod_{i<j}[Factori]\cdot [Factorj-\underline{Factorj}]\cdot \prod_{i>j}[\underline{Factori}]\right|$$
\begin{equation}\leq \sum_{j=1}^4 \prod_{i<j}|Factori|\cdot |Factorj-\underline{Factorj}|\cdot \prod_{i>j}|\underline{Factori}|.\label{4.331/2}\end{equation}
Our estimates (\ref{treceseis}), (\ref{treceocho}), \eqref{trecedieciseis} and the definition of $[Factor1]$ in (\ref{treceveintiuno}) (and its analogue for $[\underline{Factor1}]$) show that
\begin{eqnarray}
|Factor1|,\,|\underline{Factor1}|&\leq& C|\zeta-w|,\nonumber\\
|Factor1-\underline{Factor1}|&\leq& C\lambda |\zeta-w|.\label{treceveintiseis}\end{eqnarray}
Similarly, (\ref{treceseis}), (\ref{trececatorce}) yield the estimates
\begin{eqnarray}
|Factor2|,\,|\underline{Factor2}|&\leq& C|\zeta-w|,\nonumber\\
|Factor2-\underline{Factor2}|&\leq& C\lambda |\zeta-w|.\label{treceveintisiete}\end{eqnarray}
Also, (\ref{treceseis}) and (\ref{trecequince}) tell us that
\begin{eqnarray}
|Factor3|,\,|\underline{Factor3}|&\leq& C|\zeta-w|,\nonumber\\
|Factor3-\underline{Factor3}|&\leq& C\lambda |\zeta-w|.\label{treceveintiocho}\end{eqnarray}
From 3, (\ref{trecesiete}) and (\ref{treceonce}), we have
\begin{eqnarray}
|Factor4|,\,|\underline{Factor4}|&\leq& C|\zeta-w|^{-4},\nonumber\\
|Factor4-\underline{Factor4}|&\leq& C\lambda |\zeta-w|^{-4},\label{treceveintinueve}\end{eqnarray}
when $|\zeta-w|\leq 10^{-2}$.

When $||\Re (\zeta-w)||\geq c$, 3, (\ref{trecesiete}), (\ref{treceonce}) yield
\begin{eqnarray}
|Factor4|,\,|\underline{Factor4}|&\leq& C,\nonumber\\
|Factor4-\underline{Factor4}|&\leq& C\lambda,\label{trecetreinta}\end{eqnarray}

Putting (\ref{treceveintiseis}),...,(\ref{trecetreinta}) into (\ref{4.331/2}), we see that
$$\left|K(\zeta,w)-\underline{K}(\zeta,w)\right|\leq C\lambda |\zeta-w|^{-1},$$
if $|\zeta-w|\leq 10^{-2}$. Also,
$$\left|K(\zeta,w)-\underline{K}(\zeta,w)\right|\leq C\lambda,$$
if $|\zeta-w|\geq c$.
On the other hand, 2 gives
$$\left|\pa_\zeta^{4}\un_\mu(\zeta,t_0)-\pa_w^{4}\un_\mu(w,t_0)\right|\leq C|\zeta-w|.$$
Consequently,
$$\left|K(\zeta,w)-\underline{K}(\zeta,w)\right|\cdot\left|\pa_\zeta^{4}\un_\mu(\zeta,t_0)-\pa_w^{4}\un_\mu(w,t_0)\right|\leq C\lambda$$
for $\zeta$, $w\in \Gamma_+$.

Hence, by definition of $Term2(\zeta,t_0)$ in (\ref{treceveinticuatro}), we have
$$\left|Term2(\zeta,t_0)\right|\leq C\lambda $$
for $\zeta\in \Gamma_+$.

Together with (\ref{treceveinticuatro}) and (\ref{treceveinticinco}), this yields the estimate
\begin{equation*}
\int_{\zeta\in\Gamma_+}\left| Safe^\mu_2(\zeta,t_0)-\underline{Safe}^\mu_2(\zeta,t_0)\right|^2 d \Re \zeta\leq C\lambda^2.\end{equation*}
The other $Safe^\mu_j(\zeta,t_0)-\underline{Safe}^\mu_j(\zeta,t_0)$ may be handled similarly. Thus, we have the estimate
 \begin{equation}\label{trecetreintaidos}
\int_{\zeta\in\Gamma_+}\left| Safe^\mu_j(\zeta,t_0)-\underline{Safe}^\mu_j(\zeta,t_0)\right|^2 d \Re \zeta\leq C\lambda^2.\end{equation}

\subsection{Auxiliary Functions}\label{seccion14}

We assume we are in the setting of  section (\ref{seccion13}). We define
$$\underline{a}(\zeta,w,t)=$$
\begin{equation}\label{catorceuno}
\frac{\sin\left(\un_1(\zeta,t)-\un_1(w,t)\right)}{\cosh\left(\un_2(\zeta,t)-\un_2(w,t)\right)-\cos\left(\un_1(\zeta,t)-\un_1(w,t)\right)}-\frac{\pa_\zeta\un_1(\zeta,t)}{[\pa_\zeta\un_1(\zeta,t)]^2+[\pa_\zeta\un_2(\zeta,t)]^2}\cot\left(\frac{\zeta-w}{2}\right)
\end{equation}
for $\zeta$, $w\in\Omega^\text{closure}(t)$, where
$$\Gamma_+(t)=\{ \zeta\in \C\,:\, \zeta=x+ih(x,t),\,\, x\in \T\}$$
and
$$\Omega(t)=\{ \zeta\in \C\,:\, |\Im \zeta|< h(\Re \zeta,t),\,\, x\in \T\}.$$
Here $t\in[\tau^2,\tau]$. If $t\in [-\tau^2,\tau^2]$ we replace $h(x,t)$ by $\hbar(x,t)$.

Note that we refer to the unperturbed solution $\un$ in defining $\underline{a}$.

We recall the following properties  of the unpertubed solutions $\un$ (see section \ref{seccion7}):
\begin{equation}\label{catorcedos}
\text{$\un_1(\zeta,t)$, $\un_2(\zeta,t)$ are $C^3$-function on $\Omega(t)^\text{closure},\,\,t\in [-\tau^2,\tau]$ and analytic in $\zeta\in \Omega(t)$ for fixed t.}
\end{equation}
\begin{equation}\label{catorcetres}
|\un_1(\zeta,t)|,\, |\un_2(\zeta,t)|\leq C \quad \text{for $|\Im \zeta|\leq c_0,\,t\in [-\tau^2,\tau]$}.\end{equation}
\begin{equation}\label{catorcecuatro}
\text{$\un_1(\zeta,0)$, $\un_2(\zeta,0)$ are odd functions of $\zeta$ and real for real $\zeta$.}\end{equation}

We also recall the COMPLEX CHORD-ARC CONDITION
$$|\cosh\left(\underline{z}_2(\zeta,t)-\underline{z}_2(w,t)\right)-\cos\left(\underline{z}_1(\zeta,t)-\underline{z}_1(w,t)\right)|$$
\begin{equation}\label{catorcecinco}
\geq
c_{CA}\left[||\Re (\zeta-w)||+|\Im(\zeta-w)|\right]^2\end{equation}
for  $\zeta,\,w \in \Omega(t)$, $t\in [-\tau^2,\tau]$.

We define
\begin{equation}\label{catorceseis}
\underline{\tilde{a}}(\zeta,t)=\int_{w\in\Gamma_+(t)}\underline{a}(\zeta,w,t)dw.
\end{equation}
We notice that since $\underline{a}(\zeta,w,t)$ is analytic in $w\in \Omega(t)$
we can deform the contour of integration in \eqref{catorceseis} to obtain
\begin{equation}\label{catorceseisa}
\underline{\tilde{a}}(\zeta,t)=\int_{w\in\T}\underline{a}(\zeta,w,t) dw.
\end{equation}

Let us assume that
\begin{equation}\label{catorceonce}
\underline{\tilde{a}}(\zeta,t)\in C^2(\Omega(t)^{\text{closure}})\quad \text{with} \quad ||\underline{\tilde{a}}||_{C^2(\Omega(t)^{\text{closure}})}\leq C,
\end{equation}
where $C$ is a controlled constant. We will prove this fact in the next section (see \eqref{paraz}).

For $\zeta=0$, $t=0$, we have $\un_\mu(0,0)=0$; hence $\underline{a}(0,w,0)$ is an odd function of $w$, thanks to (\ref{catorcecuatro}) again.
Integrating over $w\in [-\pi,\pi]$ and comparing with (\ref{catorceseisa}), we see that
\begin{equation}\label{catorcesiete}
\underline{\tilde{a}}(0,0)=0.
\end{equation}

From (\ref{catorcesiete}), (\ref{catorceonce}) we obtain the estimate
\begin{equation}\label{catorcedoce}
|\underline{\tilde{a}}(\zeta,t)|\leq C|\zeta|+C|t|\end{equation}
for $\zeta\in \Omega(t)^\text{closure}$, $t\in[-\tau^2,\tau]$.

We take $\zeta=x+ih(x,t)$ with $t\in[\tau^2,\tau]$.

Then (\ref{catorcedoce}) gives
\begin{equation}\label{catorcetrece}
|\underline{\tilde{a}}(x+ih(x,t),t)|\leq C|x|+Ch(x,t)+Ct.\end{equation}
From section (\ref{seccion8}) follows
$$|\pa_xh(x,t)|\leq CA^{-1}|x|$$
for $x\in[-\pi,\pi]$, $t\in[0,\tau]$. Hence, (\ref{catorcetrece}) yields
$$|\pa_x h(x,t)|\cdot|\underline{\tilde{a}}(x+ih(x,t),t)|\leq CA^{-1}x^2+CA^{-1}|x|t+CA^{-1}h(x,t)$$
\begin{equation}\label{catorcecatorce}\leq CA^{-1}x^2+CA^{-1}t^2+CA^{-1}h(x,t).\end{equation}
Now
$$A^{-1}x^2\leq Ch(x,t)$$
for $x\in[-\pi,\pi]$, $t\in[\tau^2,\tau]$. Therefore, (\ref{catorcecatorce}) yields
\begin{equation}\label{catorcequince}
|\pa_xh(x,t)\underline{\tilde{a}}(x+ih(x,t),t)|\leq Ch(x,t)+CA^{-1}\tau^2\end{equation}
for $x\in  \T$, $t\in[\tau^2,\tau]$. (Compare with (\ref{ochoseis}), (\ref{ochosiete})).

Also, since we assume $\un_\mu(\zeta,t)$ are real for real $\zeta$, we see from (\ref{catorceuno}) that $\underline{a}(\zeta,w,t)$ is real for real $\zeta$ and real $w$. Hence, (\ref{catorceseisa}) shows that $\underline{\tilde{a}}(\zeta,t)$ is real for real $\zeta$. Together with (\ref{catorceonce}), this yields
\begin{equation}\label{catorcedieciseis}
|\Im \underline{\tilde{a}}(x+ih(x,t),t)|\leq Ch(x,t)\end{equation}
for $x\in \T$, $t\in [\tau^2,\tau]$.

Let us derive an analogues estimate for the function $\hbar(x,t)$ in section (\ref{seccion8}). The proof of (\ref{catorcedieciseis}) yields also
\begin{equation}\label{catorcediecisiete}
 |\Im \underline{\tilde{a}}(x+i\hbar(x,t),t)|\leq C\hbar(x,t)\end{equation}
for $x\in \T$, $t\in [-\tau^2,\tau^2]$.

Also (\ref{catorcedoce}) gives
$$|\underline{\tilde{a}}(x+i\hbar(x,t),t)|\leq C|x|+C\hbar(x,t)+C|t|$$
for $x\in [-\pi,\pi]$, $t\in[-\tau^2,\tau^2]$.

From (\ref{ochoveintidos}), we see that
$$|\pa_x \hbar(x,t)|\leq CA^{-1}|x|$$
for the same range of $x$, $t$.

Hence,
$$|\pa_x\hbar(x,t)\underline{\tilde{a}}(x+i\hbar(x,t),t)|\leq CA^{-1}x^2+CA^{-1}|x||t|+C\hbar(x,t)$$
$$\leq CA^{-1}x^2+CA^{-1}t^2+C\hbar(x,t)$$
for $x\in[-\pi,\pi]$, $t\in[-\tau^2,\tau^2]$.

Since
$$CA^{-1}x^2+CA^{-1}t^2\leq CA^{-1}x^2+CA^{-1}\tau^4\leq C\hbar(x,t)$$
for $x\in[-\pi,\pi]$, $t\in[-\tau^2,\tau^2]$ (thanks to (\ref{ochoveintidos})), it follows that
\begin{equation}\label{catorcedieciocho}
|\pa_x \hbar(x,t) \underline{\tilde{a}}(x+i\hbar(x,t),t)|\leq C\hbar(x,t)\end{equation}
for $x\in \T$, $t\in [-\tau^2,\tau^2]$.

Our basic results on the size of $\underline{\tilde{a}}$ and its imaginary part are (\ref{catorcequince}),...,(\ref{catorcedieciocho}).

\subsection{Perturbing Auxiliary Functions}\label{seccion15}

Let $z$, $\un$ be two solutions of Muskat, that satisfy the assumptions in the section (\ref{seccion13}). We also assume that $\un$ satisfies the additional assumptions made in the section (\ref{seccion14}).

We define $\underline{a}(\zeta,w,t)$ and $\underline{\tilde{a}}(\zeta,t)$ as in (\ref{catorceuno}) and (\ref{catorceseis}).

For all $t$ for which the perturbed Muskat solution $z(\cdot,t)$ is defined, we let $a(\zeta,w,t)$ and $\tilde{a}(\zeta,t)$ be the analogous expressions for $z$, i.e.,
$$a(\zeta,w,t)=$$
\begin{equation}\label{quinceuno}
\frac{\sin\left(z_1(\zeta,t)-z_1(w,t)\right)}{\cosh\left(z_2(\zeta,t)-z_2(w,t)\right)-\cos\left(z_1(\zeta,t)-z_1(w,t)\right)}-\frac{\pa_\zeta z_1(\zeta,t)}{[\pa_\zeta z_1(\zeta,t)]^2+[\pa_\zeta z_2(\zeta,t)]^2}\cot\left(\frac{\zeta-w}{2}\right)
\end{equation}
for
$$|\Im \zeta|\leq h(\Re\zeta,t),\, |\Im w |\leq h(\Re w,t)\quad \text{if } t\in[\tau^2,\tau]$$
$$|\Im \zeta|\leq \hbar(\Re\zeta,t),\, |\Im w |\leq \hbar(\Re w,t)\quad \text{if } t\in[-\tau^2,\tau^2].$$
We will show that, since $z(\zeta,t)\in C^3(\Omega(t)^{closure})$ the function

\begin{equation}\label{catorcediez}
\text{$a(\zeta,w,t)$ is $C^1$ in $w$ for $\zeta$, $w\in\Omega(t)^{closure}$ and has $C^1$-norm at most $C$}.
\end{equation}
To prove \eqref{catorcediez},  we may restrict attention to
$$\underline{\Omega}=\{\zeta,\,w \in \Omega(t)^{\text{closure}},\, t\in[0,\tau],\, |\zeta-w|\leq c\},$$
for small $c$. Outside this region there is no problem, thanks to the complex chord-arc condition (\ref{catorcecinco}).


Then
$$\sin\left(z_1(\zeta,t)-z_1(w,t)\right)=\frac{\sin\left(z_1(\zeta,t)-z_1(w,t)\right)}{z_1(\zeta,t)-z_1(w,t)}\left(z_1(\zeta,t)-z_1(w,t)\right)$$
$$=G_1(\zeta,w,t)\cdot(\zeta-w),$$
with the $C^2-$norm of $G_1(\zeta,w,t)$ at most $C$ for fixed $t$.
Examining the asymptotics as $w\to \zeta$, yields
$$G_1(\zeta,\zeta,t)=\pa_\zeta z_1(\zeta,t),$$
hence
\begin{equation}\label{catorceocho}
\sin\left(z_1(\zeta,t)-z_1(w,t)\right)=\left[\pa_\zeta z_1(\zeta,t)+G_2(\zeta,w,t)(\zeta-w)\right](\zeta-w),\end{equation}
with the $C^1-$norm of $G_2(\zeta,w,t)$ at most $C$ for fixed $t$.
Similarly,
$$\cosh\left(z_2(\zeta,t)-z_2(w,t)\right)-1=\frac{\cosh\left(z_2(\zeta,t)-z_2(w,t)\right)-1}{\left(z_2(\zeta,t)-z_2(w,t)\right)^2}\left(z_2(\zeta,t)-z_2(w,t)\right)^2$$
$$=G_3(\zeta,w,t)(\zeta-w)^2,$$
with the $C^2-$norm of $G_3(\zeta,w,t)$ at most $C$ for fixed $t$.
Examining the asymptotics as $w\to \zeta$, gives
$$G_3(\zeta,\zeta,t)=\frac{1}{2}\left(\pa_\zeta z_2(\zeta,t)\right)^2,$$
hence
$$\cosh\left(z_2(\zeta,t)-z_2(w,t)\right)-1=\left[\frac{1}{2}\left(\pa_\zeta z_2(\zeta,t)\right)^2+G_4(\zeta,w,t)(\zeta-w)\right](\zeta-w)^2,$$
with the $C^1-$norm of $G_4(\zeta,w,t)$ at most $C$ for fixed $t$.
Similarly,
$$1-\cos\left(z_1(\zeta,t)-z_1(w,t)\right)=\left[\frac{1}{2}\left(\pa_\zeta z_1(\zeta,t)\right)^2+G_5(\zeta,w,t)(\zeta-w)\right](\zeta-w)^2,$$
with the $C^1-$norm of $G_5(\zeta,w,t)$ at most $C$ for fixed $t$.
Adding, and recalling the COMPLEX CHORD-ARC CONDITION (\ref{catorcecinco}), we see that

$$\cosh\left(z_2(\zeta,t)-z_2(w,t)\right)-\cos\left(z_1(\zeta,t)-z_1(w,t)\right)$$
\begin{equation}\label{catorcenueve}=\frac{1}{2}\left[\left(\pa_\zeta z_2(\zeta,t)\right)^2+\left(\pa_\zeta z_1(\zeta,t)\right)^2\right]\cdot\left[1+G_6(\zeta,w,t)\cdot(\zeta-w)\right](\zeta-w)^2,\end{equation}
with the $C^1-$norm of $G_6(\zeta,w,t)$ at most $C$ for fixed $t$.
Dividing (\ref{catorceocho}) by (\ref{catorcenueve}), and noting that
$$\frac{1}{1+G_6\cdot(\zeta-w)}=1-\frac{G_6}{1+G_6\cdot(\zeta-w)}\cdot (\zeta-w),$$
we obtain
$$\frac{\sin\left(z_1(\zeta,t)-z_1(w,t)\right)}{\cosh\left(z_2(\zeta,t)-z_2(w,t)\right)-\cos\left(z_1(\zeta,t)-z_1(w,t)\right)}$$
$$=\left[\frac{2\pa_\zeta z_1(\zeta,t)}{\left(\pa_\zeta z_2(\zeta,t)\right)^2+\left(\pa_\zeta z_1(\zeta,t)\right)^2}+G_7(\zeta,w,t)\cdot(\zeta-w)\right]\frac{1}{(\zeta-w)}$$
$$=\frac{2\pa_\zeta z_1(\zeta,t)}{\left(\pa_\zeta z_2(\zeta,t)\right)^2+\left(\pa_\zeta z_1(\zeta,t)\right)^2}\frac{1}{(\zeta-w)}+G_7(\zeta,w,t),$$
with the $C^1-$norm of $G_7(\zeta,w,t)$ at most $C$ for fixed $t$..

Since also
$$\frac{2\pa_\zeta z_1(\zeta,t)}{\left(\pa_\zeta z_2(\zeta,t)\right)^2+\left(\pa_\zeta z_1(\zeta,t)\right)^2}\cot\left(\frac{\zeta-w}{2}\right)$$
$$=\frac{2\pa_\zeta z_1(\zeta,t)}{\left(\pa_\zeta z_2(\zeta,t)\right)^2+\left(\pa_\zeta z_1(\zeta,t)\right)^2}\frac{1}{\zeta-w}+G_8(\zeta,w,t),$$
with the $C^{1}$ norm of $G_8(\zeta,w,t)$ on $\underline{\Omega}$ bounded by $C$, we conclude that
$\underline{a}(\zeta,w,t)$ belongs to $C^{1}(\underline{\Omega})$, with norm at most $C$.

Again recalling that there is no problem outside $\underline{\Omega}$ (i.e., when $|\zeta-w|>c$ even after $\zeta$, $w$ are translated by multiples of $2\pi$), we conclude \eqref{catorcediez}.


Next, we notice that
\begin{equation}\label{quincedos}
\tilde{a}(\zeta,t)=\int_{\Im w =h(\Re w,t)}a(\zeta,w,t)dw=\int_{w\in \T}a(\zeta,w,t)dw\end{equation}
if $t\in[\tau^2,\tau]$;
\begin{equation}\label{quincetres}
\tilde{a}(\zeta,t)=\int_{\Im w =\hbar(\Re w,t)}a(\zeta,w,t)dw=\int_{w\in \T}a(\zeta,w,t)dw\end{equation}
if $t\in[-\tau^2,\tau^2].$

We shall study the regularity of the function $\tilde{a}(\zeta)$ for
$\zeta\in\Gamma_+$ and $$\Gamma_+=\{\zeta\in \C\,:\,\, \zeta=x+ih(x),\,\,\,x\in\T\}.$$ We recall that $h(x)$ is a smooth function and we will prove
that  $\tilde{a}$ has  $C^2$-norm bounded by a controlled constant. (We replace $h$ by $\hbar$ if we are in the interval $[-\tau^2,\tau^2]$.)

\begin{rem}
Notice that since $z(\zeta,t)\in C^3(\Omega(t)^\text{closure})$, from \eqref{quincedos} and \eqref{quincetres} we see that $\tilde{a}$ is  a $C^2$-function.
However, this argument is not enough to  control  the norm
$$||\tilde{a}(\cdot+i h(\cdot,t))||_{C^2(\T)}$$
by a controlled constant. We would obtain a constant which depends on $\frac{1}{\min_{x\in\T} h(x,t)}$ and we can not permit it.
\end{rem}

In fact, we will prove that $\tilde{a}$ has bounded norm in $C^{2+\delta}(\Omega(t))$ for any $\delta<1/2$. To see this, it is enough to show that
\begin{equation}\label{paraz}
\tilde{a}|_{\Gamma_+(t)}\in C^{2+\delta}\quad\text{with norm at most $C(\delta)$, for any $\delta<\frac{1}{2}$}.
\end{equation}
(The reduction from $\Omega(t)$ to $\Gamma_+(t)$ is possible, thanks to $X_1$ in  section \ref{seccionsobo}).

To prove \eqref{paraz}, we first note that
\begin{equation}\label{4.63}
\tilde{a}(\zeta,t)=P.V. \int_{\Gamma_+(t)}\frac{\sin(z_1(\zeta,t)-z_1(w,t))}{\cosh(z_2(\zeta,t)-z_2(w,t))-\cos(z_1(\zeta,t)-z_1(w,t))}d\zeta
\end{equation}
for $\zeta\in \Gamma_+(t)$, since
\begin{equation*}
P.V.\int_{\Gamma_+(t)}\cot\left(\frac{\zeta-w}{2}\right)dw=0,
\end{equation*}
as we noted in section \ref{seccion1}.

It is useful to make the change of variable\begin{align*}\zeta=&x+ih(x,t),\\
w=&x+y+i h(x+y,t),
\end{align*}
and to set
\begin{equation}\label{4.64}
\hat{z}_\mu(x)\equiv z_\mu(x+ih(x,t),t)
\end{equation}
for $\mu=1$, 2.

Then \begin{equation}\text{$\hat{z}_1(x)-x$ and $\hat{z}_2(x)$ belong to $C^{3.5}(\T)$, with norm less or equal than $C$}.\label{4.65a}\end{equation}.

Hence also
\begin{equation}
\frac{\hat{z}_\mu(x)-\hat{z}_\mu(x+y)}{y}=-\int_0^1\hat{z}'_\mu(x+sy)ds\label{4.65b}
\end{equation}
belongs to $C^{2.5}(\T\times [-\pi,\,\pi])$ with norm less or equal than $C$.

Our change of variable \eqref{4.64} gives
\begin{equation}
\tilde{a}(x+ih(x,t),t)=P.V.\int_{-\pi}^\pi \frac{\sin(\hat{z}_1(x)-\hat{z}_1(x+y))\left(1+i\pa_x h(x+y,t)\right)}{\cosh(\hat{z}_2(x)-\hat{z}_2(x+y))-\cos(\hat{z}_1(x)-\hat{z}_1(x+y))}dy
\label{4.66}\end{equation}
for $x\in\T$.

It is convenient to introduce the analytic functions
\begin{eqnarray}\label{4.67}
g_{sn}(\zeta)=\frac{\sin(\zeta)}{\zeta}, & g_{csh}(\zeta)=\frac{\cosh(\zeta)-1}{\zeta^2}, & g_{cs}(\zeta)=\frac{1-\cos(\zeta)}{\zeta^2}\end{eqnarray}
for $\zeta\in\C$ and
\begin{eqnarray}
g^\sharp_{sn}(\zeta,\underline{\zeta})=\frac{g_{sn}(\zeta)-g_{sn}(\underline{\zeta})}{\zeta-\underline{\zeta}},\nonumber\\
g^\sharp_{csh}(\zeta,\underline{\zeta})=\frac{g_{csh}(\zeta)-g_{csh}(\underline{\zeta})}{\zeta-\underline{\zeta}},\nonumber\\
g^\sharp_{cs}(\zeta,\underline{\zeta})=\frac{g_{cs}(\zeta)-g_{cs}(\underline{\zeta})}{\zeta-\underline{\zeta}}, \label{4.68}
\end{eqnarray}
for $\zeta$, $\underline{\zeta}\in\C$.

We can then write
\begin{equation}\label{4.69}
y^{-1}\sin(\hat{z}_1(x)-\hat{z}_1(x+y))=g_{sn}(\hat{z}_1(x)-\hat{z}_1(x+y))\left[\frac{\hat{z}_1(x)-\hat{z}_1(x+y)}{y}\right]
\end{equation}
and
\begin{align}
&y^{-2}\left[\cosh(\hat{z}_2(x)-\hat{z}_2(x+y))-\cos(\hat{z}_1(x)-\hat{z}_1(x+y))\right]\nonumber\\
&g_{csh}(\hat{z}_2(x)-\hat{z}_2(x+y))\left[\frac{\hat{z}_2(x)-\hat{z}_2(x+y)}{y}\right]^2+
g_{cs}(\hat{z}_1(x)-\hat{z}_1(x+y))\left[\frac{\hat{z}_1(x)-\hat{z}_1(x+y)}{y}\right]^2\label{4.70}
\end{align}
Thanks to \eqref{4.65a} and \eqref{4.65b}, the right-hand sides of \eqref{4.69} and \eqref{4.70} belong to
$C^{2.5}(\T\times[-\pi,\pi])$, with norm at most $C$. Also, the left-hand side of \eqref{4.70} has absolute value larger or equal than $c$, thanks to the chord-arc condition.

Therefore, we may express the integrand in \eqref{4.66} in the form $\frac{Q(x,y)}{y}$, where
\begin{equation*}
Q(x,y)=\frac{\left[\text{LHS of \eqref{4.69}}\right](1+i\pa_xh(x+y,t))}{\left[\text{LHS of \eqref{4.70}}\right]}
\end{equation*}
belongs to $C^{2.5}(\T\times[-\pi,\,\pi])$, with norm at most $C$.
That is,
\begin{equation}\label{4.71}
\tilde{a}(x+ih(x,t),t)=P.V.\int_{-\pi}^\pi \frac{Q(x,y)}{y}dy=\int_{-\pi}^\pi \frac{Q(x,y)-Q(x,0)}{y}dy.
\end{equation}
with
\begin{equation*}
||Q||_{C^{2.5}(\T\times[-\pi,\,\pi])}\leq C.
\end{equation*}
It  follows easily that the function $x\mapsto \tilde{a}(x+ih(x,t),t)$ belongs to $C^{2+\delta}(\T)$ with norm
less or equal than $C(\delta)$ for each $\delta<\frac{1}{2}$.
This completes the proof of \eqref{paraz}.

To conclude this section, we prove that
\begin{align}
&\frac{\sin(z_1(\zeta,t)-z_1(w,t))}{\cosh(z_2(\zeta,t)-z_2(w,t))-\cos(z_1(\zeta,t)-z_1(w,t))}-
\frac{\sin(\un_1(\zeta,t)-\un_1(w,t))}{\cosh(\un_2(\zeta,t)-\un_2(w,t))-\cos(\un_1(\zeta,t)-\un_1(w,t))}\nonumber\\
&= \lambda \left[ \Theta(\zeta,t)\cot\left(\frac{\zeta-w}{2}\right)+\Theta^\sharp(\zeta,w,t)\right],\label{quincesiete}
\end{align}
for $\zeta$, $w\in\Gamma_+(t)$, where $\Theta\in C^0(\Gamma_+(t))$ and $\Theta^\sharp\in C^0(\Gamma_+(t)\times\Gamma_+(t))$ for fixed $t$, with norms less or equal than $C$.

To see this, we again use the change of variable
\begin{align*}
\zeta=&x+ih(x,t),\\
w=& x+y+i h(x+y,t);
\end{align*}
and we introduce $\hat{z}_\mu(x)$, $\hat{\un}_\mu(x)$ for $\mu=1,$ 2, where $\hat{z}_\mu(x)$ is given by \eqref{4.64}
and $\hat{\un}_\mu(x)$ is the analogous expression with $\un_\mu$ in place of $z_\mu$.

From \eqref{4.68} and \eqref{4.69} and the analogue of \eqref{4.69} for the $\hat{\un}_\mu$, we see that
\begin{align}
&y^{-1}\left[\sin(\hat{z}_1(x)-\hat{z}_1(x+y))-\sin(\hat{\un}_1(x)-\hat{\un}_1(x+y))\right]\nonumber\\
&=\left[g_{sn}(\hat{z}_1(x)-\hat{z}_1(x+y))-g_{sn}(\hat{\un}_1(x)-\hat{\un}_1(x+y))\right]
\left[\frac{(\hat{z}_1(x)-\hat{z}_1(x+y))}{y}\right]\nonumber\\
&+\left[g_{sn}(\hat{\un}_1(x)-\hat{\un}_1(x+y))\right]\left[\frac{(\hat{z}_1(x)-\hat{z}_1(x+y))}
{y}-\frac{(\hat{\un}_1(x)-\hat{\un}_1(x+y))}{y}\right]\nonumber\\
&=\left[g^\sharp_{sn}(\hat{z}_1(x)-\hat{z}_1(x+y),\,\hat{\un}_1(x)-\hat{\un}_1(x+y))\right]
\left[(\hat{z}_1-\hat{\un}_1)(x)-(\hat{z}_1-\hat{\un}_1)(x+y)\right]\nonumber\\
&\times \left[\frac{\hat{z}_1(x)-\hat{z}_1(x+y)}{y}\right]\nonumber\\
&+\left[g_{sn}(\hat{z}_1(x)-\hat{z}_1(x+y))\right]
\left[\frac{(\hat{z}_1-\hat{\un}_1)(x)-(\hat{z}_1-\hat{\un}_1)(x+y)}{y}\right].\label{4.73}
\end{align}
Since $\hat{z}_\mu-\hat{\un}_\mu$ belongs to $C^{3.5}$ with norm less or equal than $C\lambda$, we
know that $\frac{\hat{z}_\mu(x)-\hat{\un}_\mu(x+y)}{y}$ belongs to $C^{2.5}$ with norm less or equal than $C\lambda$,
as in the proof of \eqref{4.65b}. Therefore, \eqref{4.73} shows that
\begin{align}
&y^{-1}[\sin(\hat{z}_1(x)-\hat{z}_1(x+y))-\sin(\hat{\un}_1(x)-\hat{\un}_1(x+y))]\nonumber\\
&\text{belongs to $C^{2.5}(\T\times[-\pi,\,\pi])$, with norm less or equal $C\lambda$.}\label{4.74}
\end{align}
A similar argument involving the function $g_{csh}^\sharp$ and $g_{cs}^\sharp$ in \eqref{4.68} shows that
\begin{align}
&y^{-2}\left[\left\{\cosh(\hat{z}_2(x)-\hat{z}_2(x+y))-\cos(\hat{z}_1(x)-\hat{z}_1(x+y))\right\}-
\right.\nonumber\\
&\left.\left\{\cosh(\hat{\un}_2(x)-\hat{\un}_2(x+y))-\cos(\hat{\un}_1(x)-\hat{\un}_1(x+y))\right\}\right]\nonumber\\
&\text{belongs to $C^{2.5}(\T\times[-\pi,\,\pi])$, with norm less or equal than $C\lambda$}.\label{4.75}
\end{align}
We have seen that
\begin{align}
& y^{-1}\sin(\hat{z}_1(x)-\hat{z}_1(x+y))\quad \text{and}
\quad y^{-2}[\cosh(\hat{z}_2(x)-\hat{z}_2(x+y))-\cos(\hat{z}_1(x)-\hat{z}_1(x+y))]\nonumber\\
& \text{belong to $C^{2.5}(\T\times[-\pi,\,\pi])$ with norms less or equal than $C$, and that}\nonumber\\
& \left|y^{-2}\{\cosh(\hat{z}_2(x)-\hat{z}_2(x+y))-\cos(\hat{z}_1(x)-\hat{z}_1(x+y))\}\right|\geq c \quad \text{on $\T\times[-\pi,\, \pi]$}.\label{4.76}
\end{align}
From \eqref{4.74}, \eqref{4.75}, \eqref{4.76} it now follows easily that
\begin{align}
&\frac{\sin(\hat{z}_1(x)-\hat{z}_1(x+y))}{\cosh(\hat{z}_2(x)-\hat{z}_2(x+y))-\cos(\hat{z}_1(x)-\hat{z}_1(x+y))}-\nonumber\\
&\frac{\sin(\hat{\un}_1(x)-\hat{\un}_1(x+y))}{\cosh(\hat{\un}_2(x)-\hat{\un}_2(x+y))-\cos(\hat{\un}_1(x)-\hat{\un}_1(x+y))}\nonumber\\
&=\lambda \frac{Q^\sharp(x,y)}{y}\label{4.77}
\end{align}
for a function $Q^\sharp\in C^{2.5}(\T\times[-\pi,\,\pi])$ with norm less or equal than $C$.

Our result \eqref{4.77} is much sharper than our assertion \eqref{quincesiete} In particular, we have proven \eqref{quincesiete}.

\subsection{Perturbing the Dangerous Term}\label{seccion16}
Let $z$, $\un$ be as in the section (\ref{seccion13}). Suppose that $\un$ satisfies the additional assumptions made in the section (\ref{seccion14}).

Let $a(\zeta,w,t)$, $\underline{a}(\zeta,w,t)$, $\tilde{a}(\zeta,t)$, $\underline{\tilde{a}}(\zeta,t)$ be as in section (\ref{seccion14}).

We recall from  section (\ref{seccion12}) that
$$Dangerous_{\mu}(\zeta,t)=$$
\begin{equation}\label{dieciseisuno}
\int\limits_{w\in\Gamma_+(t)}\frac{\sin\left(z_1(\zeta,t)-z_1(w,t)\right)}{\cosh\left(z_2(\zeta,t)-z_2(w,t)\right)-\cos\left(z_1(\zeta,t)-z_1(w,t)\right)}[\pa_\zeta^{4+1}z_\mu(\zeta,t)-\pa_w^{4+1}z_\mu(w,t)]dw\end{equation}
for $\zeta\in \Gamma_+(t)$.

For the Muskat solution $\un$ in place of $z$, the analogous expression is

$$\underline{Dangerous}_{\mu}(\zeta,t)=$$
\begin{equation}\label{dieciseisdos}
\int\limits_{w\in\Gamma_+(t)}\frac{\sin\left(\un_1(\zeta,t)-\un_1(w,t)\right)}{\cosh\left(\un_2(\zeta,t)-\un_2(w,t)\right)-\cos\left(\un_1(\zeta,t)-\un_1(w,t)\right)}[\pa_\zeta^{4+1}\un_\mu(\zeta,t)-\pa_w^{4+1}\un_\mu(w,t)]dw\end{equation}
for $\zeta\in \Gamma_+(t)$.

Our goal in this section is to compute $Dangerous_\mu(\zeta,t)-\underline{Dangerous}_\mu(\zeta,t)$,
modulo small errors. Let $\zeta\in\Gamma_+(t)$ be given. From (\ref{dieciseisuno}) and (\ref{dieciseisdos}), we have
$$Dangerous_\mu(\zeta,t)-\underline{Dangerous}_\mu(\zeta,t)=$$
$$\int\limits_{w\in\Gamma_+(t)}\frac{[\sin\left(z_1(\zeta,t)-z_1(w,t)\right)][\pa_\zeta^{4+1}(z_\mu-\un_\mu)(\zeta,t)-\pa_w^{4+1}(z_\mu-\un_\mu)(w,t)]}{\cosh\left(z_2(\zeta,t)-z_2(w,t)\right)-\cos\left(z_1(\zeta,t)-z_1(w,t)\right)}dw$$
$$+\int\limits_{w\in\Gamma_+(t)}\left\{\frac{\sin\left(z_1(\zeta,t)-z_1(w,t)\right)}{\cosh\left(z_2(\zeta,t)-z_2(w,t)\right)-\cos\left(z_1(\zeta,t)-z_1(w,t)\right)}\right.$$$$\left.-\frac{\sin\left(\un_1(\zeta,t)-\un_1(w,t)\right)}{\cosh\left(\un_2(\zeta,t)-\un_2(w,t)\right)-\cos\left(\un_1(\zeta,t)-\un_1(w,t)\right)}\right\}$$
$$\times[\pa_\zeta^{4+1}\un_\mu(\zeta,t)-\pa_w^{4+1}\un_\mu(w,t)]dw.$$
Substituting the definition of $a(\zeta,w,t)$ and applying (\ref{quincesiete}), we conclude that
$$Dangerous_\mu(\zeta,t)-\underline{Dangerous}_\mu(\zeta,t)=$$
$$\frac{\pa_\zeta z_1(\zeta,t)}{\left(\pa_\zeta z_1(\zeta,t)\right)^2+\left(\pa_\zeta z_2(\zeta,t)\right)^2}\int\limits_{w\in \Gamma_+(t)}
\cot\left(\frac{\zeta-w}{2}\right)\left[ \pa_\zeta^{4+1}(z_\mu-\un_\mu)(\zeta,t)-\pa_w^{4+1}(z_\mu-\un_\mu)(w,t)\right]dw$$
$$+\int\limits_{w\in \Gamma_+(t)}a(\zeta,w,t)\left[ \pa_\zeta^{4+1}(z_\mu-\un_\mu)(\zeta,t)-\pa_w^{4+1}(z_\mu-\un_\mu)(w,t)\right]dw$$
\begin{equation}\label{dieciseistres}+\int\limits_{w\in \Gamma_+(t)}\lambda \left\{\Theta(\zeta,t)\cot\left(\frac{\zeta-w}{2}\right)+\Theta^\sharp(\zeta,w,t)\right\}[\pa_\zeta^{4+1}\un_\mu(\zeta,t)-\pa_w^{4+1}\un_\mu(w,t)]dw.\end{equation}
Recall that we define
\begin{equation}\label{dieciseiscuatro}
\Lambda_{\Gamma_+(t)}F(\zeta)=-\frac{1}{2\pi}\int\limits_{w\in\Gamma_+(t)}\cot\left(\frac{\zeta-w}{2}\right)[F'(\zeta)-F'(w)]dw
\end{equation}
for holomorphic function $F$.
Also, in the terms involving $a(\zeta,w,t)$ in (\ref{dieciseistres}), we recall (\ref{quincedos}), and we integrate by parts in $w$.
Thus, (\ref{dieciseistres}) becomes
$$Dangerous_\mu(\zeta,t)-\underline{Dangerous}_\mu(\zeta,t)=$$
$$-2\pi\frac{\pa_\zeta z_1(\zeta,t)}{\left(\pa_\zeta z_1(\zeta,t)\right)^2+\left(\pa_\zeta z_2(\zeta,t)\right)^2}\left[\Lambda_{\Gamma_+(t)}(\pa_\mu^{4}z_\mu(\cdot,t)-\pa_\mu^{4}\un_\mu(\cdot,t))\right]$$
$$+\tilde{a}(\zeta,t)\pa_\zeta^{4+1}(z_\mu(\zeta,t)-\un_\mu(\zeta,t))$$
$$+\int\limits_{w\in\Gamma_+(t)}[\pa_w a(\zeta,w,t)]\pa_w^{4}(z_\mu-\un_\mu)(w,t)dw$$
\begin{equation}\label{dieciseiscinco}+\lambda \int\limits_{w\in\Gamma_+(t)}\left\{\Theta(\zeta,t)\cot\left(\frac{\zeta-w}{2}\right)+\Theta^\sharp(\zeta,w,t)\right\}[\pa_\zeta^{4+1}\un_\mu(\zeta,t)-\pa_w^{4+1}\un_\mu(w,t)]dw.\end{equation}

Recall from assumption 2 in  section (\ref{seccion13}) that
$$|\pa_\zeta^{4+1}\un_\mu(\zeta,t)-\pa_\zeta^{4+1}\un_\mu(w,t)|\leq \min\{C|\zeta-w|,C\}.$$
(Here, it is important that we refer to $\un$, not to $z$.)

Since also
$$|\Theta(\zeta,t)|,\, |\Theta^\sharp(\zeta,w,t)|\leq C$$
by (\ref{quincesiete}),
we have
\begin{equation}\label{dieciseisseis}
\left|\lambda\int\limits_{w\in \Gamma_+(t)}\left\{\Theta(\zeta,t)\cot\left(\frac{\zeta-w}{2}\right)+\Theta^\sharp(\zeta,w,t)\right\}\cdot
[\pa_\zeta^{4+1}\un_\mu(\zeta,t)-\pa_\zeta^{4+1}\un_\mu(w,t)]dw\right|\leq C\lambda,
\end{equation}
for any $\zeta\in \Gamma_+(t)$.

Also, since
$$|\pa_w a(\zeta,w,t)|\leq C$$
by (\ref{catorcediez}) and assumption 4 in section \ref{seccion13}  (recall, we assume $\lambda<1$), we have
\begin{equation}\label{dieciseissiete}
\left|\int\limits_{w\in \Gamma_+(t)}[\pa_w a(\zeta,w,t)]\pa_w^{4}(z_\mu-\un_\mu)(w,t)dw\right|\leq C\lambda\end{equation}
for all $\zeta\in \Gamma_+(t)$, thanks to assumption 4 in section (\ref{seccion13}).

Substituting (\ref{dieciseisseis}) and (\ref{dieciseissiete}) into (\ref{dieciseiscinco}), we see that
$$Dangerous_\mu(\zeta,t)-\underline{Dangerous}_\mu(\zeta,t)=$$
$$-2\pi\frac{\pa_\zeta z_1(\zeta,t)}{\left(\pa_\zeta z_1(\zeta,t)\right)^2+\left(\pa_\zeta z_2(\zeta,t)\right)^2}\left[\Lambda_{\Gamma_+(t)}(\pa_\mu^{4}z_\mu(\cdot,t)-\pa_\mu^{4}\un_\mu(\cdot,t))\right]$$
$$+\tilde{a}(\zeta,t)\pa_\zeta^{4+1}(z_\mu(\zeta,t)-\un_\mu(\zeta,t))$$
\begin{equation}\label{dieciseisocho}+Error_{PDT}(\zeta,t)\end{equation}
for $\zeta\in \Gamma_+(t)$, with
\begin{equation}\label{dieciseisnueve}
|Error_{PDT}(\zeta,t)|\leq C\lambda.\end{equation}
The basic results of this section are (\ref{dieciseisocho}), (\ref{dieciseisnueve}).

\subsection{The Main Energy Estimate}\label{seccion17}

In this subsection, we suppose that $\un$ is an unperturbed solution
of Muskat problem, while $z$ is a perturbed solution, as in the
last several sections. Let $\sigma_1^0$, $\sigma_1$,
$\sigma_{\pm}$, $h$ and $\hbar$ as in the sections \ref{seccion7}, \ref{seccion8} and \ref{seccion9}.

Our goal is to estimate (from below) the quantity
$$\frac{d}{dt}\int_{\zeta\in\Gamma_+(t)}\left|\pa_\zeta^{4}(z_\mu-\un_\mu)(\zeta,t)\right|^2d\Re \zeta.$$

To do so, we first recall (I) from  lemma (\ref{lemma6}). Applying that result to both our Muskat solutions, we see that
$$\pa_t[\pa_\zeta^{4}(z_\mu-\un_\mu)(\zeta,t)]=$$
$$[Dangerous_\mu(\zeta,t)-\underline{Dangerous}_\mu(\zeta,t)]+\sum_{j=1}^6[Safe^\mu_j(\zeta,t)-\underline{Safe}^\mu_j(\zeta,t)]$$
$$+\sum_{1\leq\nu\leq\nu_{max}}c^\mu_\nu [Easy^\mu_\nu(\zeta,t)-\underline{Easy}^\mu_\nu(\zeta,t)]$$
for $\zeta\in\Gamma_+(t)$.

According to (\ref{trecediecinueve}), (\ref{trecetreintaidos}) and (\ref{dieciseisocho}), (\ref{dieciseisnueve}), this implies that
$$\pa_t[\pa_\zeta^{4}(z_\mu-\un_\mu)(\zeta,t)]=$$
$$-2\pi\frac{\pa_\zeta z_1(\zeta,t)}{\left(\pa_\zeta z_1(\zeta,t)\right)^2+\left(\pa_\zeta z_2(\zeta,t)\right)^2}\left[\Lambda_{\Gamma_+(t)}(\pa_\mu^{4}z_\mu(\cdot,t)-\pa_\mu^{4}\un_\mu(\cdot,t))\right]+\tilde{a}(\zeta,t)\pa_\zeta^{4+1}(z_\mu(\zeta,t)-\un_\mu(\zeta,t))$$
\begin{equation}\label{diecisieteuno}+Error(\zeta,t)\end{equation}
for $\zeta\in \Gamma_+(t)$, where
\begin{equation}\label{diecisietedos}
\int_{w\in\Gamma_+(t)}|Error(\zeta,t)|^2d \Re \zeta \leq C\lambda^2.
\end{equation}
Recall our assumption
\begin{equation}\label{diecisietetres}
\int_{\zeta\in \Gamma_{+}(t)}\left|\pa_\zeta^{4}
(z_\mu-\underline{z}_\mu)(\zeta,t)\right|^2 d\Re\zeta\leq
\lambda^2
\end{equation}
and note that
$$\frac{1}{2}\frac{d}{dt}\int_{\zeta\in \Gamma_{+}(t)}\left|\pa_\zeta^{4} (z_\mu-\underline{z}_\mu)(\zeta,t)\right|^2 d\Re\zeta$$
$$=\Re\int\limits_{\zeta\in \Gamma_{+}(t)}\overline{[\pa_\zeta^{4}(z_\mu-\un_\mu)(\zeta,t)]}\cdot\left\{\pa_t[\pa_\zeta^{4}(z_\mu-\un_\mu)(\zeta,t)]+ih_t(\Re\zeta,t)\cdot \pa_\zeta^{4+1}(z_\mu-\un_\mu)(\zeta,t)\right\}d\Re \zeta,$$
where
$$\Gamma_+(t)=\{x+ih(x)\,:\, x\in\T\}$$
if $t\in[\tau^2,\tau]$ and
$$\Gamma_+(t)=\{x+i\hbar(x)\,:\, x\in\T\}$$
if $t\in[-\tau^2,\tau^2]$. (In this range we use $\hbar$ in place of $h$
in the previous expression.)

We conclude that
$$\frac{1}{2}\frac{d}{dt}\int_{\zeta\in \Gamma_{+}(t)}\left|\pa_\zeta^{4} (z_\mu-\underline{z}_\mu)(\zeta,t)\right|^2 d\Re\zeta$$
$$=\Re\int\limits_{\zeta\in\Gamma_+(t)}\overline{F(\zeta)}\cdot\left\{\frac{-2\pi z'_1(\zeta,t)}{(z'_1(\zeta,t))^2+(z'_2(\zeta,t))^2}\la_{\Gamma_+(t)}F(\zeta)+\left(i\pa_th(\Re \zeta,t)+\tilde{a}(\zeta,t)\right)F'(\zeta)\right\}d\Re\zeta$$
\begin{equation}\label{diecisietecuatro} +Error ,\end{equation}
where
\begin{equation}\label{diecisietecinco}
F(\zeta)=\pa_\zeta^{4}(z_\mu-\un_\mu)(\zeta,t),\end{equation} and
\begin{equation}\label{diecisieteseis}
|Error|\leq C\lambda^2.
\end{equation}

Recall (\ref{catorcequince}),...,(\ref{catorcedieciocho}), (\ref{catorceonce}) and (\ref{paraz}). Our assumptions on $A$, $\lambda$, $\tau$, $\kappa$ in section \ref{seccion8} imply
\begin{equation}\label{diecisietesiete}
0< \lambda <\tau^2.
\end{equation}
Recall from section \ref{seccion14}, equations \eqref{catorceonce}, \eqref{catorcedieciseis}, \eqref{catorcediecisiete} and \eqref{catorcedieciocho}, that
\begin{align}
&\tilde{a}(\zeta,t)\quad\text{is a $C^{2}-$ function of $\zeta\in\Gamma_+(t)$, with norm at most $C$;}\label{diecisieteocho}\\
&|\Im \tilde{a}(\zeta,t)|\leq C h(\Re \zeta,t)\quad \text{for $\zeta\in \Gamma_+(t)$};\label{diecisietenueve}\\
&|\pa_xh(x,t)\tilde{a}(\zeta,t)|\leq Ch(x,t)+CA^{-1}\tau^2\quad \text{for $\zeta=x+ih(x,t)\in\Gamma_+(t)$.}\label{diecisietediez}
\end{align}
Estimates (\ref{diecisietenueve}), (\ref{diecisietediez}) hold for $t\in[\tau^2,\tau]$. When $t\in[-\tau^2,\tau^2]$, we use $\hbar$ in place of $h$.

Note that
$$\Im\{\tilde{a}(x+ih(x,t),t)\cdot(1+i\pa_xh(x,t))^{-1}\}$$
$$=[\Im \tilde{a}(x+ih(x,t),t)-\pa_xh(x,t)\cdot \Re \tilde{a}(x+ih(x,t),t)]\cdot\left(1+(\pa_xh(x,t))^2\right)^{-1}.$$
Hence (\ref{diecisieteocho}),...,(\ref{diecisietediez}) show that
\begin{equation}\label{diecisieteonce}
\tilde{a}(x+ih(x,t),t)\cdot(1+i\pa_xh(x,t))^{-1}=a_{R}(x)+ia_{I}(x),
\end{equation}
where
\begin{equation}\label{diecisietedoce}
a_{R}(x),\,a_{I}(x)\quad\text{are real, $C^{2}-$functions on
$\T$, with norm $\leq C$}\end{equation} and
\begin{equation}\label{diecisietetrece}
|a_{I}(x)|\leq \frac{Ch(x,t)+CA^{-1}\tau^2}{1+(\pa_x
h(x,t))^2}\quad \text{for all $x\in\T$}.
\end{equation}
If $t\in[-\tau^2,\tau^2]$, then we use $\hbar$ in place of $h$.

Let
\begin{equation}\label{diecisietecatorce}
f(x)=F(x+ih(x,t))=\pa_\zeta^{4}(z_\mu-\un_\mu)(\zeta,t)\Big|_{\zeta=x+ih(x,t)}.
\end{equation}
Thus, (\ref{diecisietetres}) gives
\begin{equation}\label{diecisietequince}
\int_{x\in\T}|f(x)|^2dx\leq C\lambda^2.
\end{equation}
Note that
\begin{equation}\label{diecisietedieciseis}
f'(x)=(1+i\pa_xh(x,t))F'(x+ih(x,t)).
\end{equation}
We have
\begin{align}
&\Re \int_{\zeta\in\Gamma_+(t)}\overline{F(\xi)}\tilde{a}(\xi)F'(\xi)d\Re \xi\nonumber\\
&=\Re \int_{\zeta\in\Gamma_+(t)}\overline{f}(x)\cdot\left\{[a_{R}(x)+ia_{I}(x)]\cdot(1+i\pa_x h(x,t))\right\}\left\{(1+i\pa_xh(x,t))^{-1}f'(x)\right\}dx\nonumber\\
&=\Re\int_{x\in\T}\overline{f}(x)\cdot[a_{R}(x)+ia_{I}(x)]f'(x)dx\nonumber\\
&=\int_{x\in\T}a_{R}(x)\cdot\frac{1}{2}\frac{d}{dx}\left(|f(x)|^2\right)dx
+\Re\int_{x\in\T}\overline{f(x)}a_{I}(x)if'(x)dx\nonumber\\
\label{diecisietediecisiete}&=-\int_{x\in\T}a'_{R}(x)\cdot\frac{1}{2}\left(|f(x)|^2\right)dx
+\Re\int_{x\in\T}\overline{f(x)}a_{I}(x)if'(x)dx.\end{align}
The first term on the right is dominated by $C\lambda^2$ in absolute value, thanks to (\ref{diecisietedoce}), (\ref{diecisietequince}). By (\ref{diecisietedoce}), (\ref{diecisietetrece}) and the sharp G{\aa}rding inequality, the second term on the right is bounded below by
$$-\Re\int\limits_{x\in\T}\overline{f(x)}\left[Ch(x,t)+CA^{-1}\tau^2\right](1+(\pa_x h(x,t))^2)^{-1}\la f(x)dx-C\int\limits_{x\in\T}|f(x)|^2.$$
Putting the above remarks and (\ref{diecisietequince}) into (\ref{diecisietediecisiete}), we see that
$$\Re \int_{\zeta\in\Gamma_+(t)}\overline{F}(\zeta)\tilde{a}(\zeta,t)F'(\zeta)d\zeta$$
\begin{equation}\label{diecisietedieciocho}
\geq -C\lambda^2-C\Re \int_{x\in\T}\overline{f}(x)\left[h(x,t)+A^{-1}\tau^2\right](1+(\pa_xh(x,t))^2)^{-1}\la f(x)dx.\end{equation}
It is convenient to use (\ref{six}) to conclude that
$$\left| \Re \int_{x\in\T}\overline{f}(x)\left[h(x,t)+A^{-1}\tau^2\right]\frac{i\pa_xh(x,t)}{1+(\pa_xh(x,t))^2}\la f(x) dx\right|$$
$$\leq C\int_{x\in\T}|f(x)|^2dx\leq C'\lambda^2.$$
Together with (\ref{diecisietedieciocho}), this yields the estimate
$$\Re \int_{\zeta\in\Gamma_+(t)}\overline{F}(\zeta)\tilde{a}(\zeta,t)F'(\zeta)d\Re\zeta$$
$$\geq -C\lambda^2-C\Re\int_{x\in\T}\overline{f}(x)\left[h(x,t)+A^{-1}\tau^2\right]\frac{(1-i\pa_xh(x,t))}{1+(\pa_xh(x,t))^2}\la f(x) dx.$$
Recalling that
$$\la f(x)=(1+i\pa_xh(x,t))\la_{\Gamma_+(t)}F(x+ih(x,t))+Error(x),$$
with
$$||Error(x)||_{L^2(\T)}\leq C||f||_{L^2(\T)}$$
(see (\ref{dossiete})), we conclude that
$$\Re \int_{\zeta\in\Gamma_+(t)}\overline{F}(\zeta)\tilde{a}(\zeta,t)F'(\zeta)d\Re\zeta$$
$$\geq -C\lambda^2-C\Re \int_{x\in\T}\overline{f}(x)\left[h(x,t)+A^{-1}\tau^2\right]\frac{(1-i\pa_xh(x,t))(1+i\pa_xh(x,t))}{1+(\pa_xh(x,t))^2}\la_{\Gamma_+(t)} F(x+ih(x,t)) dx,$$
where again we use (\ref{diecisietequince}) to control the term involving $Error(x)$. Recalling the definition (\ref{diecisietecatorce}) of $f$, we can conclude that
$$\Re \int_{\zeta\in\Gamma_+(t)}\overline{F}(\zeta)\tilde{a}(\zeta,t)F'(\zeta)d\Re\zeta$$
$$\geq- C\lambda^2-C\Re \int_{\zeta\in\Gamma_+(t)}\overline{F(\zeta)}\left[h(x,t)+A^{-1}\tau^2\right]\la_{\Gamma_+(t)}F(\zeta)d\Re\zeta.$$
Substituting this inequality into (\ref{diecisietecuatro}), we see that
$$\frac{1}{2}\frac{d}{dt}\int_{\zeta\in \Gamma_{+}(t)}\left|\pa_\zeta^{4} (z_\mu-\underline{z}_\mu)(\zeta,t)\right|^2 d\Re\zeta$$
$$\geq \Re\int\limits_{\zeta\in\Gamma_+(t)}\overline{F(\zeta)}\cdot\left\{\left[\sigma_+(\Re\zeta,t)-Ch(x,t)-CA^{-1}\tau^2\right]\la_{\Gamma_+(t)}F(\zeta)\right.$$
\begin{equation}\label{diecisietediecinueve}+\left(i\pa_th(\Re \zeta,t)\right)F'(\zeta)\Big\}d\Re\zeta-C\lambda^2,\end{equation}
where we recall that
$$\sigma_+(\Re(\zeta),t)=\frac{-2\pi\pa_\zeta z_1(\zeta,t)}{\left(\pa_\zeta z_1(\zeta,t)\right)^2+\left(\pa_\zeta z_2(\zeta,t)\right)^2}\quad \text{for $\zeta\in\Gamma_+(t)$}.$$

\subsubsection{The main energy estimate I}
Now suppose $t\in[\tau^2,\tau]$. We prove a lower bound for the right-hand side of (\ref{diecisietediecinueve}).

We set up the partition of the unity $1=\theta_{in}+\theta_{out}$ on $\T$ with $\theta_{in}$, $\theta_{out}\geq 0$,
\begin{align*}
supp\, \theta_{in}  \subset & \{ ||x||\leq 20A^{-1}t^{\frac{1}{2}}\},\\
supp\, \theta_{out}  \subset & \{ ||x||\geq 10A^{-1}t^{\frac{1}{2}}\},
\end{align*}
and
\begin{equation*}
\left|\left(\frac{d}{dx}\right)^j \theta_{in,\,out}(x)\right|\leq C \left(A^{-1}t^{\frac{1}{2}}\right)^{-j}
\end{equation*}
for $0\leq j\leq 3$.

Of course $\theta_{in}$ and $\theta_{out}$ depend on $t$, but we have fixed $t$ above.

Our plan is to insert a factor $\theta_{in}$ or $\theta_{out}$ in the integrand in  (\ref{diecisietediecinueve}) and study the two resulting integrals.

First we study the term involving $\theta_{out}$. We prepare to apply lemma \ref{lemma2}. Here, we take $h(x)=h(x,t)$, $h_t(x)=\pa_t h(x,t)$, $\theta=\theta_{out}$, $\tilde{\sigma}_1=\sigma_+-Ch-CA^{-1}\tau^2$ and $S=6A^2$. We check the hypotheses of lemma \ref{lemma2}, with constants of the form $C(A)$. We write $C(A)$ to denote constants depending only on $A$.

Hypothesis (1) merely says that $h(x,t)>0$ for $t\in[\tau^2,\tau]$, which is immediate from \eqref{ochouno}.

Hypothesis (2) asserts that
$$|\pa_t h(x,t)|\leq 6A^2 h(x,t)$$
on $supp\, \theta_{out}$. This assertion follows at once from \eqref{ochocinco}.

Hypothesis (3) asserts that the $C^3-$norm of $h(x,t)$ and $\pa_t h(x,t)$ as functions of $x$ for fixed $t$, are dominated by $C(A)$.  This assertion follows at once from \eqref{ochoocho} and \eqref{ochonueve}.

Hypothesis (4) asserts that the $C^2$-norm of the function $\tilde{\sigma}_1$ is bounded by $C(A)$. We know that the $C^2-$norm of $Ch(x,t)+CA^{-1}\tau^2$ (for fixed $t$) is bounded by $C A^{-1}$ (see \eqref{ochoocho}). Also the $C^2$-norm of $\sigma_{+}$ is at most $C$ by \eqref{nuevetrece}. Therefore, hypothesis 4 holds.

Hypothesis (5) asserts that the $C^2-$norm of the function
$$x\mapsto \theta_{out}(x)(1+ih'(x,t))^{-1}(\tilde{\sigma}_1(x)+\pa_th(x,t))$$
is at most $C(A)$.

Since the $C^3-$norm of $h(x,t)$ (for fixed $t$) is at most $CA^{-1}$,  it is enough to show that
\begin{equation}\label{diecisieteveintiuno}
x\mapsto \theta_{out}(x)(\tilde{\sigma}_1(x)+\pa_th(x,t))\end{equation}
has $C^2-$ norm at most $C(A)$.

We estimate the derivatives of this function up to $2^{\underline{nd}}$ order.

For $||x||\leq 10A^{-1}t^\frac{1}{2}$, we have $\theta_{out}(x)=0$. Hence, the function in question vanishes in that interval.

For $||x||\geq 20A^{-1}t^\frac{1}{2}$, we have $\theta_{out}(x)$=1.

Since we have already checked that the $C^2-$norm of $\tilde{\sigma}_1(x)$ is at most $C$, and since the $C^2-$norm of $\pa_t h$ (for fixed $t$) is at most $C(A)$, we see that the derivatives up to order 2 of the function in (\ref{diecisieteveintiuno}) are at most $C(A)$ for $||x||\geq 20A^{-1}t^\frac{1}{2}.$

It remains to examine the region $\{10A^{-1}t^\frac{1}{2}\leq ||x||\leq 20A^{-1}t^\frac{1}{2}$\}.

First we recall that
\begin{equation}\label{diecisieteveintidos}
\left|\pa_x^j \sigma_+(x)\right|\leq CA(t^\frac{1}{2})^{2-j}\end{equation}
for $0\leq j\leq 2$, $||x||\leq 20A^{-1}t^\frac{1}{2}$, by \eqref{nueveveinte},
and
\begin{equation}\label{diecisieteveintitres}
\left|\pa_x^j\left\{\pa_th(x,t)\right\}\right|\leq CA(t^\frac{1}{2})^{2-j}\end{equation}
for $0\leq j\leq 2$, $||x||\leq 20A^{-1}t^\frac{1}{2}$, by \eqref{ochodiecinueve},
and also
\begin{equation}\label{diecisieteveinticuatro}
\left|\pa_x^j\theta_{out}(x)\right|\leq C(A^{-1}t^\frac{1}{2})^{-j}\end{equation}
for $0\leq j\leq 2$, $||x||\leq 20A^{-1}t^\frac{1}{2}$.

From (\ref{ochodieciseis}), we have
\begin{equation}\label{diecisieteveinticinco}
\left|\pa_x^jh(x,t)\right|\leq C(t^\frac{1}{2})^{2-j}\end{equation}
for $0\leq j\leq 2$, $||x||\leq 20A^{-1}t^\frac{1}{2}$.
And also, obviously
\begin{equation}\label{diecisieteveintiseis}
\left|\pa_x^jA^{-1}\tau^2\right|\leq C(t^\frac{1}{2})^{2-j}\end{equation}
for $0\leq j\leq 2$, $||x||\leq 20A^{-1}t^\frac{1}{2}$, since $t\in[\tau^2,\tau]$.

From (\ref{diecisieteveintidos}), (\ref{diecisieteveinticinco}), (\ref{diecisieteveintiseis}), we see that
$$\left|\pa_x^j\left\{ \sigma_+(x)-Ch(x,t)-CA^{-1}\tau^2\right\}\right|\leq CA(t^\frac{1}{2})^{2-j}$$
for $0\leq j\leq 2$, $||x||\leq 20A^{-1}t^\frac{1}{2}$,
i.e.,
$$\left|\pa_x^j\tilde{\sigma}_1(x)\right|\leq CA(t^\frac{1}{2})^{2-j}$$
for $0\leq j\leq 2$, $||x||\leq 20A^{-1}t^\frac{1}{2}$.

Together with (\ref{diecisieteveintitres}), this yields
$$\left|\pa_x^j\left\{\tilde{\sigma}_1(x)+\pa_th(x,t)\right\}\right|\leq CA(t^\frac{1}{2})^{2-j}$$
for $0\leq j\leq 2$, $||x||\leq 20A^{-1}t^\frac{1}{2}$.

In turn, this estimate and (\ref{diecisieteveinticuatro}) yield
$$\left|\pa_x^j\left[\theta_{out}(x)\left\{\tilde{\sigma}_1(x)+\pa_th(x,t)\right\}\right]\right|\leq CA^3(t^\frac{1}{2})^{2-j}$$
for $0\leq j\leq 2$, $||x||\leq 20A^{-1}t^\frac{1}{2}$, completing our verification of hypothesis 5.

Hypothesis 6 asserts that
$$\Re\left\{ \theta_{out}(x)(1+ih'(x,t))^{-1}(\tilde{\sigma}_1(x,t)+\pa_t h(x,t))\right\}\geq0$$
on $\T$, i.e.,
$$\frac{\theta_{out}(x)}{1+(h'(x,t))^2} \left\{\Re\tilde{\sigma}_1(x,t)+\pa_t h(x,t)+\pa_x h(x,t)\Im \tilde{\sigma}_1(x,t)\right\}\geq 0$$
on $\T$, which will follow if we can prove that
\begin{equation}\label{diecisieteveintisiete}
\Re\tilde{\sigma}_1(x,t)+\pa_t h(x,t)+\pa_x h(x,t)\Im \tilde{\sigma}_1(x,t)\geq 0
\end{equation}
on $\T$. By (\ref{nueveveintiuno}), (\ref{nueveveintitres}), we have
\begin{equation}\label{diecisieteveintisietea}
\Re \tilde{\sigma}_1(x,t)\geq \sigma_1^0(x,t)-Ch(x,t)-CA^{-1}\tau^2\end{equation}
for $x\in\T$, and
\begin{equation}\label{diecisieteveintisieteb}
|\Im \tilde{\sigma}_1(x,t)|\leq Ch(x,t)+CA^{-1}\tau^2\end{equation}
for $x\in\T$. We recall that
$$\sigma_1^0(x)=\frac{-2\pi \un_1'(x)}{(\un_1'(x))^2+(\un_2'(x))^2}.$$
Also,
\begin{equation}\label{diecisieteveintisietec}
|\pa_x h(x,t)|\leq CA^{-1}\end{equation}
by (\ref{ochoocho}).

The above remarks tell us that
$$\Re \tilde{\sigma}_1(x,t)+\pa_t h(x,t)+\pa_xh(x,t)\Im \tilde{\sigma}_1(x,t)$$
\begin{equation}\label{diecisieteveintiocho}
\geq \sigma_1^0(x,t)+\pa_t h(x,t)-Ch(x,t)-CA^{-1}\tau^2\end{equation}
for all $x\in\T$.

The right-hand side of (\ref{diecisieteveintiocho}) is positive for all $x\in\T$; see (\ref{ochoseis}). Therefore, (\ref{diecisieteveintisiete}) holds, completing the verification of hypothesis 6.

Note also that $0\leq\theta_{out}\leq 1$, hence  hypothesis 7 is satisfied.

Thus, hypotheses 1,...,7 hold here, with constants $C(A)$. Applying  lemma (\ref{lemma2}), we obtain the following inequality.
$$\Re\int\limits_{\zeta\in\Gamma_+(t)}\overline{F(\zeta)}\theta_{out}(\Re\zeta)\cdot\left\{\left[\frac{-2\pi z'_1(\zeta,t)}{(z'_1(\zeta,t))^2+(z'_2(\zeta,t))^2}-Ch(\Re \zeta,t)-CA^{-1}\tau^2\right]\la_{\Gamma_+(t)}F(\zeta)\right.$$
$$+\left(i\pa_th(\Re \zeta,t)\right)F'(\zeta)\Big\}d\Re\zeta$$
\begin{equation}\label{diecisieteveintinueve}
\geq -C(A)\sum_{\pm}\int_{\zeta\in\Gamma_{\pm}(t)}|F(\zeta)|^2d\Re\zeta.\end{equation}

Next, we prepare to apply lemma (\ref{lemma3}), taking $\theta=\theta_{in}$, $h(x)=h(x,t)$, $h_t=\pa_th(x,t)$, $\tilde{\sigma}_1=\sigma_+ -Ch-CA^{-1}\tau^2$  and $\delta=20A^{-1}t^\frac{1}{2}$.
We now check that the hypotheses of lemma (\ref{lemma3}) hold here, with constants $C(A)$.

Hypothesis (1) simply asserts that $h>0$. Since $t\in [\tau^2,\tau]$, this follows at once from \eqref{ochouno}.

Hypothesis (2) asserts that (for fixed $t$) the $C^2-$norm of $\tilde{\sigma}_1$ and the $C^3-$norms $h(x,t)$ and $\pa_t h(x,t)$ are bounded by $C(A)$. This is contained in our previous results \eqref{nuevetrece}, \eqref{ochoocho} and \eqref{ochonueve}.

Hypothesis 3 asserts that the $C^2-$ norms of $\theta_{in}(x)\tilde{\sigma}_1(x)$ and $\theta_{in}\pa_t h(x,t)$ are at most $C(A)$.

Recall that $supp\, \theta_{in}\subset\{||x||<20A^{-1}t^{\frac{1}{2}}\}$, where we have the estimates
\begin{align*}
&\left|\left(\frac{d}{dx}\right)^j\theta_{in}(x)\right|\leq C\left(A^{-1}t^\frac{1}{2}\right)^{-j}\quad \text{for $0\leq j\leq 2$},\\
&\left|\left(\frac{d}{dx}\right)^j\pa_t h(x,t)\right|\leq C(A)\left(A^{-1}t^\frac{1}{2}\right)^{2-j}\quad \text{for $0\leq j\leq 2$ (see \eqref{ochodiecinueve})},\\
&\left|\left(\frac{d}{dx}\right)^j \sigma_+(x,t)\right|\leq C(A)\left(A^{-1}t^\frac{1}{2}\right)^{2-j}\quad \text{for $0\leq j\leq 2$ (see \eqref{nueveveinte})},\\
&\left|\left(\frac{d}{dx}\right)^jh(x,t)\right|\leq C\left(t^\frac{1}{2}\right)^{2-j}\quad \text{for $0\leq j\leq 2$ (see \eqref{ochodieciseis})},\\
&\left|\left(\frac{d}{dx}\right)^jA^{-1}\tau^2\right|\leq C\left(t^\frac{1}{2}\right)^{2-j}\quad \text{for $0\leq j\leq 2$, $t\in [\tau^2,\tau$)}.
\end{align*}
Hence,
\begin{align}
&\left|\left(\frac{d}{dx}\right)^j\theta_{in}(x)\pa_t h(x,t)\right|\leq C(A)\left(t^\frac{1}{2}\right)^{2-j}\quad \text{for $0\leq j\leq 2$} ,\label{manolita}\\
&\left|\left(\frac{d}{dx}\right)^j\theta_{in}(x)\tilde{\sigma}_1(x)\right|\leq C(A)\left(t^\frac{1}{2}\right)^{2-j}\quad \text{for $0\leq j\leq 2$.}\label{manolito}
\end{align}
In particular, the $C^2-$norm of $\theta_{in}\tilde{\sigma}_1(x)$ and $\theta_{in}(x)\pa_t h(x,t)$ are at most $C(A)$, completing the verification of hypothesis 3.

Hypothesis 4 is satisfied since $\supp\,\,\theta_{in}\subset\{||x||\leq 20 A^{-1}t^\frac{1}{2}\}$.

Hypothesis 5 asserts that
\begin{equation}\label{diecisietetreintaiseis}
\left|\left(\frac{d}{dx}\right)^j\Re\left\{\theta_{in}(x)(1+i\pa_xh(x,t))^{-1}\tilde{\sigma}_1(x)\right\}\right|\leq C(A)\delta^{1-j}\end{equation}
and
\begin{equation}\label{diecisietetreintaisiete}
\left|\left(\frac{d}{dx}\right)^j\Re\left\{\theta_{in}(x)(1+i\pa_xh(x,t))^{-1}\pa_th(x,t)\right\}\right|\leq
C(A)\delta^{1-j}\end{equation} for $0\leq j\leq 2$, $x\in\T$.

From \eqref{ochoocho} we see that
$$\left|\left(\frac{d}{dx}\right)^j\left\{(1+i\pa_xh(x,t))^{-1}\right\}\right|\leq C\leq C\delta^{-j}$$
for $0\leq j\leq 2$, $x\in \T$. Together with \eqref{manolita} and \eqref{manolito}, this implies the estimates
$$\left|\left(\frac{d}{dx}\right)^j\left\{\theta_{in}(x)(1+i\pa_xh(x,t))^{-1}\tilde{\sigma}_1(x)\right\}\right|\leq C(A)\delta^{2-j}$$
and
$$\left|\left(\frac{d}{dx}\right)^j\Re\left\{\theta_{in}(x)(1+i\pa_xh(x,t))^{-1}\pa_th(x,t)\right\}\right|\leq
C(A)\delta^{2-j},$$
which are  stronger than  \eqref{diecisietetreintaiseis} and \eqref{diecisietetreintaisiete}.
This completes the verification of hypothesis 5.

Hypothesis 6 asserts that, for all $x\in\T$, we have
$$|\Re\left\{\theta_{in}(x)(1+i\pa_xh(x,t))^{-1}\pa_t h(x,t)\right\}|\leq \Re \left\{\theta_{in}(x)(1+i\pa_xh(x,t))^{-1}\tilde{\sigma}_1(x)\right\}$$
i.e.,
$$\frac{\theta_{in}(x)}{1+(\pa_xh(x,t))^2}\left|\Re \left\{(1-i\pa_xh(x,t))\pa_th(x,t)\right\}\right|\leq \frac{\theta_{in}(x)}{1+(\pa_xh(x,t))^2}\Re \left\{(1-i\pa_xh(x,t))\tilde{\sigma}_1(x)\right\},$$
i.e.,
\begin{equation}\label{diecisietetreintaiocho}
|\pa_t h(x,t)|\leq \Re \tilde{\sigma}_1(x)+\pa_xh(x,t)\Im \tilde{\sigma}_1(x)
\end{equation}
in $supp\,\theta_{in}$.

In view of (\ref{diecisieteveintisietea}), (\ref{diecisieteveintisieteb}), (\ref{diecisieteveintisietec}), estimate (\ref{diecisietetreintaiocho}) will follow, if we can prove that
\begin{equation}\label{diecisietetreintainueve}
|\pa_th(x,t)|\leq \sigma_1^0(x,t)-C'h(x,t)-C'A^{-1}\tau^2
\end{equation}
in $supp\, \theta_{in}$.

However, (\ref{diecisietetreintainueve}) follows at once from (\ref{ochosiete}). This completes the verification of hypothesis 6.

Note also that $0\leq \theta_{in}\leq 1$, hence  hypothesis 7 is satisfied.

Thus, all the hypotheses (1),...,(7) of the lemma (\ref{lemma3}) hold here, with constants $C(A)$. Applying that lemma, we obtain the following inequality.
$$\Re\int\limits_{\zeta\in\Gamma_+(t)}\overline{F(\zeta)}\theta_{in}(\Re\zeta)\cdot\left\{\left[\frac{-2\pi z'_1(\zeta,t)}{(z'_1(\zeta,t))^2+(z'_2(\zeta,t))^2}-Ch(\Re \zeta,t)-CA^{-1}\tau^2\right]\la_{\Gamma_+(t)}F(\zeta)\right.$$
$$+\left(i\pa_th(\Re \zeta,t)\right)F'(\zeta)\Big\}d\Re\zeta$$
$$\geq -C(A)\int_{\zeta\in\Gamma_{+}(t)}|F(\zeta)|^2d\Re\zeta.$$
Adding this to (\ref{diecisieteveintinueve}), we obtain our basic lower bound:
$$\Re\int\limits_{\zeta\in\Gamma_+(t)}\overline{F(\zeta)}\cdot\left\{\left[\frac{-2\pi z'_1(\zeta,t)}{(z'_1(\zeta,t))^2+(z'_2(\zeta,t))^2}-Ch(\Re \zeta,t)-CA^{-1}\tau^2\right]\la_{\Gamma_+(t)}F(\zeta)\right.$$
$$+\left(i\pa_th(\Re \zeta,t)\right)F'(\zeta)\Big\}d\Re\zeta$$
$$\geq -C(A)\sum_{\pm}\int_{\zeta\in\Gamma_{\pm}(t)}|F(\zeta)|^2d\Re\zeta.$$
This estimate holds for $t\in[\tau^2,\tau]$.

Together with (\ref{diecisietediecinueve}) and (\ref{diecisietecinco}), this gives the estimate
$$\frac{1}{2}\frac{d}{dt}\int_{\zeta\in \Gamma_{+}(t)}\left|\pa_\zeta^{4} (z_\mu-\underline{z}_\mu)(\zeta,t)\right|^2 d\Re\zeta$$
$$\geq -C\lambda^2-C(A)\sum_{\pm}\int_{\zeta\in\Gamma_{\pm}(t)}|\pa_\zeta^{4}(z_\mu-\un_\mu)(\zeta)|^2d\Re\zeta$$
for $t\in[\tau^2,\tau]$, under the assumptions made on $z$, $\un$.

Recalling from assumption 4 in  section \ref{seccion13} that the integral on the right is at most $\lambda^2$, we conclude that
\begin{equation}\label{diecisietecuarenta}
\frac{1}{2}\frac{d}{dt}\int_{\zeta\in
\Gamma_{+}(t)}\left|\pa_\zeta^{4}
(z_\mu-\underline{z}_\mu)(\zeta,t)\right|^2 d\Re\zeta\geq
-C(A)\lambda^2,\end{equation} if $t\in[\tau^2,\tau]$.

This estimate holds provided $z_\mu$ is a Muskat solutions satisfying
\begin{equation}\label{diecisietecuarentaiuno}
\sum_{\pm}\int_{\zeta\in\Gamma_{\pm}(t)}|(z_\mu-\un_\mu)(\zeta)|^2d\Re\zeta+\sum_{\pm}\int_{\zeta\in\Gamma_{\pm}(t)}|\pa_\zeta^{4}(z_\mu-\un_\mu)(\zeta)|^2d\Re\zeta\leq
\lambda^2
\end{equation}
and
\begin{equation}\label{diecisietecuarentaidos}
\lambda<\tau^2.
\end{equation}
Here, $$\Gamma_{\pm}(t)=\{x\pm ih(x,t)\, :\, x\in\T\}.$$

\subsubsection{The main energy estimate II}

Next, we fix $t\in[-\tau^2,\tau^2]$. We prepare to apply lemma (\ref{lemma2}) for $\theta(x)=1$, $h(x)=\hbar(x,t)$, $h_t(x)=\pa_t \hbar(x,t)$, $S=\tau^{-1}$, and
\begin{equation}\label{diecisietecuarentaitres}
\tilde{\sigma}_1(x)= \frac{-2\pi z'_1(\zeta,t)}{(z'_1(\zeta,t))^2+(z'_2(\zeta,t))^2}\Big|_{\zeta=x+i\hbar(x,t)}-C\hbar(x,t)-CA^{-1}\tau^2.\end{equation}

We will check that the hypotheses of the lemma (\ref{lemma2}) hold, with constants of the form $C(A)$.

Hypothesis 1 asserts that $\hbar(x,t)>0$ for all $x\in\T$, which is immediate from \eqref{ochoveintidos}.

Hypothesis 2 asserts that $|\pa_t\hbar(x,t)|\leq \tau^{-1}\hbar(x,t)$ for all $x\in\T$, which is precisely \eqref{ochoveintiseis}.

Hypothesis 3 asserts that, for fixed $t$, the $C^3-$norms of $\hbar(x,t)$ and $\pa_t\hbar(x,t)$ are bounded by $C(A)$. This is immediate from \eqref{ochoveintiocho} and \eqref{ochoveintinueve}.

Hypothesis 4 asserts that the $C^2-$norm of the function $\tilde{\sigma}_1$ is at most $C(A)$. The $C^2-$norm of the function $\sigma_+(x,t)$ is at most $C$; see (\ref{nuevetrece}). The $C^2-$norm of the function
$$x\mapsto C\hbar(x,t)+CA^{-1}\tau^2$$
is at most $CA^{-1}$; see (\ref{ochoveintiocho}).

These remarks and (\ref{diecisietecuarentaitres}) show that the $C^2-$norm of $\tilde{\sigma}_1$ is at most $C$, proving hypothesis 4.

Hypothesis 5 asserts that the $C^2-$norm of the function
\begin{equation}\label{diecisietecuarentaicuatro}
(1+i\pa_x\hbar(x,t))^{-1}(\tilde{\sigma}_1(x)+\pa_t\hbar(x,t))\end{equation}
is at most $C(A)$.
We have just seen that the $C^2-$norm of $\tilde{\sigma}_1$ is at most $C$. Moreover, the $C^2-$norm of $\pa_t\hbar(x,t)$ is at most $CA$, by (\ref{ochoveintinueve}). Since also the $C^2-$norm of $(1+i\pa_x\hbar(x,t))^{-1}$ is at most $C$, thanks to (\ref{ochoveintiocho}), we conclude that the function (\ref{diecisietecuarentaicuatro}) has $C^2-$norm at most $C(A)$, thus proving hypothesis 5.

Hypothesis 6 asserts that
$$\Re\{(1+i\pa_x\hbar(x,t))^{-1}(\tilde{\sigma}_1(x)+\pa_t\hbar(x,t))\}\geq 0$$
on $\T$, i.e.,
$$\Re \tilde{\sigma}_1(x)+\pa_x\hbar(x,t)\Im\tilde{\sigma}_1(x)+\pa_t \hbar (x,t)\geq 0$$
on $\T$, i.e.,
$$\Re \sigma_+(x,t)+\pa_t\hbar(x,t)-C\hbar(x,t)-CA^{-1}\tau^2$$
\begin{equation}\label{diecisietecuarentaicinco}+\pa_x\hbar(x,t)\Im \sigma_+(x,t)\geq 0\end{equation}
on $\T$.

From (\ref{nueveveintidos}) we have
$$\Re \sigma_+(x,t)\geq  \sigma_1^0(x,t)-C\hbar(x,t)-C\tau^{10}$$
for all $x\in\T$; and from (\ref{nueveveinticuatro}), we have
$$\left|\Im \sigma_{+}(x,t)\right|\leq C\hbar(x,t)+C\tau^{10}$$
for all $x\in \T$.

Also, from (\ref{ochoveintiocho}), we have
$$|\pa_x \hbar(x,t)|\leq CA^{-1}$$
for all $x\in\T$.

The above remarks imply that, for all $x\in\T$, we have
$$\Re \sigma_+(x,t)+\pa_t\hbar(x,t)-C\hbar(x,t)-CA^{-1}\tau^2$$
$$+\pa_x\hbar(x,t)\Im \sigma_+(x,t) $$
\begin{equation}\label{diecisietecuarentaiseis}
\geq \sigma_1^0(x,t)+\pa_t\hbar(x,t)-C\hbar(x,t)-CA^{-1}\tau^2.\end{equation}
The right-hand side of (\ref{diecisietecuarentaiseis}) is positive; see (\ref{ochoveintisiete}). Hence, (\ref{diecisietecuarentaiseis}) implies (\ref{diecisietecuarentaicinco}), completing the proof of hypothesis 6.

Thus, hypotheses 1,...,6 hold, with constants $C(A)$. Hypothesis 7 holds trivially, since $\theta=1.$

We may now apply the lemma in section (\ref{seccion4}). We obtain the estimate
$$\Re\int\limits_{\zeta\in\Gamma_+(t)}\overline{F(\zeta)}\cdot\left\{\left[\frac{-2\pi z'_1(\zeta,t)}{(z'_1(\zeta,t))^2+(z'_2(\zeta,t))^2}-C\hbar(\Re\zeta,t)-CA^{-1}\tau^2\right]\la_{\Gamma_+(t)}F(\zeta)\right.$$
$$+\left(i\pa_t\hbar(\Re \zeta,t)\right)F'(\zeta)\Big\}d\Re\zeta$$
\begin{equation}\label{diecisietecuarentaisiete}\geq -C(A)\tau^{-1}\sum_{\pm}\int_{\zeta\in\Gamma_{\pm}(t)}|F(\zeta)|^2d\Re\zeta.\end{equation}
This estimate holds for $t\in[-\tau^2,\tau^2]$.
Moreover, we have estimate (\ref{diecisietediecinueve}), in which $h$ should be replaced by $\hbar$ since we are working with $t\in[-\tau^2,\tau^2]$. Hence, (\ref{diecisietediecinueve}) and (\ref{diecisietecuarentaisiete}) yield the estimate
$$\frac{1}{2}\frac{d}{dt}\int_{\zeta\in \Gamma_{+}(t)}\left|\pa_\zeta^{4} (z_\mu-\underline{z}_\mu)(\zeta,t)\right|^2 d\Re\zeta$$
\begin{equation}\label{diecisietecuarentaiocho}\geq -C\lambda^2-C(A)\tau^{-1}\sum_{\pm}\int_{\zeta\in\Gamma_{\pm}(t)}|\pa_\zeta^{4}(z_\mu-\un_\mu)(\zeta)|^2d\Re\zeta\end{equation}
for $t\in[-\tau^2,\tau^2]$, with
$$\Gamma_{\pm}(t)=\{x\pm i\hbar(x,t)\, :\, x\in\T\}.$$
Recalling from (\ref{diecisietecinco}), and from assumption 4 in section (\ref{seccion13}) that
$$\sum_{\pm}\int_{\zeta\in\Gamma_{\pm}(t)}|\pa_\zeta^{4}(z_\mu-\un_\mu)(\zeta)|^2d\Re\zeta\leq \lambda^2,$$
we conclude from (\ref{diecisietecuarentaiocho}) that
\begin{equation}\label{diecisietecuarentainueve}
\frac{1}{2}\frac{d}{dt}\int_{\zeta\in
\Gamma_{+}(t)}\left|\pa_\zeta^{4}
(z_\mu-\underline{z}_\mu)(\zeta,t)\right|^2 d\Re\zeta\geq
-C(A)\tau^{-1}\lambda^2,\end{equation} if $t\in[-\tau^2,\tau^2]$.

This estimate holds provided $z_\mu$ is a Muskat solution satisfying
\begin{equation}\label{diecisietecincuenta}
\sum_{\pm}\int_{\zeta\in\Gamma_{\pm}(t)}|(z_\mu-\un_\mu)(\zeta)|^2d\Re\zeta+\sum_{\pm}\int_{\zeta\in\Gamma_{\pm}(t)}|\pa_\zeta^{4}(z_\mu-\un_\mu)(\zeta)|^2d\Re\zeta\leq
\lambda^2
\end{equation}
and
\begin{equation}\label{diecisietecincuentaiuno}
\text{ $A$, $\tau$, $\lambda$, $\kappa$ are as assumed in section \ref{seccion8}}.
\end{equation}
Here, $$\Gamma_{\pm}(t)=\{x\pm i\hbar(x,t)\,:\, x\in\T\}.$$
Our basic energy estimates are (\ref{diecisietecuarenta}),..., (\ref{diecisietecuarentaidos}) for $t\in[\tau^2,\tau]$ and (\ref{diecisietecuarentainueve}),...,(\ref{diecisietecincuentaiuno}) for $t\in[-\tau^2,\tau^2]$.

\subsection{Conclusion}\label{subseccion1}

The previous estimates allows us to achieve the following result.

\begin{thm}\label{mean}
Let $z(x,t)$ be a solution of the Muskat equation in the interval
$t\in[t_{\text{least}},\tau]\subset [-\tau^2,\tau]$ and let $\un(x,t)$ be the unperturbed solution.  Assume
that $z(x,t)$ satisfies
\begin{itemize}
\item  $z_1(x,t)-x$ and $z_2(x,t)$ are periodic with period
$2\pi$. \item $z(\zeta,t)$ is real for $\zeta$ real. \item
$z(\zeta,t)$ is analytic in $\zeta\in\Omega(t)$. \item
$z(\zeta,t)-(\zeta,0)\in H^4(\Omega(t))$.
 \item Complex Arc-Chord condition.
$$|\cosh(z_2(\zeta,t)-z_2(w,t))-\cos(z_1(\zeta,t)-z_1(w,t))|\geq c_{CA} [||\Re (\zeta-w)||+|\Im(\zeta-w)|]^2,$$
for $\zeta$, $w\in\overline{\Omega}(t)$.
\end{itemize}
Here in the definition of $\Omega(t)$ we use
$h(x,t)$ if $t\in[\tau^2,\tau]$ and $\hbar(x,t)$ if
$t\in[-\tau^2,\tau^2].$ Then
$$\frac{d}{dt}||z(\cdot,t)-\un(\cdot,t)||^2_{H^4(\Omega(t))}\geq - C(A)\lambda^2,$$

if $t\in [\tau^2,\tau]\cap [t_{\text{least}},\tau]$, and if also
$$||z(\cdot,t)-\un(\cdot,t)||_{H^4(\Omega(t))}\leq \lambda.$$

 In addition,
 $$\frac{d}{dt}||z(\cdot,t)-\un(\cdot,t)||^2_{H^4(\Omega(t))}\geq -C(A)\tau^{-1}\lambda^2,$$

if $t\in [-\tau^2,\tau^2]\cap [t_{\text{least}},\tau]$, and if also
$$||z(\cdot,t)-\un(\cdot,t)||_{H^4(\Omega(t))}\leq \lambda$$
and $\lambda\leq  \tau^2$.
\end{thm}
Proof: We have to estimate the quantity:
\begin{equation}\label{dnor}\frac{d}{dt}||z(\cdot,t)-\un(\cdot,t)||^2_{H^4(\Omega(t))},\end{equation}
where
\begin{equation}\label{nor}||z(\cdot,t)-\un(\cdot,t)||^2_{H^4(\Omega(t))}=\sum_{\pm}\int_{\zeta\in\Gamma_{\pm}(t)}|z(\zeta,t)-\un(\zeta,t)|^2d\Re\zeta+\sum_{\pm}\int_{\zeta\in\Gamma_{\pm}(t)}|\pa^4_\zeta (z(\zeta,t)-\un(\zeta,t))|^2d\Re\zeta.\end{equation}
The term in \eqref{dnor} coming from the second term on the right hand side on \eqref{nor} is bounded using \eqref{diecisietecuarenta} and \eqref{diecisietecuarentainueve}.
The rest of terms are lower order terms which are easy to control.

\section{Galerkin approximations and main theorem}

We  prove  existence of solution of the Muskat equation in $H^4(\Omega(t))$ and  the main theorem \ref{conclusion}.

\subsection{Conformal Maps Depending on a Parameter}\label{conformalmap}
Let $m\geq 100$ and $h_0>0$ be given, and let $h(x,t)$ be a
real-valued function on $\T\times[0,1]$ with
$$|\partial_t^l\partial_x^j h(x,t)|\leq 1$$ for $l,j\leq 100
m+100$.

We write $c$, $C$, etc. to denote constants determined by $m$
and $h_0$ alone.

Let $V_0(x+iy)$ be the solution of the Dirichlet problem
\begin{align*}
\Delta V_0(z)=0& \quad \text{for $-2h_0<\Im z<0$}\\
V_0(x)=0 & \quad \text{for $x\in\R$}\\
V_0(x-2ih_0)=\frac{1}{2\pi}\Im
\left\{\cot\left(\frac{x-2ih_0}{2}\right)\right\} & \quad
\text{for $x\in\R$}.
\end{align*}
Then $V_0$ is $2\pi-$periodic and real-analytic; hence we may
regard $V_0$ as harmonic, real analytic and $2\pi-$periodic on a
strip
$$H_0=\{-2h_0-c_0<\Im z<c_0\}.$$

We will solve the Dirichlet problem on the region
$$\Omega(t;\delta)=\{x+iy\,:\, -2h_0-\delta h(x,t)<y<\delta
h(x,t)\}$$ for $0<\delta<c_1$ (small enough $c_1$).

As an approximate Poisson kernel, we try
$$P(x+iy,x_+)=\frac{1}{2\pi}\Im
\left[\cot\left(\frac{\zeta}{2}\right)\right]\cdot (1+(\delta
h'(x_+,t))^2)-V_0(\zeta)\cdot(1+(\delta h'(x_+,t))^2),$$ for
$x+iy\in \Omega(t;\delta)$, $x_+\in\R$, where
$$\zeta=(x+iy)-(x_++i\delta
h(x_+,t))-i\delta h'(x_+,t)\sin[(x+iy)-(x_++i\delta h(x_+,t))].$$

Since $\delta<c_1$ for small enough $c_1$, we know that $\zeta\in
H_0$, so $P(x+iy,x_+)$ is well-defined.

Suppose $f\in C^m(\T)$. Then
$$u_t(x+iy)=\int_{x_+\in\T}P(x+iy,x_+)f(x_+)dx_+$$ is harmonic in
$\Omega(t,\delta)$, and one checks that
\begin{align}
&\lim_{\ep\to 0^+} u_t(x+i\delta h(x,t)-i\ep)=
f(x)+\int_{\T}K_1(x,x_+,t)dx_+\label{con1}\\
&\lim_{\ep\to 0^+} u_t(x+i\delta h(x,t)-2h_0i-i\ep)=\int_\T
K_2(x,x_+,t)dx_+\label{con2}\end{align}

where

\begin{equation}\label{con3}|\pa_t^l\pa_x^j\pa_{x_+}^k K_i(x,x_+,t)|\leq C\delta \quad \text{for
$l$, $j$, $k\leq m$ and $i=1$, $2$},\end{equation} and
\begin{equation}|\pa_x^j\pa_{x_+}^k K_i(x,x_+,t)|\leq C\delta \quad \text{for
 $j$, $k\leq m$ and $i=1$, $2$}.\label{con4}\end{equation}

We indicate the main steps in verifying the above.

We introduce a new variable $p$, and we study how $\zeta$ behaves when we take $y=\delta h(x,t)-p$. We work on the domain
\begin{equation*}
\{(x,x_+,p,\delta,t)\:\,|x-x_+|\leq \pi,\,\delta\in[0,c],\,t\in[0,1]\}
\end{equation*}
for a small constant $c$.

Let us write $G_1$, $G_2$, etc. to denote functions on this domain having at least $10m$ derivatives of absolute value at most $C$.

Then
\begin{equation}\label{confo1}
\zeta=Z-i\delta h'(x,t)\sin(Z)
\end{equation}
with
\begin{equation}\label{confo2}
Z=[x+i\delta h(x,t)-ip]-[x_++i\delta h(x_+,t)].
\end{equation}
First of all, note that
\begin{align*}
&\zeta=Z\cdot\left(1-i\delta h'(x_+,t)\frac{\sin(Z)}{Z}\right)\\
&=[(x-x_+)-ip]\cdot\left[1+i\delta\frac{h(x,t)-h(x_+,t)}{(x-x_+)+ip}\right]
\cdot\left[1-i\delta h'(x_+,t)\frac{\sin(Z)}{Z}\right],
\end{align*}
and therefore
\begin{equation}\label{confo3}
|\zeta|=|(x-x_+)-ip|\cdot (1+O(\delta))
\end{equation}
uniformly in $x$, $x_+$, $p$, $\delta$, $t$.

Next, by writing
\begin{equation*}
\sin(Z)=Z+Z^3\left[\frac{\sin(Z)-Z}{Z^3}\right],
\end{equation*}
we find that
\begin{equation*}
\zeta=Z-i\delta h'(x_+,t)Z+\delta Z^3G_1;
\end{equation*}
moreover,
\begin{equation*}
Z=(x-x_+)\cdot[1+i\delta h'(x_+,t)]-ip+\delta(x-x_+)^2G_2.
\end{equation*}
Hence,
\begin{equation}\label{confo4}
\zeta=[1+(\delta h'(x_+,t))^2]\cdot(x-x_+)-\delta ph'(x_+,t)-ip+\delta(x-x_+)^2G_3+\delta p^2G_4.
\end{equation}
Next, we refine \eqref{confo3} in the case $|x-x_+|<p^{\frac{2}{3}}$. In that case, we have
\begin{equation*}
|\delta(x-x_+)^2G_3+\delta p^2G_4|\leq C\delta p^\frac{4}{3}
\end{equation*}
and
\begin{equation*}
|\delta p h'(x_+,t)|\leq C\delta p,
\end{equation*}
while
\begin{equation*}
|[1+(\delta h'(x_+,t))^2](x-x_+)-ip|\geq cp.
\end{equation*}
Therefore, \eqref{confo4} implies that
\begin{equation}\label{confo5}
|\zeta|=|[1+(\delta h'(x_+,t))^2](x-x_+)-\delta p h'(x_+,t)-ip|\cdot (1+O(\delta p^\frac{1}{3})),
\end{equation}
uniformly in $x$, $x_+$, $p$, $\delta$, $t$, under the assumptions $|x-x_+|< p^\frac{2}{3}$.

Now we are ready to analyze $\Im\cot(\zeta/2)$. For $|\Re \zeta|\leq \frac{3}{2}\pi$, we have
\begin{equation*}
\frac{1}{2}\cot\left(\frac{\zeta}{2}\right)=\frac{1}{\zeta}+H(\zeta)
\end{equation*}
with $H$ analytic and real-valued for real $\zeta$.

On the other hand, $\zeta$ is a $C^{90m+90}-$smooth function of $x$, $x_+$, $p$, $\delta$, $t$, equal to $(x-x_+)$ when $p=\delta=0$. Therefore,
\begin{equation}\label{confo6}
\frac{1}{2}\cot\left(\frac{\zeta}{2}\right)=-\frac{\Im \zeta}{|\zeta|^2}+p G_5+\delta G_6.
\end{equation}

From \eqref{confo3},...,\eqref{confo6} and the properties of $V_0$ discussed above, it is now a routine task to check
\eqref{con1},...,\eqref{con4}.

Thus, $u_t(x+iy)$ is approximately equal to $f$ on the upper
boundary of $\Omega(t;\delta)$, approximately zero on the lower
boundary of $\Omega(t;\delta)$ and (exactly) harmonic in $\Omega(t,\delta)$.

Moreover, $K_i(x,x_+,t)$ ($i=1,2$) are given by explicit formulas.

There are analogous explicit formulas for a harmonic function on
$\Omega(t;\delta)$ approximately equal to a given function $g$ on
the lower boundary of $\Omega(t;\delta)$ and approximately equal
to zero on the upper boundary.

Our explicit formulas, together with a Neumann  series, show that
the Dirichlet problem

\begin{align*}
&V(z,t)\quad \text{harmonic in $\Omega(t;\delta)$ for each fixed $t$}\\
&V(x+i\delta h(x,t),t)=1\\
&V(x-i\delta h(x,t)-2ih_0,t)=-1
\end{align*}
admits a solution $V$ that is $C^m$ on the region
$$\Xi(\delta)=\{(z,t)\,:\,z\in \Omega(t;\delta)^{\text{closure}},\,
t\in[0,1]\}.$$

In particular, the above Neumann series converges, since we take
$\delta<c_1$ for a small enough $c_1$.

Taking a harmonic conjugate $U$ of $V$, we may suppose $U$ also
belongs to $C^m$ on the region $\Xi(\delta)$.

The analytic function $\Phi_t=U+iV$ is a conformal map of
$\Omega(t\,;\, \delta)$ to the strip $\{|\Im z|<1\}$.

Moreover, $\Phi_t$ depends $C^m-$smoothly on $t$, and
$\Phi_t(z+2\pi)=\Phi_t(z)+\lambda(t)$ (all $z\in \Omega(t;\delta)$
for a smooth, real-valued function $\lambda(t)$.

A trivial rescaling now shows that $\Omega(t\,;\,\delta)$ maps
conformally to a strip $\{|\Im z|\leq h^\sharp(t)\}$ by a
conformal map $\Phi_t^\sharp$, that depends $C^m-$smoothly on $t$
and satisfies $\Phi_t^\sharp(z+2\pi)=\Phi_t^\sharp(z)+2\pi$;
moreover, $h^\sharp(t)$ is $C^m-$smooth.

We can now easily remove our small $\delta$ assumptions.

Let
$$\Omega(t)=\{x+iy\,:\,h_-(x,t)<y<h_+(x,t)\}$$
for $t\in[0,1]$; where $h_+$ and $h_-$ are real, smooth, and
$2\pi$-periodic.

We assume that $h_-(x,t)<h_+(x,t)$ for all $x\in\R$, $t\in[0,1]$;
and we suppose that $h_-(x,t)$ and $h_+(x,t)$ are real-analytic in
$x$ for each fixed $t$. (The real-analyticity assumptions ought to
be removed, but our $h_+$ and $h_-$ happen to be real-analytic.
See section \ref{seccion8}.)

Given $t_0\in[0,1]$, a conformal mapping defined on a neighborhood
of $\Omega(t_0)^{\text{closure}}$ carries $\Omega(t_0)$ to a
strip. This same conformal transformation carries $\Omega(t)$ to a
domain close to that strip, whenever $t$ is close enough to $t_0$.
Rescaling in the $t-$ variable we reduce matters to the case
considered above. Hence, $\Omega(t)$ is mapped
conformally to $\{|\Im z|<h^\sharp(t)\}$ by a conformal map
$z\mapsto \Phi(z,t)$, that depends $C^m-$smoothly on $(z,t)$ for  $z\in\Omega(t)^{\text{closure}}$, $t$ near $t_0$. We have
also $\Phi(z+2\pi,t)=\Phi(z,t)+2\pi$.

Thus, we have our smooth family of conformal maps, defined in a
small neighborhood of any given $t_0$. These maps carry the upper
boundary (respectively, the lower boundary) of $\Omega(t)$  to the
upper boundary (respectively, the lower boundary) of $H(t)$. Since
such $\Phi_t$ are uniquely determined up to translations, it is
trivial to patch together our results on small $t$-intervals.

Thus, we obtain the following result.

\begin{lemma}\label{conformal}
Let $m\geq 100$ be given. For $t\in[0,1]$, let
$$\Omega(t)=\{x+iy\,:\, h_-(x,t)<y<h_+(x,t)\}$$
where $h_+$, $h_-$ are smooth, real-valued functions on
$\R\times[0,1]$. We suppose that $h_-(x,t)<h_+(x,t)$ for all
$(x,t)\in \R\times [0,1]$, and that $h_+(x,t)$, $h_-(x,t)$ are $2\pi-$periodic and
real-analytic in $x$ for each fixed $t$.

Then there exist a $C^m-$smooth, positive function $h^\sharp(t)$
defined on $[0,1]$, and a $C^m-$smooth map
$$\Phi\,:\, \{(z,t)\in \C\times [0,1]\,:\, z\in
\Omega^{\text{closure}}(t)\}\rightarrow \C,$$
such that, for each fixed $t\in[0,\,1]$, $z\mapsto\Phi(z,t)$ maps $\Omega(t)$ conformally to $\{|\Im \zeta|\leq h^\sharp(t)\}$.

 Moreover,
$\Phi(z+2\pi,t)=\Phi(z,t)+2\pi$ for each $(z,t)$ in the domain of
$\phi$.
\end{lemma}

\subsection{Changing Coordinates}\label{seccion18}

Let $h(x,t)$ be a real-analytic, positive function on $\T\times [t_0-\delta,t_0+\delta]$, for $h(x+2\pi,t)=h(x,t)$ for all $x$, $t$.
Let
\begin{eqnarray}
\Gamma_{\pm}(t)&=&\{\zeta\in\C\,:\, \Im \zeta=\pm h(\Re \zeta,t)\}\nonumber\\
\Omega(t)&=&\{\zeta\in\C\,:\, |\Im \zeta|<h(\Re \zeta,t)\}.\label{18.1}
\end{eqnarray}

Then, by lemma \ref{conformal} there exist a smooth positive function $\hh(\that)$ on $[t_0-\delta,t_0+\delta]$ and a smooth map
\begin{equation}\label{18.2}
\Phi\,:\,\{(\zeta,t)\,:\, \zeta\in \Omega^{closure}(t),\, |t-t_0|\leq \delta\}\mapsto \{(\hat{\zeta},\hat{t})\,:\, |\Im \hat{\zeta}|\leq \hat{h}(\hat{t}),\, |\hat{t}-t_0|\leq \delta\}
\end{equation} with the following properties:
\begin{equation}\label{18.3}
\Phi(\zeta,t)=(\hz,\that)=(\phi(\zeta,t),t)
\end{equation}
for each $(\zeta,t)$.
The map
\begin{equation}\label{18.4}
\zeta\mapsto \phi(\zeta,t)\end{equation}
is a biholomorphic map from $\Omega(t)$ to the strip $\{|\Im \hz|<\hh(\that)\}$, for each $t$.
\begin{equation}\label{18.5}
\phi(\zeta+2\pi,t)=\phi(\zeta,t)+2\pi\end{equation}
\begin{equation}\label{18.6}
\phi\,:\,\Gamma_+(t)\mapsto \{\Im \hz=\hh(t)\}\end{equation}
for each fixed $t$.
\begin{equation}\label{18.6a}
\phi(\zeta,t)\quad \text{is real for real $\zeta$}.\end{equation}
The inverse of $\Phi$ is the map
\begin{equation}\label{18.7}
\Psi(\hz,\that)=(\psi(\hz,\that),\that)=(\zeta,t).\end{equation}
We would like to solve the Muskat equation
\begin{equation}\label{18.8}
\pa_t z_\mu (\zeta,t)=\int\limits_{w\in\Gamma_+(t)}\frac{\sin(z_1(\zeta,t)-z_1(w,t))[\pa_\zeta z_\mu(\zeta,t)-\pa_wz_\mu(w,t)]}{\cosh(z_2(\zeta,t)-z_2(w,t))-\cos(z_1(\zeta,t)-z_1(w,t))}dw\end{equation}
for $\zeta\in\Omega(t)$, $t\in[t_0-\delta',t_0]$, with
\begin{equation}\label{18.9}
z_\mu(\zeta,t)\quad\text{holomorphic in $\zeta\in\Omega(t)$ for fixed $t$},\end{equation}
with
\begin{equation}\label{18.10}
z_\mu(\zeta,t)-\zeta \delta_{1,\,\mu}\in H^{4}\left(\Gamma_{\pm}(t)\right)\quad \text{for fixed $t$},\end{equation}
where $\delta_{i,\,j}$ is the Kronecker  delta, with
\begin{equation}\label{18.11}
\text{$z_1(\zeta,t)-\zeta$ and $z_2(\zeta,t)$ periodic with period $2\pi$, for fixed $t$.}\end{equation}
\begin{equation}\label{18.12}
\text{$z_\mu(\zeta,t)$ is real for real $\zeta$, any $t\in[t_0-\delta',t_0]$,}
\end{equation}
and satisfying the initial condition
\begin{equation}\label{18.13}
z_\mu(\zeta,t_0)=z^0_\mu(\zeta),\end{equation}
where $z^0_\mu(\cdot)$ satisfy (\ref{18.9}),...,(\ref{18.12}), as well as the chord-arc condition
$$|\cosh(z_2(\zeta,t_0)-z_2(w,t_0))-\cos(z_1(\zeta,t_0)-z_1(w,t_0))|$$
\begin{equation}\label{18.14}
\geq c_{CA}\left[||\Re(\zeta-w)||+|\Im (\zeta-w)|\right]^2\end{equation}
for $\zeta$, $w\in \Omega^{closure}(t_0)$.

We assume that $z^0_\mu(\zeta)$ satisfy a Rayleigh-Taylor condition, to be explained below. Our goal is to transform the Muskat problem (\ref{18.8}),...,(\ref{18.14}) to an equivalent problem in $(\hz,\that)$ coordinates, and to show that the Rayleigh-Taylor condition for (\ref{18.8}),...,(\ref{18.14}) implies an analogous Rayleigh-Taylor condition in $(\hz,\that)$ coordinates.

For (\ref{18.8}),...,(\ref{18.14}), we introduce the 'Generalized Rayleigh-Taylor function'
$$RT(\zeta,t)=$$
$$\Re\left(\frac{-2\pi\pa_\zeta z_1(\zeta,t)}{(\pa_\zeta z_1(\zeta,t))^2+(\pa_\zeta z_2(\zeta,t))^2}[1+i\pa_x h(\Re \zeta,t)]^{-1}\right)$$
$$+\Im\left(\left\{P.V.\int\limits_{w\in \Gamma_+(t)}\frac{\sin(z_1(\zeta,t)-z_1(w,t))}{\cosh(z_2(\zeta,t)-z_2(w,t))-\cos(z_1(\zeta,t)-z_1(w,t))})dw\right.\right.$$\begin{equation}\label{18.15}+i\pa_th(\Re \zeta,t)\Big\}[1+i\pa_x h(\Re \zeta,t)]^{-1}\Big)\end{equation}
for $\zeta\in \Gamma_+(t)$.

Our Rayleigh-Taylor assumption regarding (\ref{18.8}),...,(\ref{18.14}) is as follows
\begin{equation}\label{18.16}
RT(\zeta,t_0)>0\end{equation}
(strict inequality) for $\zeta\in \Gamma_+(t)$.
Let us check that assumption (\ref{18.16}) holds for $z_\mu(\zeta,t)$ as in the previous sections, using either of our previous functions $h(x,t)$ or $\hbar(x,t)$ for the function h of this section. (If we use $h(x,t)$, then we assume that $t\in[\tau^2,\tau]$; if we use $\hbar(x,t)$, then we assume that $t\in[-\tau^2,\tau^2]$). We will also assume that $\lambda\leq A^{-1}\tau^2$. After we check this, we return to the task of transforming (\ref{18.8}),...,(\ref{18.16}) into $(\hz,\that)-$coordinates.

Let's start with the case of the function $h(x,t)$ (rather than $\hbar(x,t)$).

Then, for $\zeta=x+ih(x,t)\in \Gamma_+(t)$, we have
$$\left|\frac{-2\pi\pa_\zeta z_1(\zeta,t)}{(\pa_\zeta z_1(\zeta,t))^2+(\pa_\zeta z_2(\zeta,t))^2}-\frac{-2\pi\pa_x z_1(x,t)}{(\pa_x z_1(x,t))^2+(\pa_x z_2(x,t))^2}\right|\leq Ch(x,t),$$
and
$$\left|\frac{-2\pi\pa_x z_1(x,t)}{(\pa_x z_1(x,t))^2+(\pa_x z_2(x,t))^2}-\frac{-2\pi\pa_x \un_1(x,t)}{(\pa_x \un_1(x,t))^2+(\pa_x \un_2(x,t))^2}\right|\leq CA^{-1}\tau^2.$$
Therefore,
$$\Re\left(\frac{-2\pi\pa_\zeta z_1(\zeta,t)}{(\pa_\zeta z_1(\zeta,t))^2+(\pa_\zeta z_2(\zeta,t))^2}[1+i\pa_x h(\Re \zeta,t)]^{-1}\right)$$
$$\geq -Ch(x,t)-CA^{-1}\tau^2+\Re\left(\frac{-2\pi\pa_x \un_1(x,t)}{(\pa_x \un_1(x,t))^2+(\pa_x \un_2(x,t))^2}[1+i\pa_x h(x,t)]^{-1}\right)$$
\begin{equation}\label{18.17}
=-Ch(x,t)-CA^{-1}\tau^2+\frac{-2\pi\pa_x \un_1(x,t)}{(\pa_x \un_1(x,t))^2+(\pa_x \un_2(x,t))^2}(1+(\pa_xh(x,t))^2)^{-1}.\end{equation}
Recall that
$$\tilde{a}(\zeta,t)\equiv \int\limits_{w\in \Gamma_+(t)}\left\{\frac{\sin(z_1(\zeta,t)-z_1(w,t))}{\cosh(z_2(\zeta,t)-z_2(w,t))-\cos(z_1(\zeta,t)-z_1(w,t))}\right.$$
\begin{equation}\label{18.19}
-\frac{\pa_\zeta \un_1(\zeta,t)}{(\pa_\zeta \un_1(\zeta,t))^2+(\pa_\zeta \un_2(\zeta,t))^2}\cot\left(\frac{\zeta-w}{2}\right)\Big\}dw.\end{equation}
Thanks to our observation in section (\ref{seccion3}) that
$$P.V.\int_{w\in\Gamma_+(t)}\cot\left(\frac{\zeta-w}{2}\right)dw=0,$$
we have that
\begin{equation}\label{18.18}
\tilde{a}(\zeta,t)\equiv P.V.\int\limits_{w\in \Gamma_+(t)}\frac{\sin(z_1(\zeta,t)-z_1(w,t))}{\cosh(z_2(\zeta,t)-z_2(w,t))-\cos(z_1(\zeta,t)-z_1(w,t))}dw\end{equation}

Comparing (\ref{18.19}) with (\ref{quinceuno}) and (\ref{quincedos}), we see that $\tilde{a}$ is as in section (\ref{seccion15}). Hence, (\ref{quincesiete}) together with (\ref{catorcequince}), (\ref{catorcedieciseis}), tell us that
$$|\pa_x h(x,t)\tilde{a}(\zeta,t)|\leq Ch(x,t)+CA^{-1}\tau^2$$
with $x=\Re \zeta$, $\zeta\in \Gamma_+(t)$.

Thanks to (\ref{catorcediecisiete}), (\ref{catorcedieciocho}) and (\ref{quincesiete}), analogous estimates hold also for $\hbar(x,t)$ in place of $h(x,t)$.

Consequently,
$$\Im \left\{ \tilde{a}(\zeta,t)[1+i\pa_x h(x,t)]^{-1}\right\}\geq -Ch(x,t)-CA^{-1}\tau^2.$$
Recalling (\ref{18.18}), we conclude that
$$
\Im\left\{P.V.\int\limits_{w\in \Gamma_+(t)}\frac{\sin(z_1(\zeta,t)-z_1(w,t))}{\cosh(z_2(\zeta,t)-z_2(w,t))-\cos(z_1(\zeta,t)-z_1(w,t))})dw [1+i\pa_x h(x,t)]^{-1}\right\}$$\begin{equation}\label{18.20}\geq -Ch(x,t)-CA^{-1}\tau^2\end{equation}
for $x=\Re\zeta$, $\zeta\in \Gamma_+(t)$.

Analogous estimates hold for $\hbar(x,t)$ in place of $h(x,t)$.

We note that
\begin{equation}\label{18.21}
\Im\left\{i\pa_t h(x,t)[1+i\pa_xh(x,t)]^{-1}\right\}=\frac{\pa_t h(x,t)}{1+(\pa_x h(x,t))^2}.\end{equation}
Putting (\ref{18.17}), (\ref{18.18}), (\ref{18.20}), (\ref{18.21}) into definition (\ref{18.15}), we see that, for $\zeta\in \Gamma_+(t)$, $x=\Re \zeta$, we have:
\begin{align}RT(\zeta,t) & \geq (1+(\pa_x h(x,t))^2)^{-1}\nonumber\\
&\times\left[ -Ch(x,t)-CA^{-1}\tau^2+\frac{-2\pi\pa_x \un_1(x,t)}{(\pa_x \un_1(x,t))^2+(\pa_x \un_2(x,t))^2}+\pa_t h(x,t)\right]\label{18.22}\end{align}
The term in square brackets in (\ref{18.22}) is strictly positive, thanks to (\ref{ochoseis}).

Thus, (\ref{18.16}) holds for $t\in[\tau^2,\tau]$, under the assumptions of the previous sections.

Similarly, (\ref{18.16}) holds also for $t\in[-\tau^2,\tau^2]$, with $\hbar(x,t)$ in place of $h(x,t)$, under the assumptions of the previous sections.

We have succeeded in verifying the Rayleigh-Taylor condition (\ref{18.16}) for our Muskat solutions from section \ref{seccionM}.
We now return to the business of transforming our Muskat problem (\ref{18.8}),...,(\ref{18.14}), into $(\hz,\that)-$coordinates.

We make the change of variable
\begin{equation}\label{18.23}
\text{$\hz=\phi(\zeta,t)$, $\hat{w}=\phi(w,t)$, for fixed $t$},\end{equation}
so that also
\begin{equation}\label{18.24}
\text{$\zeta=\psi(\hz,\that)$, $w=\psi(\hat{w},\that)$, for fixed $\that$}.\end{equation}
Recall, $\that=t$. We define
\begin{equation}\label{18.25}
\hze_\mu(\hz,\that)=z_\mu(\zeta,t)\quad \text{for $(\zeta,t)=\Psi(\hz,\that)$}.\end{equation}
Note
\begin{equation}\label{18.26}
\pa_{\hz} \hze_\mu(\hz,\that)=\pa_{\hz}\psi(\hz,t)\pa_{\zeta}z_\mu(\zeta,t)\quad\text{for $(\zeta,t)=\Psi(\hz,\that)$}.\end{equation}
Similarly,
\begin{equation}\label{18.27}
\pa_{\hw} \hze_\mu(\hw,\that)=\pa_{\hw}\psi(\hw,t)\pa_wz_\mu(w,t)\quad\text{for $(w,t)=\Psi(\hw,\that)$}.\end{equation}
Thus
\begin{equation}\label{18.28}
\pa_\zeta z_\mu(\zeta,t)=\pa_\zeta \phi(\zeta,t)\pa_{\hz}\hze_\mu(\hz,\that)\quad\text{for $(\hz,\that)=\Phi(\zeta,t)$};\end{equation}
and
\begin{equation}\label{18.29}
\pa_w z_\mu(w,t)=\pa_w \phi(w,t)\pa_{\hw}\hze_\mu(\hw,\that)\quad\text{for $(\hw,\that)=\Phi(\zeta,t)$}.\end{equation}
Also,
\begin{equation}\label{18.30}
dw=\pa_{\hw}\psi(\hw,\that)d\hw\end{equation}
for fixed $t$.

The chain rule shows that for $(\hz,\that)=\Phi(\zeta,t)$, we have
$$\pa_t z_\mu(\zeta,t)=\pa_t\left [\hze_\mu\circ\Phi\right](\zeta,t)=\pa_t\{\hze_\mu(\phi(\zeta,t),t)\}=\pa_{\that}\hze_\mu(\hz,\that)+[\pa_t\phi(\zeta,t)]\pa_{\hz}\hze_\mu(\hz,\that)$$
\begin{equation}\label{18.31}
=\pa_{\that}\hze_\mu(\hz,\that)+[(\pa_t\phi)\circ\Psi(\hz,\that)]\pa_{\hz}\hze_\mu(\hz,\that).\end{equation}
In view of the above remarks, our Muskat equation (\ref{18.8}) is equivalent to the following equation in $(\hz,\that)-$coordinates
$$\pa_{\that}\hze_\mu(\hz,\that)+[(\pa_t\phi)\circ\Psi(\hz,\that)]\pa_{\hz}\hze_\mu(\hz,\that)=$$
$$\int\limits_{\Im \hw=\hh(\that)}\frac{\sin(\hze_1(\hz,\that)-\hze_1(\hw,\that))}{\cosh(\hze_2(\hz,\that)-\hze_2(\hw,\that))-\cos(\hze_1(\hz,\that)-\hze_1(\hw,\that))}$$
\begin{equation}\label{18.32}
\times\left\{[(\pa_\zeta\phi)\circ\Psi(\hz,\that)]\pa_{\hz}\hze_\mu(\hz,\that)-[(\pa_w\phi)\circ\Psi(\hw,\that)]\pa_{\hw}\hze_\mu(\hw,\that)\right\}\pa_{\hw}\psi(\hw,\that)d\hw\end{equation}
Since $\phi(\cdot,t)$ and $\psi(\cdot,t)$ are inverse functions, we see that
$$[(\pa_w\phi)\circ\Psi(\hw,\that)]\cdot[\pa_{\hw}\psi(\hw,\that)]=1.$$
Hence, defining
\begin{equation}\label{18.33}
\hat{B}(\hz,\hw,\that)\equiv [(\pa_\zeta\phi)\circ\Psi(\hz,\that)][\pa_{\hw}\psi(\hw,\that)]-1,\end{equation}
we see that
\begin{equation}\label{18.34}
\text{$\hat{B}(\hz,\hw,\that)$ is smooth on $\{|\Im \hz|,\,|\Im \hw|\leq \hh(\that),\, \that\in [t_0-\delta,t_0+\delta]\}$}\end{equation}
\begin{equation}\label{18.35}
\text{$\hat{B}(\hz,\hw,\that)$ is holomorphic in $(\hz,\hw)\in\{|\Im \hz|,\,|\Im \hw|< \hh(\that)\}$ for fixed $\that$}\end{equation}
and
\begin{equation}\label{18.36}
\hat{B}(\hz,\hz,\that)=0.\end{equation}
\begin{rem}
Note that $\hat{B}$ does not depend on our Muskat solution; it is entirely determined by the function $h(x,t)$ and the coordinate change $(\hz,\that)=\Phi(\zeta,t)$.
\end{rem}
Thanks to (\ref{18.33}),...,(\ref{18.36}), we can write (\ref{18.32}) in the equivalent form
$$\pa_{\that}\hze_\mu(\hz,\that)=$$
$$\left\{-[(\pa_t\phi)\circ\Psi(\hz,\that)]+\int\limits_{\Im \hw=\hh(\that)}\frac{\sin(\hze_1(\hz,\that)-\hze_1(\hw,\that))\hat{B}(\hz,\hw,\that)}{\cosh(\hze_2(\hz,\that)-\hze_2(\hw,\that))-\cos(\hze_1(\hz,\that)-\hze_1(\hw,\that))}d\hw\right\}\pa_{\hz}\hze_\mu(\hz,\that)$$
\begin{equation}\label{18.37}
+\int\limits_{\Im \hw=\hh(\that)}\frac{\sin(\hze_1(\hz,\that)-\hze_1(\hw,\that))}{\cosh(\hze_2(\hz,\that)-\hze_2(\hw,\that))-\cos(\hze_1(\hz,\that)-\hze_1(\hw,\that))}[\pa_{\hz}\hze_\mu(\hz,\that)-\pa_{\hw}\hze_\mu(\hw,\that)]d\hw.\end{equation}
This holds for $|\Im \hz|\leq \hh(\that)$.

We define

\begin{equation}\label{18.38}
\hat{A}(\hz,\that)=-[(\pa_t\phi)\circ\Psi(\hz,\that)]\end{equation}
for $|\Im \hz|\leq \hh(\that)$, $\that\in [t_0-\delta,t_0+\delta]$.

Thus
\begin{equation}\label{18.39}
\text{$\hat{A}(\hz,\that)$ is smooth in $(\hz,\that)$, and holomorphic in $\hz$ for fixed $\that$.}\end{equation}
\begin{rem}
Note that $\hat{A}(\hz,\that)$ does not depend on our Muskat solution. It is determined entirely by $h(x,t)$ and $\Phi(\zeta,t)$.\end{rem}
We may rewrite (\ref{18.37}) in the equivalent form
$$\pa_{\that}\hze_\mu(\hz,\that)=$$
$$\left\{\hat{A}(\hz,\that)+\int\limits_{\Im \hw=\hh(\that)}\frac{\sin(\hze_1(\hz,\that)-\hze_1(\hw,\that))\hat{B}(\hz,\hw,\that)}{\cosh(\hze_2(\hz,\that)-\hze_2(\hw,\that))-\cos(\hze_1(\hz,\that)-\hze_1(\hw,\that))}d\hw\right\}\pa_{\hz}\hze_\mu(\hz,\that)$$
\begin{equation}\label{18.40}
+\int\limits_{\Im \hw=\hh(\that)}\frac{\sin(\hze_1(\hz,\that)-\hze_1(\hw,\that))}{\cosh(\hze_2(\hz,\that)-\hze_2(\hw,\that))-\cos(\hze_1(\hz,\that)-\hze_1(\hw,\that))}[\pa_{\hz}\hze_\mu(\hz,\that)-\pa_{\hw}\hze_\mu(\hw,\that)]d\hw\end{equation}
We call (\ref{18.40}) the TRANSFORMED MUSKAT EQUATION or 'TME'.

Our additional conditions (\ref{18.9}),...,(\ref{18.14}) are easily seen to be equivalent to the following:
\begin{equation}\label{18.41}
\text{$\hze_\mu(\hz,\that)$ is holomorphic in $\{|\Im\hz|<\hh(\that)\}$ for fixed $\that$}.\end{equation}
\begin{equation}\label{18.42}
\text{$\hze_\mu(\hz,\that)-\hat{\zeta}\delta_{1,\, \mu}\in H^{4}\left(\{\Im \hz=\pm \hh(\that)\}\right)$ for fixed $\that$.}\end{equation}
\begin{equation}\label{18.43}
\text{$\hze_1(\hz,\that)-\hz$ and $\hze_2(\hz,\that)$ are $2\pi-$periodic for fixed $\that$.}
\end{equation}
\begin{equation}\label{18.44}
\text{$\hze_\mu(\hz,\that)$ is real for real $\hz$}.
\end{equation}
The initial data $\hze_\mu(\hz,t_0)=\hze^0_\mu(\hz)$ are given and satisfy (\ref{18.41}),...,(\ref{18.44}), as well as the chord-arc condition
$$|\cosh(\hze_2(\hz,t_0)-\hze_2(\hw,t_0))-\cos(\hze_1(\hz,t_0)-\hze_1(\hw,t_0))|$$
\begin{equation}\label{18.45}
\geq c \left[||\Re (\hz-\hw)||+|\Im(\hz-\hw)|\right]^2 \end{equation}
for $|\Im \hz|$, $|\Im \hw|\leq \hh(t_0)$.

Here, of course,
\begin{equation}\label{18.46}
\hze_\mu^0(\hz)=z_\mu^0\circ \Psi(\hz,t_0).\end{equation}
To deduce the equivalence of $(\ref{18.12})$ and $(\ref{18.44})$, we use the fact that
\begin{eqnarray}
&\text{$\hat{A}(\hz,\that)$ is real for real $\hz$; and}\nonumber\\
&\text{$\hat{B}(\hz,\hw,\that)$ is real for real $\hz$, $\hw$}\label{18.47}
\end{eqnarray}
and we change the contour in (\ref{18.40}) to $\{\Im \hw=0\}$, i.e., our contour may be taken to be the circle $\T$.

We use also the fact that
\begin{eqnarray}
&\text{$\hat{A}(\hz,\that)$ is $2\pi-$periodic in $\hz$ for fixed $\that$; and that}\nonumber\\
&\text{$\hat{B}(\hz,\hw,\that)$ is $2\pi-$periodic  in both $\hz$ and $\hw$ for fixed $\that$.}\label{18.48}
\end{eqnarray}
Observe that (\ref{18.40}),...,(\ref{18.45})  are equivalent to (\ref{18.8}),...,(\ref{18.14}).

It remains to formulate a 'Generalized Rayleigh-Taylor Condition' for (\ref{18.40}),...,(\ref{18.45}), and to prove the equivalence of that condition to (\ref{18.16}).

We define (for $\Im \hz=\hh(\that)$)
$$\hat{RT}(\hz,\that)\equiv$$
$$\Re\left\{\frac{-2\pi\pa_{\hz}\hze_1(\hz,\that)}{(\pa_{\hz}\hze_1(\hz,\that))^2+(\pa_{\hz}\hze_2(\hz,\that))^2}\right\}$$
$$+\Im\left\{ P.V.\int\limits_{\Im \hw=\hh(\that)}\frac{\sin(\hze_1(\hz,\that)-\hze_1(\hw,\that))}{\cosh(\hze_2(\hz,\that)-\hze_2(\hw,\that))-\cos(\hze_1(\hz,\that)-\hze_1(\hw,\that))}(\hat{B}(\hz,\hw,\that)+1)d\hw\right\}$$
\begin{equation}\label{18.49}
+\Im\{\hat{A}(\hz,\that)\}+\hh'(\that).\end{equation}

We will see that (\ref{18.49}) is the natural Rayleigh-Taylor function for the problem (\ref{18.40}),...,(\ref{18.45}).

Our 'Generalized Rayleigh-Taylor Condition' for (\ref{18.40}),...,(\ref{18.45}) is
\begin{equation}\label{18.50}
\hat{RT}(\hz,t_0)>0\end{equation}
(strict positivity) for all $\hz$ such that $\Im \hz=\hh(\that)$.

We hope to produce a solution to (\ref{18.40}),...,(\ref{18.45}), provided our initial data $\hze_\mu^0$ satisfy the chord-arc condition (\ref{18.45}) and the 'Generalized Rayleigh-Taylor condition' (\ref{18.50}).

This will be done in a later section; for the moment, we simply note that our solution to (\ref{18.40}),...,(\ref{18.45}) will later be shown to exist in a short time interval $[t_0-\delta',t_0]$ where $\delta'$ will be a small enough positive number determined by the constants in our assumptions on the initial data.

For the rest of this section, our task is to prove that condition (\ref{18.16}) implies condition (\ref{18.50}).

To do so, we first make a few elementary remarks on the maps $\Phi$, $\Psi$, $\phi$, $\psi$ in (\ref{18.2}),...,(\ref{18.7}). Then we return to (\ref{18.16}) and (\ref{18.50}).

First of all, since $\zeta\mapsto\phi(\zeta,t)$ maps $\Gamma_+(t)$ to $\{\Im \hz=\hh(t)\}$ for each $t$, we have
\begin{equation}\label{18.51}
\pa_{\zeta}\phi(\zeta,t)=\rho(\zeta,t)[1+i\pa_x h(x,t)]^{-1}\end{equation}
for $x=\Re\zeta$, $\Im \zeta=h(x,t)$; here,
\begin{equation}\label{18.51a}
\rho(\zeta,t)>0\end{equation}
for $x=\Re\zeta$, $\Im \zeta=h(x,t)$.
The operator
\begin{equation}\label{18.52}
\pa_t+[\pa_th(x,t)]\frac{\pa}{\pa \Im \zeta}\quad \text{corresponds under $\Phi$ to}\quad
 \pa_{\that}+[\hh'(t)]\frac{\pa}{\pa \Im \hz}+\hat{b}(\hz,\that)\frac{\pa}{\pa \Re \hz},\end{equation}
 \begin{equation}\label{18.53}\text{for $x=\Re\zeta$, $\Im \zeta=h(x,t)$, $(\hz,\that)=\Phi(\zeta,t)$;}\end{equation}
here
\begin{equation}\label{18.54}
\text{$\hat{b}(\hz,\that)$ is real whenever (\ref{18.53}) holds.}\end{equation}
Indeed, (\ref{18.51}), (\ref{18.51a}) hold because $\Gamma_+(t)$ maps to $\{\Im \hz=\hh(\that)\}$ with positive orientation (for each fixed $t$); and (\ref{18.52}),...,(\ref{18.54}) hold because a moving particle that stays on the upper boundary $\Gamma_+(t)$ as $t$ varies will stay on the upper boundary $\Im \hz=\hh(\that)$ when viewed in the $(\hz,\that)$ coordinate system.

Let $F(\zeta,t)$ be a holomorphic function of $\zeta$ for fixed $t$; let $\hat{F}=F\circ \Psi$. Then (\ref{18.52}) gives
\begin{equation}\label{18.55}
\pa_tF(\zeta,t)+[\pa_t h(x,t)]i\pa_\zeta F(\zeta,t)=\pa_{\that}\hat{F}(\hz,\that)+\hh'(\that)i\pa_{\hz}\hat{F}(\zeta,\that)+\hat{b}(\hz,\that)\pa_{\hz}\hat{F}(\hz,\that),\end{equation}
whenever (\ref{18.53}) holds.

On the other hand, since $F(\zeta,t)=\hat{F}(\phi(\zeta,t),t)$, the chain rule gives
\begin{equation}\label{18.56}
\pa_t F(\zeta,t)=\pa_{\that}\hat{F}(\hz,\that)+[\pa_t\phi(\zeta,t)]\pa_{\hz}\hat{F}(\hz,t)\end{equation}
whenever $(\hz,\that)=\Phi(\zeta,t)$, in particular, whenever (\ref{18.53}) holds.

Subtracting (\ref{18.56}) from (\ref{18.55}), we find that
\begin{equation}\label{18.57}
[i\pa_t h(x,t)]\pa_\zeta F(\zeta,t)=[i\hh'(\that)+\hat{b}(\hz,\that)-\pa_t\phi(\zeta,t)]\pa_{\hz}\hat{F}(\hz,\that),\end{equation}
whenever (\ref{18.53}) holds.

The chain rule gives
$$\pa_\zeta F(\zeta,t)=[\pa_\zeta\phi(\zeta,t)]\pa_{\hz}\hat{F}(\hz,\that)$$
whenever (\ref{18.53}) holds. Putting this into (\ref{18.57}) and recalling that $F(\zeta,t)$ may be any function holomorphic in $\zeta$ for fixed $t$, we see that
$$i\pa_th(x,t)[\pa_\zeta\phi(\zeta,t)]=i\pa_{\hat{t}}\hh(\that)-\hat{b}(\hz,\that)-\pa_t\phi(\zeta,t)$$
whenever (\ref{18.53}) holds. Recalling (\ref{18.51}), we now see that
\begin{equation}\label{18.58}
\frac{i\pa_th(x,t)\rho(\zeta,t)}{1+i\pa_xh(x,t)}=i\pa_{\hat{t}}\hh(\that)+\hat{b}(\hz,\that)-\pa_t\phi(\zeta,t).\end{equation}

Next, we note that since $z_\mu(\zeta,t)=\hze_\mu(\phi(\zeta,t),t)$, the chain rule gives
$$\pa_\zeta z_\mu(\zeta,t)=[\pa_\zeta \phi(\zeta,t)]\pa_{\hz}\hze _\mu(\hz,\that)$$
whenever $(\hz,\that)=\Phi(\zeta,t)$, in particular, whenever (\ref{18.53}) holds.

Consequently
$$\frac{-2\pi \pa_\zeta z_1(\zeta,t)[\pa_\zeta\phi(\zeta,t)]}{(\pa_\zeta z_1(\zeta,t))^2+(\pa_{\zeta} z_2(\zeta,t))^2}=\frac{-2\pi \pa_{\hz} \hat{z}_1(\hz,t)}{(\pa_{\hz} \hat{z}_1(\hz,t))^2+(\pa_{\hz} \hat{z}_2(\hz,t))^2}$$
whenever (\ref{18.53}) holds. Recalling (\ref{18.51}), we conclude that
\begin{equation}\label{18.59}
\frac{-2\pi \pa_\zeta z_1(\zeta,t)}{(\pa_\zeta z_1(\zeta,t))^2+(\pa_{\zeta} z_2(\zeta,t))^2}\frac{\rho(\zeta,t)}{1+i\pa_x h(x,t)}=\frac{-2\pi \pa_{\hz} \hat{z}_1(\hz,t)}{(\pa_{\hz} \hat{z}_1(\hz,t))^2+(\pa_{\hz} \hat{z}_2(\hz,t))^2}\end{equation}
whenever (\ref{18.53}) holds.

Next, substituting the definition (\ref{18.33}) of $\hat{B}(\hz,\hw,\that)$, and changing variable by $\hz=\phi(\zeta,t)$, $\hw=\phi(w,t)$, $\that=t$, we find that
$$P.V.\int\limits_{\Im \hw=\hh(\that)}\frac{\sin(\hze_1(\hz,\that)-\hze_1(\hw,\that))}{\cosh(\hze_2(\hz,\that)-\hze_2(\hw,\that))-\cos(\hze_1(\hz,\that)-\hze_1(\hw,\that))}(\hat{B}(\hz,\hw,\that)+1)d\hw$$
$$=P.V.\int\limits_{\Im \hw=\hh(\that)}\frac{\sin(\hze_1(\hz,\that)-\hze_1(\hw,\that))}{\cosh(\hze_2(\hz,\that)-\hze_2(\hw,\that))-\cos(\hze_1(\hz,\that)-\hze_1(\hw,\that))}[(\pa_\zeta\phi)\circ \Psi(\hz,\that)][\pa_{\hw}\psi(\hw,\that)]d\hw$$
$$=P.V.\int\limits_{w\in \Gamma_+(t)}\frac{\sin(z_1(\zeta,t)-z_1(w,t))}{\cosh(z_2(\zeta,t)-z_2(w,t))-\cos(z_1(\zeta,t)-z_1(w,t))})dw\pa_\zeta \phi(\zeta,t).$$
Recalling (\ref{18.51}), we conclude that
$$P.V.\int\limits_{\Im \hw=\hh(\that)}\frac{\sin(\hze_1(\hz,\that)-\hze_1(\hw,\that))}{\cosh(\hze_2(\hz,\that)-\hze_2(\hw,\that))-\cos(\hze_1(\hz,\that)-\hze_1(\hw,\that))}(\hat{B}(\hz,\hw,\that)+1)d\hw$$
\begin{equation}\label{18.60}
=\frac{\rho(\zeta,t)}{1+i\pa_xh(x,t)}P.V.\int\limits_{w\in \Gamma_+(t)}\frac{\sin(z_1(\zeta,t)-z_1(w,t))}{\cosh(z_2(\zeta,t)-z_2(w,t))-\cos(z_1(\zeta,t)-z_1(w,t))})dw\end{equation}
whenever (\ref{18.53}) holds.

Next, we apply (\ref{18.38}) and (\ref{18.58}). Thus, whenever (\ref{18.53}) holds, we have
\begin{equation}\label{18.61}
\hat{A}(\hz,\that)=-\pa_t\phi(\zeta,t)=\frac{i\pa_th(x,t)\rho(\zeta,t)}{1+i\pa_xh(x,t)}-i\hh'(\that)-\hat{b}(\hz,\that).\end{equation}
Now, substituting (\ref{18.59}), (\ref{18.60}) and (\ref{18.61}) into the definition (\ref{18.49}) of $\hat{RT}(\hz,\that)$, we find that, whenever (\ref{18.53}) holds, we have
$$\hat{RT}(\hz,\that)=$$
$$\Re\left\{ \frac{-2\pi \pa_\zeta z_1(\zeta,t)}{(\pa_\zeta z_1(\zeta,t))^2+(\pa_{\zeta} z_2(\zeta,t))^2}\frac{\rho(\zeta,t)}{1+i\pa_x h(x,t)}\right\}$$
$$+\Im\left\{\frac{\rho(\zeta,t)}{1+i\pa_x h(x,t)}P.V.\int\limits_{w\in \Gamma_+(t)}\frac{\sin(z_1(\zeta,t)-z_1(w,t))}{\cosh(z_2(\zeta,t)-z_2(w,t))-\cos(z_1(\zeta,t)-z_1(w,t))})dw\right\}$$
\begin{equation}\label{18.62}
+\Im\left\{\frac{i\pa_th(x,t)\rho(\zeta,t)}{1+i\pa_xh(x,t)}-i\hh'(\that)-\hat{b}(\hz,\that)\right\}+\hh'(\that).\end{equation}
Note that the $\hh'(\that)$ terms cancel, and that $\hat{b}(\hz,\that)$ does not  contribute to the right-hand side of (\ref{18.62}), thanks to (\ref{18.54}).

Recalling also (\ref{18.51a}), we see that (\ref{18.62}) may be rewritten in the equivalent form
 $$\hat{RT}(\hz,\that)=$$
$$\rho(\zeta,t)\left[\Re\left\{ \frac{-2\pi \pa_\zeta z_1(\zeta,t)}{(\pa_\zeta z_1(\zeta,t))^2+(\pa_{\zeta} z_2(\zeta,t))^2}\frac{1}{1+i\pa_x h(x,t)}\right\}\right.$$
$$+\Im\left\{\frac{1}{1+i\pa_x h(x,t)}P.V.\int\limits_{w\in \Gamma_+(t)}\frac{\sin(z_1(\zeta,t)-z_1(w,t))}{\cosh(z_2(\zeta,t)-z_2(w,t))-\cos(z_1(\zeta,t)-z_1(w,t))})dw\right\}$$
\begin{equation*}
+\left.\Re\left\{\frac{\pa_th(x,t)}{1+i\pa_xh(x,t)}\right\}\right]\end{equation*}
whenever (\ref{18.53}) holds.

Comparing this equation with the definition (\ref{18.15}) of $RT(\zeta,t)$, we conclude that
\begin{equation}\label{18.63}
\hat{RT}(\hz,\that)=\rho(\zeta,t)RT(\zeta,t)\end{equation}
whenever (\ref{18.53}) holds.

Given $\hz$, $\that$ with $\Im \hz =\hh(\that)$, we set $(\zeta,t)=\Psi(\hz,\that)$ and $x=\Re \zeta$. Then $\Im \zeta=h(x,t)$, $x=\Re \zeta$, and $(\hz,\that)=\Phi(\zeta,t)$, i.e., (\ref{18.53}) holds. Therefore, (\ref{18.63}) applies. We recall from (\ref{18.51a}) that $\rho(\zeta,t)$ is strictly positive. If (\ref{18.16}) holds for $z_\mu(\zeta,t)$, then $RT(\zeta,t)$ is also strictly positive. Therefore, by (\ref{18.63}), we have $\hat{RT}(\hz,\that)>0$ (strict inequality). This is the desired 'Generalized Rayleigh-Taylor Condition' (\ref{18.50}). Thus, as claimed, (\ref{18.50}) follows from (\ref{18.16}).

\subsection{Generalized Muskat}\label{seccion19}
In this section we prove local existence for the 'Transformed Muskat equation' (\ref{18.40}),...,(\ref{18.44}). Let's start with some definitions.
\begin{equation}\label{19.1}
\text{Let $h(t)$ be a positive, smooth function defined for $t\in I_{time}=[t_0-\delta,t_0]$.}\end{equation}
Define
\begin{equation}\label{19.2}
\Gamma_\pm(t)=\{\zeta\in\C\,:\, \Im \zeta = \pm h(t)\}.\end{equation}
\begin{equation}\label{19.3}
\Omega(t)=\{\zeta\in\C\,:\,|\Im\zeta|<h(t)\}\end{equation}
\begin{equation}\label{19.4}
\Omega=\{(\zeta,t)\,:\,\zeta\in \Omega(t)^{closure},\,t\in I_{time}\}\end{equation}
\begin{equation}\label{19.5}
\Omega^+=\{(\zeta,w,t)\,:\,\zeta\in \Omega(t)^{closure},\,w\in\Omega(t)^{closure},\,t\in I_{time}\}\end{equation}
\begin{equation}\label{19.6}
\text{Let $A(\zeta,t)$ be a smooth, complex-valued function on $\Omega$.}\end{equation}
Assume:
\begin{equation}\label{19.7}
\text{$A(\zeta,t)$ is real when $\zeta$ is real;}\end{equation}
\begin{equation}\label{19.8}
A(\zeta+2\pi,t)=A(\zeta,t); \end{equation}
and
\begin{equation}\label{19.9}
\text{$\zeta\mapsto A(\zeta,t)$ is holomorphic on $\Omega(t)$ for fixed $t$.}\end{equation}
\begin{equation}\label{19.10}
\text{Let $B(\zeta,w,t)$ be a smooth, complex-valued function on $\Omega^+$.}\end{equation}
Assume that
\begin{equation}\label{19.11}
\text{$B(\zeta,w,t)$ is real when $\zeta$, $w$ are real;}\end{equation}
\begin{equation}\label{19.12}
B(\zeta+2\pi,w,t)=B(\zeta,w+2\pi,t)=B(\zeta,w,t);\end{equation}
\begin{equation}\label{19.13}
\text{$(\zeta,w)\mapsto B(\zeta,w,t)$ is holomorphic on $\Omega(t)\times \Omega(t)$ for fixed $t$};\end{equation}
\begin{equation}\label{19.14}
B(\zeta,\zeta,t)=0.\end{equation}
\begin{equation}\label{19.15}
\text{For $\mu=1,2$, let $z^0_\mu(\zeta)$ be a holomorphic function on $\Omega(t_0)$.}\end{equation}
Assume that
\begin{equation}\label{19.16}
z_1^0(\zeta+2\pi)=z_1^0(\zeta)+2\pi\quad\text{and}\quad z_2^0(\zeta+2\pi)=z_2^0(\zeta);\end{equation}
\begin{equation}\label{19.17}
\text{$z_\mu^0(\zeta)$ is real whenever $\zeta$ is real.}\end{equation}
Assume that
\begin{equation}\label{19.18}
\text{$z_\mu^0(\zeta)$ extends to a continuous function on $\Omega(t_0)^{closure}$,}\end{equation}
such that
\begin{equation}\label{19.19}
\text{$(z_1^0(\zeta)-\zeta)\Big|_{\Gamma_{\pm}(t_0)}$, $z_2^0(\zeta)\Big|_{\Gamma_{\pm}(t_0)}$ belong to $H^{4}$}.\end{equation}
We assume also that $z_\mu^0$ satisfy the chord-arc condition:
\begin{equation}\label{19.20}
|\cosh(z_2^0(\zeta)-z_2^0(w))-\cos(z_1^0(\zeta)-z_1^0(w))|\geq c_{CA} \left[||\Re (\zeta-w)||+|\Im(\zeta-w)|\right]^2 \end{equation}
for $\zeta$, $w\in \Omega(t_0)^{closure}$.

Also, we define the Rayleigh-Taylor function for the initial data $z_\mu^0$ as follows.

For $\zeta\in \Gamma_+(t_0)$, we define
$$RT^0(\zeta)=$$
$$\Re\left\{\frac{-2\pi(z_1^0)'(\zeta)}{((z_1^0)'(\zeta))^2+((z_2^0)'(\zeta))^2}\right\}$$
$$+\Im \left\{P.V.\int\limits_{w\in\Gamma_+(t)}\frac{\sin(z_1^0(\zeta)-z_1^0(w))}{\cosh(z_2^0(\zeta)-z_2^0(w))-\cos(z_1^0(\zeta)-z_1^0(w))}[B(\zeta,w,t_0)+1]dw\right\}$$
\begin{equation}\label{19.21}
+\Im A(\zeta,t_0)+h'(t).\end{equation}

We assume that the following Rayleigh-Taylor condition holds for the initial data:
\begin{equation}\label{19.22}
RT^0(\zeta)>0\quad \text{for all $\zeta\in\Gamma_+(t_0)$} \end{equation}
Note that we assume strict inequality here; we do not allow $RT^0(\zeta)$ to vanish at any point $\zeta\in\Gamma_+(t_0)$.

Next, we introduce the notation of a 'controlled constant'.

A controlled constant is a positive real number determined entirely by the following data.

An upper bound for the $C^{4+100}-$norm of $h$ on $I_{time}$.

A positive lower bound for $h(t_0)$.

An upper bound for the $C^{4+10}-$norms of $A(\zeta,t)$ and $B(\zeta,w,t)$.

An upper bound for the $H^{4}\,-$norm of $z_1^0(\zeta)-\zeta$ and $z_1^0(\zeta)$ on $\Gamma_{\pm}(t_0)$.

The constant $c_{CA}$ in the Chord-Arc condition.

A positive lower bound for the Rayleigh-Taylor function $RT^0(\zeta)$ on $\Gamma_+(t_0)$.

We write $c$, $C$, $C'$, etc. to denote controlled constants.

These symbols may denote different constants in different occurrences.
$$\text{We assume that the length $\delta$ of the time interval $I_{time}$ above, is less than a small enough}$$\begin{equation}\label{19.23}\text{controlled constant.}\end{equation}
Under the above assumptions, we would like to prove the existence of a solution to the following GENERALIZED MUSKAT PROBLEM.
\begin{equation}\label{19.24}
\text{Find $C^1$ functions $z_\mu$ on $\Omega$ satisfying the following conditions:}\end{equation}
\begin{equation}\label{19.25}
\text{$z_\mu(\zeta,t)$ is real whenever $\zeta$ is real.}\end{equation}
\begin{equation}\label{19.26}
\text{$z_1(\zeta,t)-\zeta$ and $z_2(\zeta,t)$ are $2\pi-$periodic function in $\zeta$ for fixed $t$.}\end{equation}
\begin{equation}\label{19.27}
\text{$z_\mu(\zeta,t)$ is holomorphic on $\Omega(t)$ for fixed $t$.}\end{equation}
\begin{equation}\label{19.28}
\text{$z_1(\zeta,t)-\zeta$ and $z_2(\zeta,t)$ belong to $H^{4}(\Gamma_\pm(t))$ for fixed $t$,}\end{equation}
with norm bounded by a controlled constant.
\begin{equation}\label{19.29}
z_\mu(\zeta,t_0)=z_\mu^0(\zeta)\end{equation}
\begin{equation}\label{19.30}
|\cosh(z_2(\zeta,t)-z_2(w,t))-\cos(z_1(\zeta,t)-z_1(w,t))|\geq C'_{CA} \left[||\Re (\zeta-w)||+|\Im(\zeta-w)|\right]^2\end{equation}
for $(\zeta,w,t)\in \Omega^+$.

For $t\in I_{time}$, $\zeta\in \Gamma_+(t)$, define
$$RT(\zeta)=$$
$$\Re\left\{\frac{-2\pi(\pa_\zeta z_1)(\zeta,t)}{((\pa_\zeta z_1)(\zeta,t))^2+((\pa_\zeta z_2)(\zeta,t))^2}\right\}$$
$$+\Im \left\{P.V.\int\limits_{w\in\Gamma_+(t)}\frac{\sin(z_1(\zeta,t)-z_1(w))}{\cosh(z_2w(\zeta,t)-z_2(w))-\cos(z_1(\zeta,t)-z_1(w))}[B(\zeta,w,t)+1]dw\right\}$$
\begin{equation}\label{19.31}
+\Im A(\zeta,t)+h'(t).\end{equation}

Then
\begin{equation}\label{19.32}
\text{$RT(\zeta,t)>0$ for $t\in I_{time}$, $\zeta\in\Gamma_+(t)$.}\end{equation}

For $(\zeta,t)\in \Omega$, the following equation holds:

$$\pa_{t}z_\mu(\zeta,t)=$$
$$\left\{A(\zeta,t)+\int\limits_{w\in\Gamma_+(t)}\frac{\sin(z_1(\zeta,t)-z_1(w,t))B(\zeta,w,t)}{\cosh(z_2(\zeta,t)-z_2(\hw,t))-\cos(z_1(\zeta,t)-z_1(w,t))}dw\right\}\pa_{\zeta}z_\mu(\zeta,t)$$
\begin{equation}\label{19.33}
+\int\limits_{w\in\Gamma_+(t)}\frac{\sin(z_1(\zeta,t)-z_1(w,t))}{\cosh(z_2(\zeta,t)-z_2(w,t))-\cos(z_1(\zeta,t)-z_1(w,t))}[\pa_{\zeta}z_\mu(\zeta,t)-\pa_{w}z_\mu(w,t)]dw.\end{equation}
In the above integrals, the integrands are holomorphic in $w\in\Omega(t)$ for fixed $\zeta$ and $t$; in particular, the integrands have no poles. Hence, we may change the contours in these integrals from $\Gamma_+(t)$ to $\T$, without changing the values of the integrals.

In the next section, we prove existence of solutions to the above GENERALIZED MUSKAT PROBLEM.

Thanks to the results in  section (\ref{seccion18}) this yields the desired local existence of solutions of our original Muskat equation, with desired holomorphic and $H^{4}$ properties, etc.

\subsection{Galerkin Approximation}\label{seccion20}
\begin{equation}\label{20.1}
\text{Let $h(t)$, $\Gamma_\pm(t)$, $t_0$, $\delta$, $\Omega(t)$, $\Omega$, $\Omega^+$, $A(\zeta,t)$, $B(\zeta,w,t)$, $z_\mu^0(\zeta)$ be as in section (\ref{seccion19})}.\end{equation}
For a large enough positive integer $N$, we define $z_\mu^{[N]}(\zeta)$ from $z_\mu^0(\zeta)$, by using the projection
\begin{equation}\label{20.2}
\Pi_N\,:\, \sum_{-\infty}^\infty A_k e^{ik\zeta}\mapsto \sum_{-N}^N A_k e^{ik\zeta}.\end{equation}
We define $z_\mu^{[N]}(\zeta)$ by stipulating that
\begin{equation}\label{20.3}
z_1^{[N]}(\zeta)-\zeta=\Pi_N[z_1^0(\zeta)-\zeta]\end{equation}
and
\begin{equation}\label{20.4}
z_2^{[N]}(\zeta)=\Pi_N[z_2^0(\zeta)].\end{equation}
For $N$ large enough, the functions $z_\mu^{[N]}$ satisfy the chord-arc and Rayleigh-Taylor condition (\ref{19.30}), (\ref{19.31}), (\ref{19.32}) with controlled constants. Also, since $z_2^{[N]}(\zeta)$ and $z_2^{[N]}(\zeta)-\zeta$ are trigonometric polynomials, those function are holomorphic on $\Omega(t_0)$, they extend to smooth functions on the closure of $\Omega(t_0)$, and their $H^{4}-$ norms on $\Gamma_{t_0}$ are at most the $H^{4}-$norms of $z_2^0$ and $z_1^0(\zeta)-\zeta$, respectively, and are thus dominated by a controlled constant.

Consequently, the conditions imposed in section (\ref{seccion19}) for $z_\mu^0$ hold also for $z_\mu^{[N]}$, with controlled constants, provided $N$ is large enough. From now on, we suppose that $N$ is large enough that this holds.

In place of the Generalized Muskat equation (\ref{19.33}), we will solve its Galerkin approximation. That is, we will solve the equation
\begin{equation}\label{20.5}
\pa_t z_\mu(\zeta,t)=\Pi_N[J(\cdot,t)](\zeta,t),\end{equation}
where
$$J(\zeta,t)=$$
$$\left\{A(\zeta,t)+\int\limits_{w\in\Gamma_+(t)}\frac{\sin(z_1(\zeta,t)-z_1(w,t))B(\zeta,w,t)}{\cosh(z_2(\zeta,t)-z_2(w,t))-\cos(z_1(\zeta,t)-z_1(w,t))}dw\right\}\pa_{\zeta}z_\mu(\zeta,t)$$
\begin{equation}\label{20.6}
+\int\limits_{w\in\Gamma_+(t)}\frac{\sin(z_1(\zeta,t)-z_1(w,t))}{\cosh(z_2(\zeta,t)-z_2(w,t))-\cos(z_1(\zeta,t)-z_1(w,t))}[\pa_{\zeta}z_\mu(\zeta,t)-\pa_{w}z_\mu(w,t)]dw\end{equation}
We impose the initial condition
\begin{equation}\label{20.7}
z_\mu(\zeta,t_0)=z_\mu^{[N]}(\zeta).\end{equation}
Here, $z_\mu(\zeta,t)$ depend on $N$, though our notation does not indicate that.

We hope to solve (\ref{20.5}), (\ref{20.6}), (\ref{20.7}) for $t\in I_{time}$.

Note that (\ref{20.5}), (\ref{20.6}), (\ref{20.7}) is simply a system of ODE's. We solve these ODE's, starting at time $t_0$, and proceeding backwards in time, until either we reach time $t_0-\delta$ (the left endpoint of the interval $I_{time}$), or else the conditions (\ref{19.24}),...,(\ref{19.32}) (with suitable controlled constants) are not satisfied.

We will see that we get all the way back to time $t_0-\delta$. Thus, we obtain a solution of (\ref{20.5}), (\ref{20.6}), (\ref{20.7}), satisfying conditions (\ref{19.24}),...,(\ref{19.32}) with controlled constants; these controlled constants are independent of $N$, so by \eqref{19.28}, we may pick a sequence $N_\nu\to \infty$ for which the corresponding solutions of (\ref{20.5}), (\ref{20.6}), (\ref{20.7}) tend weakly to a limit $z_\mu(\zeta,t)$ in $H^{4}\left(\Gamma_\pm(t)\right)$ for fixed $t$ (rational).

The resulting $z_\mu(\zeta,t)$ is then easily seen to solve the GENERALIZED MUSKAT PROBLEM in section (\ref{seccion19}).

It all comes down to an a priori estimate
\begin{equation}\label{20.8}
\frac{d}{dt}\int\limits_{\zeta\in \Gamma_{\pm}(t)}\left|\pa_{\zeta}^{4}z_\mu(\zeta,t)\right|^2d\zeta\geq -C',
\end{equation}
to be proven as long as our solution to (\ref{20.5}), (\ref{20.6}), (\ref{20.7}) satisfies (\ref{19.24}),...,(\ref{19.32}). Note that $z_\mu(\overline{\zeta},t)=\overline{z_\mu(\zeta,t)}$ so the integrals over $\Gamma_+(t)$ and $\Gamma_-(t)$ in \eqref{20.8} are equal.

We sketch a proof of (\ref{20.8}). We have
$$\frac{1}{2}\frac{d}{dt}\int\limits_{\zeta\in \Gamma_{\pm}(t)}\left|\pa_{\zeta}^{4}z_\mu(\zeta,t)\right|^2d\zeta$$
$$=\Re \int\limits_{\zeta\in\Gamma_+(t)} \overline{\pa_\zeta^{4}z_\mu(\zeta,t)}\left\{\pa_t \pa_\zeta^{4}z_\mu(\zeta,t)+ih'(t)\pa_\zeta^{4+1}z_\mu(\zeta,t)\right\}d\zeta$$
$$=\Re \int\limits_{\zeta\in\Gamma_+(t)} \overline{\pa_\zeta^{4}z_\mu(\zeta,t)}\left\{\pa_\zeta^{4}\Pi_N[J(\cdot,t)](\zeta,t)+ih'(t)\pa_\zeta^{4+1}z_\mu(\zeta,t)\right\}d\zeta$$
$$=\Re \int\limits_{\zeta\in\Gamma_+(t)} \overline{\pa_\zeta^{4}z_\mu(\zeta,t)}\left\{\Pi_N[\pa_\zeta^{4}J(\cdot,t)](\zeta,t)+ih'(t)\pa_\zeta^{4+1}z_\mu(\zeta,t)\right\}d\zeta$$
\begin{equation}\label{20.9}
=\Re \int\limits_{\zeta\in\Gamma_+(t)} \overline{\pa_\zeta^{4}z_\mu(\zeta,t)}\left\{\pa_\zeta^{4}J(\zeta,t)+ih'(t)\pa_\zeta^{4+1}z_\mu(\zeta,t)\right\}d\zeta\end{equation}
since $\pa_\zeta^ 4 z_\mu(\zeta,t)$, is a trigonometric polynomial in the range of $\Pi_N$.

Thus, $\Pi_N$ has disappeared from the right-hand side of (\ref{20.9}).

The above manipulations are made possible by the fact that the contours $\Gamma_\pm(t)$ are horizontal straight lines. That is why we took the trouble to make the conformal mapping in sections \ref{conformalmap} and \ref{seccion18}.

Using ideas from section \ref{seccion12}, we can compute $\pa_\zeta^4 J(\zeta,t)$ modulo errors with controlled $L^2-$norm on $\Gamma_+(t)$. To do so, we first recall our assumptions that the $H^4-$norms of $z_1(\zeta,t)-\zeta$ and $z_2(\zeta,t)$ on $\Gamma_+(t)$ are bounded by a controlled constant $C$, and that the chord-arc condition holds with a controlled constant $c$. Consequently, for $\zeta\in \Gamma_+(t)$ and $y\in[-\pi,\pi]$, the following estimates hold:
\begin{align}
|\sin(z_1(\zeta,t)-z_1(\zeta+y,t))|\leq C|y|, &\quad |\cos(z_1(\zeta,t)-z_1(\zeta+y,t)|\leq C\label{5.108}\\
|\sinh(z_1(\zeta_2,t)-z_2(\zeta+y,t))|\leq C|y|, & \quad |\cosh(z_2(\zeta,t)-z_2(\zeta+y,t))|\leq C \label{5.109}\\
|\pa_\zeta^kz_\mu(\zeta,t)-\pa_\zeta^k z_\mu(\zeta+y,t)|\leq C|y|& \quad \text{for $0\leq k\leq 2$, $\mu=1,2$}\label{5.110}\\
|\pa_\zeta^3z_\mu(\zeta,t)-\pa_\zeta^3 z_\mu(\zeta+y,t)|\leq C|y|M_\mu(\zeta)\ &\quad \text{for $\mu=1,2$}\label{5.111}
\end{align}
where $M_\mu(\zeta)$ is the Hardy-Littlewood maximal function of $\pa_\zeta^4z_\mu|_{\Gamma_+(t)}$, and consequently
\begin{equation}\label{5.112}
\int_{\Gamma_+(t)}|M_\mu(\zeta)|^2d\zeta\leq C.
\end{equation}
Also, our assumptions on the function $B(\zeta,w,t)$ tell us that
\begin{equation}\label{5.113}
\tilde{B}(\zeta,y,t)\equiv B(\zeta,\zeta+y,t)\quad \text{for $\zeta\in\Gamma_+(t)$, $y\in[-\pi,\pi]$,}
\end{equation}
satisfies
\begin{equation}\label{5.114}
|\pa_\zeta^l \tilde{B}(\zeta,y,t)|\leq C|y|\quad \text{for $0\leq l\leq 4$}.
\end{equation}

We now proceed to the computation of $\pa_\zeta^4 J(\zeta,t)$, for $\zeta\in\Gamma_+(t)$.

By the proof of lemma \ref{lemma6}, we have
\begin{align}
&\pa_\zeta^4\int_{w\in\Gamma_+(t)}\frac{\sin(z_1(\zeta,t)-z_1(w,t))}{\cosh(z_2(\zeta,t)-z_2(w,t))-\cos(z_1(\zeta,t)-z_1(w,t))}\left[\pa_\zeta z_\mu(\zeta,t)-\pa_w z_\mu(w,t)\right]dw\nonumber\\
&=\text{Dangerous}_\mu(\zeta,t)+\sum_{j=1}^6\text{Safe}_j^\mu(\zeta,t)+\sum_{\nu=1}^{\nu_{max}}c_\nu^\mu \text{Easy}_\nu^\mu(\zeta,t), \label{5.115}
\end{align}
with Dangerous$_\mu$, Safe$_j^\mu$ and Easy$_\nu^\mu$ defined as in lemma \ref{lemma6}, with $\underline{k}=4$.

Using our estimates \eqref{5.108},...,\eqref{5.112}, one shows easily that
\begin{equation}\label{5.116}
|\text{Easy}_\nu^\mu(\zeta,t)|\leq C \quad \text{for each $\mu$, $\nu$}.\end{equation}

To see how to estimate the "Safe" terms, we take for example

\begin{align}
&\text{Safe}_2^\mu(\zeta,t)=\nonumber\\
&(\text{coeff})\int_{-\pi}^\pi \frac{\left[\pa_\zeta z_2(\zeta,t)-\pa_\zeta z_2(\zeta+y,t)\right]\sinh(z_2(\zeta,t)-z_2(\zeta+y,t))\sin(z_1(\zeta,t)-z_1(\zeta+y,t))}
{\left(\cosh(z_2(\zeta,t)-z_2(\zeta+y,t),t)-\cos(z_1(\zeta,t)-z_1(\zeta+y,t))\right)^2}\nonumber\\
&\times \left(\pa_\zeta^4z_\mu(\zeta,t)-\pa_\zeta^4z_\mu(\zeta+y,t)\right)dy\label{5.117}
\end{align}
Since $z_1(\zeta,t)-\zeta$ and $z_2(\zeta,t)$ have $H^4-$norms at most $C$, and since the chord-arc
condition holds with a controlled constant $c$, one shows easily that the expression in curly brackets in
\eqref{5.117} has the form $\frac{\theta(\zeta,y,t)}{y}$, with estimates
$$|\theta(\zeta,y,t)|\leq C$$
and
$$|\pa_y \theta(\zeta,y,t)|\leq C.$$
Consequently, \eqref{5.117} and the $L^2-$norm boundedness of the Hilbert transform tell us that
\begin{equation*}
\int_{\Gamma_+(t)}\left|\text{Safe}_2^\mu(\zeta,t)\right|^2d\zeta\leq C.
\end{equation*}
We can treat the other "Safe" terms in \eqref{5.115} in the same way.

Together with \eqref{5.116} and the definition of the "Dangerous" term in \eqref{5.115}, this now implies that
\begin{align}
&\pa_\zeta^4\int_{w\in\Gamma_+(t)}\frac{\sin(z_1(\zeta,t)-z_1(w,t))}{\cosh(z_2(\zeta,t)-z_2(w,t))-\cos(z_1(\zeta,t)-z_1(w,t))}
\left[\pa_\zeta z_\mu(\zeta,t)-\pa_\zeta z_\mu (w,t)\right]dw\nonumber\nonumber\\
&\int_{w\in\Gamma_+(t)}\frac{\sin(z_1(\zeta,t)-z_1(w,t))}{\cosh(z_2(\zeta,t)-z_2(w,t))-\cos(z_1(\zeta,t)-z_1(w,t))}
\left[\pa_\zeta^5 z_\mu(\zeta,t)-\pa_\zeta^5 z_\mu (w,t)\right]dw\nonumber\\
&+ \text{Error}_\mu(\zeta,t),\label{5.118}
\end{align}
with
\begin{equation}
\int_{\Gamma_+(t)}|\text{Error}_\mu(\zeta,t)|^2d\zeta\leq  C.\label{5.119}
\end{equation}
Thus, modulo errors with bounded $L^2-$norms, we have computed the fourth $\zeta-$derivative of the last integral in
\eqref{20.6}.

We now turn our attention to the integral involving $B(\zeta,w,t)$ in \eqref{20.6}. Let $0\leq k\leq 4$. As in the proof of lemma \ref{lemma6}, we find that
\begin{equation}\label{5.120}
\pa_\zeta^k\int_{w\in \Gamma_+(t)}\frac{\sin(z_1(\zeta,t)-z_1(w,t))}{\cosh(z_2(\zeta,t)-z_2(w,t))-\cos(z_1(\zeta,t)-z_1(w,t))}B(\zeta,w,t)dw
\end{equation}
may be written as a linear combination (with harmless coefficients) of terms of the form

\begin{align}
\int_{-\pi}^\pi&\left[\sin(z_1(\zeta,t)-z_1(\zeta+y,t))\right]^{k'}\left[\cos(z_1(\zeta,t)-z_1(\zeta+y,t))\right]^{m'}\nonumber\\
&\left[\sinh(z_1(\zeta,t)-z_1(\zeta+y,t))\right]^{k''}\left[\cosh(z_1(\zeta,t)-z_1(\zeta+y,t))\right]^{m''}\nonumber\\
& \pa_\zeta^l \tilde{B}(\zeta,y,t)\prod_{\nu=1}^{\nu'_{max}}\left[\pa_\zeta^{k'_\nu}z_1(\zeta,t)-\pa_\zeta^{k'_\nu}z_1(\zeta+y,t)\right]
\prod_{\nu=1}^{\nu_{max}''}\left[\pa_\zeta^{k''_\nu}z_2(\zeta,t)-\pa_\zeta^{k''_\nu}z_2(\zeta+y,t)\right]\nonumber\\
&\left[\cosh(z_2(\zeta,t)-z_2(\zeta+y,t))-\cos(z_1(\zeta,t)-z_1(\zeta+y,t))\right]^{-m}dy\label{5.121}
\end{align}
where $m\geq1$, each $k'_\nu$, $k''_\nu\geq 1$, and also
$$k'+k''+1+\nu'_{max}+\nu''_{max}-2m\geq 0,$$
and
$$l+\sum_{\nu}k'_{\nu}+\sum_{\nu}k''_{\nu}=k.$$

These terms \eqref{5.121} can be estimated using \eqref{5.108},...,\eqref{5.114}, as we have just done for the
terms Easy$_\nu^\mu(\zeta,t)$ and Safe$_2^\mu$ in \eqref{5.115}.

If each $k'_\nu$ and $k''_\nu$ is at most 3, then we can prove that the integral \eqref{5.121} has absolute value at most $C$ (this is the analogue of the "Easy" terms). If one of the $k'_\nu$ or $k_\nu''$ is equal to 4, then we can prove that \eqref{5.121} has $L^2-$norm at most $C$ as a function of $\zeta\in\Gamma_+$ (this is the analogue of the "Safe" terms). Note that the analogue of the $Danguerous$ term from lemma \ref{lemma6} is an $Easy$ term in the present context, because four derivatives fall on the known, harmless function $B$.

Consequently, for \eqref{5.120}, we learn that
\begin{equation*}
\left|\pa_\zeta^k \int_{w\in\Gamma_+(t)}\frac{\sin(z_1(\zeta,t)-z_1(w,t))}{\cosh(z_2(\zeta,t)-z_2(w,t))-\cos(z_1(\zeta,t)-z_1(w,t))}B(\zeta,w,t)dw
\right|\leq C
\end{equation*}
for $\zeta\in \Gamma_+(t)$, $0\leq k\leq 3$; and
\begin{equation*}
\pa_\zeta^4\int_{w\in\Gamma_+(t)}\frac{\sin(z_1(\zeta,t)-z_1(w,t))}{\cosh(z_2(\zeta,t)-z_2(w,t))-\cos(z_1(\zeta,t)-z_1(w,t))}B(\zeta,w,t)dw
\end{equation*}
has $L^2-$norm at most $C$ on $\Gamma_+(t)$.

Therefore,
\begin{align}
&\pa_\zeta^4\left\{\int_{w\in\Gamma_+(t)}\frac{\sin(z_1(\zeta,t)-z_1(w,t))B(\zeta,w,t)}
{\cosh(z_2(\zeta,t)-z_2(w,t))-\cos(z_1(\zeta,t)-z_1(w,t))}dw\pa_\zeta z
_\mu(\zeta,t)\right\}\nonumber\\
&=\int_{w\in\Gamma_+(t)}\frac{\sin(z_1(\zeta,t)-z_1(w,t))B(\zeta,w,t)dw}
{\cosh(z_2(\zeta,t)-z_2(w,t))-\cos(z_1(\zeta,t)-z_1(w,t))}\pa_\zeta^5 z_\mu(\zeta,t)+\tilde{\text{Error}}_\mu(\zeta,t)\label{5.122}
\end{align}
with
\begin{equation}\label{5.123}
\int_{\Gamma_+(t)}\left|\tilde{\text{Error}}_\mu(\zeta,t)\right|^2 d\zeta\leq C.
\end{equation}
We have succeeded in computing the fourth $\zeta-$derivative of another of the terms in \eqref{20.6}.

Finally, since the $H^4$-norms of $z_1(\zeta,t)-\zeta$ and $z_2(\zeta,t)$ on $\Gamma_+(t)$ are less than a controlled constant $C$, and since $A(\zeta,t)$ has 4 derivatives bounded by $C$, we have
\begin{align}
&\pa_\zeta^4\left\{A(\zeta,t)\pa_\zeta z_\mu(\zeta,t)\right\}=\nonumber\\
& A(\zeta,t)\pa_\zeta^5 z_\mu(\zeta,t)+\text{Error}^\sharp_\mu(\zeta,t),\label{5.124}
\end{align}
for $\zeta\in\Gamma_+(t)$, with

\begin{equation}
\int_{\Gamma_+(t)}\left|\text{Error}_\mu^\sharp(\zeta,t)\right|^2 d\zeta\leq C\label{5.125}.
\end{equation}

Now, comparing \eqref{5.118}, \eqref{5.119}, \eqref{5.122}, \eqref{5.123}, and \eqref{5.124}, \eqref{5.125} with the definition \eqref{20.6} of $J$, we learn that
\begin{align}
&\pa_\zeta^4 J(\zeta,t)=\nonumber\\
&\left\{A(\zeta,t)+\int_{w\in\Gamma_+(t)}\frac{\sin(z_1(\zeta,t)-z_1(w,t))B(\zeta,w,t)}
{\cosh(z_2(\zeta,t)-z_2(w,t))-\cos(z_1(\zeta,t)-z_1(w,t))}dw\right\}\pa_\zeta^5 z_\mu(\zeta,t)\nonumber\\
&+\int_{w\in\Gamma_+(t)}\frac{\sin(z_1(\zeta,t)-z_1(w,t))}
{\cosh(z_2(\zeta,t)-z_2(w,t))-\cos(z_1(\zeta,t)-z_1(w,t))}\left[\pa_\zeta^5 z_\mu(\zeta,t)-\pa_\zeta^5 z_\mu(w,t) \right]dw\nonumber\\
&+\text{Err}_\mu(\zeta,t)\quad \text{for $\zeta\in\Gamma_+(t)$}, \label{5.126}
\end{align}
where

\begin{equation}
\int_{\Gamma_+(t)}\left|\text{Err}_\mu(\zeta,t)\right|^2d\zeta \leq C.\label{5.127}
\end{equation}

We have succeeded in computing $\pa_\zeta^4 J(\zeta,t)$ modulo an error with $L^2-$norm at most $C$ on $\Gamma_+(t)$.

We may rewrite $\pa_\zeta^4 J(\zeta,t)$ in the form

\begin{align}
&\pa_\zeta^4 J(\zeta,t)=\nonumber\\
&\left\{A(\zeta,t)+P.V.\int_{w\in\Gamma_+(t)}\frac{\sin(z_1(\zeta,t)-z_1(w,t))\left(B(\zeta,w,t)+1\right)}
{\cosh(z_2(\zeta,t)-z_2(w,t))-\cos(z_1(\zeta,t)-z_1(w,t))}dw\right\}\pa_\zeta^5 z_\mu(\zeta,t)\nonumber\\
&+\int_{w\in\Gamma_+(t)}\left\{\frac{\pa_\zeta z_1(\zeta,t)}{\left(\pa_\zeta z_1(\zeta,t)\right)^2+\left(\pa_\zeta z_2(\zeta,t)\right)^2}\cot\left(\frac{\zeta-w}{2}\right)\right.\nonumber\\
&\quad\quad\quad\quad\quad\left . -\frac{\sin(z_1(\zeta,t)-z_1(w,t))}
{\cosh(z_2(\zeta,t)-z_2(w,t))-\cos(z_1(\zeta,t)-z_1(w,t))}\right\}\pa_w^5z_\mu(w,t)dw\nonumber\\
&-P.V.\int_{w\in\Gamma_+(t)}\frac{\pa_\zeta z_1(\zeta,t)}{\left(\pa_\zeta z_1(\zeta,t)\right)^2+\left(\pa_\zeta z_2(\zeta,t)\right)^2}\cot\left(\frac{\zeta-w}{2}\right)\pa_w^5z_\mu(w,t)dw\nonumber\\
&+\text{Err}_\mu(\zeta,t)\quad\qquad \text{for $\zeta\in\Gamma_+(t)$}, \label{5.128}
\end{align}
with Err$_\mu(\zeta,t)$ as in \eqref{5.127}.

Let us estimate the second-to-last integral in \eqref{5.128}.

The expression in curly brackets in that integral has $C^1-$norm at most $C$, as a function of $w\in\Gamma_+(t)$; this follows by arguments in sections \ref{seccion14} and \ref{seccion15}.

Therefore, integrating by parts in $w$, we obtain an integral of the form
\begin{equation*}
\int_{w\in\Gamma_+(t)}\left[\text{Quantity of absolute value at most $C$}\right]\cdot \pa_w^4z_\mu(w,t)dw;
\end{equation*}
such an integral has absolute value at most $C$.

Consequently, \eqref{5.128} may be rewritten in the form
\begin{align}
&\pa_\zeta^4 J(\zeta,t)=\nonumber\\
&\left\{A(\zeta,t)+P.V.\int_{w\in\Gamma_+(t)}\frac{\sin(z_1(\zeta,t)-z_1(w,t))\left(B(\zeta,w,t)+1\right)dw}{\cosh(z_2(\zeta,t)-z_2(w,t)-\cos(z_1(\zeta,t)-z_1(w,t)))}\right\}
\pa_\zeta^5z_\mu(\zeta,t)\nonumber\\
&-\frac{\pa_\zeta z_1(\zeta,t)}{\left(\pa_\zeta z_1(\zeta,t)\right)^2+\left(\pa_\zeta z_2(\zeta,t)\right)^2}
P.V. \int_{w\in\Gamma_+(t)}\cot\left(\frac{\zeta-w}{2}\right)\pa^5_w z_\mu(w,t)\nonumber\\
&+\text{Err}_\mu^{\sharp\sharp}(\zeta,t)\quad \qquad\ \text{for $\zeta\in\Gamma_+(t)$},\label{5.129}
\end{align}
with
\begin{equation}
\int_{\Gamma_+(t)}\left|\text{Err}_\mu^{\sharp\sharp}\right|^2d\zeta\leq C.\label{5.130}
\end{equation}
Now let
\begin{equation}
F_\mu(\zeta)=\pa_\zeta^4 z_\mu(\zeta,t)\label{5.131}
\end{equation}
Then
\begin{equation}
\int_{\zeta\in\Gamma_+(t)}\left|F_\mu(\zeta)\right|^2d\zeta\leq C.\label{5.132}
\end{equation}
We have
\begin{align}
&P.V.\int_{w\in\Gamma_+(t)}\cot\left(\frac{\zeta-w}{2}\right)\pa_w^5z_\mu(w,t)dw\nonumber\\
&=-P.V.\int_{w\in\Gamma_+(t)}\cot\left(\frac{\zeta-w}{2}\right)\left[F'_\mu(\zeta)-F'_\mu(w)\right]dw\nonumber\\
&=2\pi\Lambda F_\mu(\zeta)\quad\qquad \text{for $\zeta\in\Gamma_+(t)$
},\label{5.1321/2}
\end{align}
thanks to \eqref{ladefi}. Also, applying \eqref{dostrece} (with our present $h(t)$ in place of $h(x)$ in \eqref{dostrece}) we see that
\begin{equation}
\pa_\zeta^5 z_\mu(\zeta,t)=F'_\mu(\zeta,t)=-i \Lambda F_\mu(\zeta)+ Goof_\mu(\zeta)\label{5.133}
\end{equation}
for $\zeta\in \Gamma_+(t)$, where
\begin{equation}
\int_{\Gamma_+(t)}\left| Goof_\mu(\zeta)
\right|^2d\zeta\leq C\frac{1}{h(t)^2}\sum_{\pm}\int_{\Gamma_\pm(t)}
\left|F_\mu(\zeta)\right|^2d\zeta. \label{5.134}
\end{equation}
We recall that $h(t)>c$; see our definition of controlled constants in section \ref{seccion19}. Also, since $F_\mu(\overline{\zeta})=\overline{F_\mu(\zeta)}$, we have
\begin{equation*}
\int_{\Gamma_+(t)}\left|F_\mu(\zeta)\right|^2d\zeta=\int_{\Gamma_-(t)}\left|F_\mu(\zeta)\right|^2d\zeta\leq C
\end{equation*}
by \eqref{5.132}.

Hence, \eqref{5.134} yields
\begin{equation}
\int_{\Gamma_+(t)}\left|Goof_{\mu}(\zeta)\right|^2d\zeta\leq C. \label{5.135}
\end{equation}
Now, putting \eqref{5.1321/2}, \eqref{5.133} and \eqref{5.135} into \eqref{5.129}, and recalling \eqref{5.130},
we see that
\begin{align}
&\pa_\zeta^4 J(\zeta,t)=\nonumber\\
&-i\left\{A(\zeta,t)+P.V.\int_{w\in\Gamma_+(t)}\frac{\sin(z_1(\zeta,t)-z_1(w,t))\left(B(\zeta,w,t)+1\right)dw}{\cosh(z_2(\zeta,t)-z_2(w,t)-\cos(z_1(\zeta,t)-z_1(w,t)))}\right\}
\Lambda F_\mu(\zeta)\nonumber\\
&-\frac{2\pi\pa_\zeta z_1(\zeta,t)}{\left(\pa_\zeta z_1(\zeta,t)\right)^2+\left(\pa_\zeta z_2(\zeta,t)\right)^2}
\Lambda F_\mu(\zeta)\nonumber\\
&+\text{Erreur}_\mu(\zeta,t)\quad \qquad\ \text{for $\zeta\in\Gamma_+(t)$},\label{5.136}
\end{align}
with
\begin{equation}\label{5.137}
\int_{\Gamma_+(t)}\left|\text{Erreur}_\mu(\zeta,t)\right|^2d\zeta\leq C.
\end{equation}
(Here, we have also used the fact that the quantity in curly brackets in \eqref{5.129} is at most $C$, thanks to our assumptions on $A(\zeta,t)$, $B(\zeta,t)$ and on $z_\mu(\zeta,t)$.)

Also, from \eqref{5.133} and \eqref{5.135}, we have
\begin{equation}
ih'(t)\pa_\zeta^5 z_\mu(\zeta,t)=i h'(t)F_\mu'(\zeta)=h'(t)\Lambda F_\mu(\zeta)+Goof_\mu^{extra}(\zeta,t)\label{5.138}
\end{equation}
for $\zeta\in\Gamma_+(t)$, with
\begin{equation}\label{5.139}
\int_{\Gamma_+(t)}\left| Goof_\mu^{extra}(\zeta,t)
\right|^2d\zeta\leq C.
\end{equation}
From \eqref{5.136},..., \eqref{5.139}, we see that
\begin{align}
&\pa_\zeta^4J(\zeta,t)+ih'(t)\pa_\zeta^5 z_\mu(\zeta,t)\nonumber\\
&=\left\{-iA(\zeta,t)-iP.V.\int_{w\in\Gamma_+(t)}
\frac{\sin(z_1(\zeta,t)-z_1(w,t))\left(B(\zeta,w,t)+1\right)dw}
{\cosh(z_2(\zeta,t)-z_2(w,t)-\cos(z_1(\zeta,t)-z_1(w,t)))}\right.\nonumber\\
&\left.-\frac{2\pi\pa_\zeta z_1(\zeta,t)}{\left(\pa_\zeta z_1(\zeta,t)\right)^2+\left(\pa_\zeta z_2(\zeta,t)\right)^2}
+h'(t)\right\}\Lambda F_\mu(\zeta)+\text{Errr}_\mu(\zeta) \label{5.140}\end{align}
for $\zeta\in\Gamma_+(t)$, with
\begin{equation}\label{5.141}
\int_{\Gamma_+(t)}\left |\text{Errr}_\mu(\zeta,t)\right|^2d\zeta\leq C.
\end{equation}
From \eqref{5.131}, \eqref{5.132} and \eqref{5.140}, \eqref{5.141}, we find that
\begin{align}
&\Re \int_{\zeta\in\Gamma_+(t)}\overline{\pa_\zeta^4 z_\mu(\zeta,t)}\left[\pa_\zeta^4 J(\zeta,t)+ih'(t)\pa_\zeta^5 z_\mu(\zeta,t)\right]d\zeta\nonumber\\
&\geq \Re \int_{\zeta\in\Gamma_+(t)}\overline{F_\mu(\zeta)}R(\zeta)\Lambda F_\mu(\zeta)d\zeta-C,\label{5.142}
\end{align}
where
\begin{align}
R(\zeta)=&-iA(\zeta,t)-i P.V. \int_{w\in\Gamma_+(t)}
\frac{\sin(z_1(\zeta,t)-z_1(w,t))\left(B(\zeta,w,t)+1\right)dw}{\cosh(z_2(\zeta,t)-z_2(w,t)-\cos(z_1(\zeta,t)-z_1(w,t)))}\nonumber\\
&-\frac{2\pi\pa_\zeta z_1(\zeta,t)}{\left(\pa_\zeta z_1(\zeta,t)\right)^2+\left(\pa_\zeta z_2(\zeta,t)\right)^2}+h'(t).\label{5.143}
\end{align}
We note that $R(\zeta)$ has $C^2-$norm at most $C$ on $\Gamma_+(t)$.
(Compare with \eqref{quinceuno}, \eqref{quincedos} and \eqref{paraz}.

From  \eqref{six} and our present estimate \eqref{5.132}, we see that
\begin{equation}
\Re\int_{\zeta\in\Gamma_+(t)} \overline{F_\mu(\zeta)}R(\zeta)\Lambda F_\mu(\zeta)d\zeta\geq \Re\int_{\zeta\in\Gamma_+(t)}\overline{F_\mu(\zeta)}\Re \left(R(\zeta)\right)\Lambda F_\mu(\zeta)d\zeta-C.\label{5.144}
\end{equation}
Moreover, $\Re(R(\zeta))\geq 0$ on $\Gamma_+(t)$, thanks to the Rayleigh-Taylor hypothesis \eqref{19.22} and \eqref{19.23}.
(Recall: We are solving backwards in time, and we stop as soon as $RT$ fails.) Since also the $C^2-$norm of $\Re(R(\zeta))$ on $\Gamma_+(t)$ is at most $C$, the  inequality \eqref{six} now tell us that
\begin{equation}\Re\int_{\Gamma_+(t)}\overline{F_\mu(\zeta)}\Re \left(R(\zeta)\right)\Lambda F_\mu(\zeta)d\zeta\geq
-C\int_{\Gamma_+(t)}\left|F_\mu(\zeta)\right|^2d\zeta\geq -C,\label{5.145}\end{equation}
by \eqref{5.132}.

From \eqref{5.142}, \eqref{5.144}, \eqref{5.145}, we learn that
\begin{equation}
\Re\int_{\Gamma_+(t)}\overline{\pa_\zeta^4 z_\mu(\zeta,t)}\left[\pa_\zeta^4 J(\zeta,t)+ih'(t)\pa_\zeta^5z_\mu(\zeta,t)\right]d\zeta\geq -C.
\end{equation}
Recalling \eqref{20.9}, we conclude that
\begin{equation}
\frac{d}{dt}\int_{\zeta\in\Gamma_+(t)}\left|\pa_\zeta^4z_\mu(\zeta,t)\right|^2d\zeta\geq -C.\label{5.146}
\end{equation}
Since we are solving backwards in time, this is our desired energy inequality.

We conclude that our Galerkin problem can be solved in a suitable time interval $[t_0-\delta,t_0]$, maintaining control of the $H^4-$norm, the chord-arc constant, and a positive lower bound for the Rayleigh-Taylor function. This control is uniform in the parameter N used to define our Galerkin approximation. Here, $\delta$ may be taken to be a small enough controlled constant.

It is now a routine exercise to deduce the following local existence theorem for the Generalized Muskat Problem,

\begin{lemma}\label{5.X}
Let $h(t)$, $I_{time}$, $\Omega(t)$, $\Omega$, $\Omega^+$, $A(\zeta,t)$, $B(\zeta,w,t)$, $z_1^0(\zeta)$, $z_2^0(\zeta)$, and $RT^0(\zeta)$ be as in \eqref{19.1},...,\eqref{19.23}. Then there exist functions $z_1(\zeta,t)$, $z_2(\zeta,t)$ on $\Omega$, such that \eqref{19.24},...,\eqref{19.33} hold.
\end{lemma}
This lemma provides only local existence, since hypothesis \eqref{19.23} demands that the length of our time interval is less than a small enough controlled constant.

\subsection{Local Existence for Muskat Solutions in the Complex Domain}\label{seccion5.5}
In this section, we give the following local-in-time result.
\begin{thm}\label{5XX}
\begin{enumerate}
\item Let $h(x,t)$ be a positive real-analytic function on $\T\times[t_0-\delta,t_0+\delta]$; set
\begin{align*}
\Gamma_{\pm}(t)=&\left\{\zeta\in\C\,:\, \Im \zeta=\pm h(\Re\zeta,t) \right\}
\end{align*}
and
\begin{align*}
\Omega(t)=&\left\{\zeta\in \C\,:\,|\Im\zeta|<h(\Re\zeta,t)\right\}.
\end{align*}
Let $z^0(\zeta)=(z^0_1(\zeta),\,z^0_2(\zeta))$ be continuous on $\Omega(t_0)^{\text{closure}}$ and
holomorphic on $\Omega(t_0)$. Suppose that $z_1^0(\zeta)-\zeta$ and $z^0_2(\zeta)$ are $2\pi-$periodic and take real values for real $\zeta$. Also, suppose that  $z_1^0(\zeta)-\zeta$ and $z^0_2(\zeta)$ belong to $H^4(\Gamma_\pm(t_0))$.
\item Assume the chord-arc condition
\begin{align*}
\left|\cosh(z_2^0(\zeta)-z_2^0(w))-\cos(z_1^0(\zeta)-z_1^0(w))\right|\geq & c\left[||\Re(\zeta-w)||+|\Im(\zeta-w)|\right]^2,
\end{align*}
for $\zeta$, $w\in \Omega(t_0)^{\text{closure}}$.
\item Finally, assume the Rayleigh-Taylor condition $RT(\zeta)>0$ for all $\zeta\in \Gamma_+(t)$, where $RT$ is defined by \eqref{18.15} in terms of $(z^0_1(\zeta),\,z^0_2(\zeta))$.
\end{enumerate}
Then there exist positive numbers $\delta'$, $C'$, $c'$ and, functions $(z_1(\zeta,t),\,z_2(\zeta,t))$, continuous
on
\begin{equation*}
\left\{(\zeta,t)\in \C\times [t_0-\delta',t_0]\,:\,\zeta\in\Omega(t)^{\text{closure}}\right\},
\end{equation*}
for which the following hold.
\begin{description}
\item A. For each fixed $t\in[t_0-\delta',t_0]$, the functions $z_1(\zeta,t)-\zeta$ and $z_2(\zeta,t)$ are $2\pi-$periodic, holomorphic on $\Omega(t)$, real-valued for real $\zeta$, and have $H^4-$norm at most $C'$ on $\Gamma_\pm(t)$. Moreover,
    \begin{align*}
    \left|\cosh(z_2(\zeta,t)-z_2(w,t))-\cos(z_1(\zeta,t)-z_1(w,t))\right|\geq c'\left[||\Re(\zeta-w)||+|\Im(\zeta-w)|\right]^2,
    \end{align*}
    for $\zeta$, $w\in \Omega^\text{closure}$, $t\in[t_0-\delta',t_0]$.

    Also, $RT(\zeta,t)>0$ (strict inequality) for all $\zeta\in \Gamma_+(t)$, $t\in[t_0-\delta',t_0]$.

    (Here, $RT(\zeta,t)$ is defined from $z_1(\zeta,t)$, $z_2(\zeta,t)$ by \eqref{18.15}).
\item B. The functions $z_1(\zeta,t)$, $z_2(\zeta,t)$ satisfy the Muskat equation
\begin{equation*}
\pa_t z_\mu(\zeta,t)=\int_{w\in\Gamma_+}\frac{\sin(z_1(\zeta,t)-z_1(w,t))[\pa_\zeta z_\mu(\zeta,t)-\pa_w z_\mu(w,t)]}
{\cosh(z_2(\zeta,t)-z_2(w,t))-\cos(z_1(\zeta,t)-z_1(w,t))}dw,
\end{equation*}
for $\mu=1,\,2$, $\zeta\in \Omega(t)^\text{closure}$, $t\in[t_0-\delta',t_0]$, with initial conditions
\begin{equation*}
z_\mu(\zeta,t_0)=z_\mu^0(\zeta)
 \end{equation*}
 for $\mu=1,\,2$ and $\zeta\in \Omega(t_0)^\text{closure}$.
\item C. The constants $\delta'$, $C'$, $c'$ above may be taken to depend only on the following:
\begin{itemize}
\item $\delta$ (see 1.);
\item The function $h(x,t)$;
\item Bounds for the $H^4-$norms of $z_1^0(\zeta)-\zeta$ and $z_2^0(\zeta)$ on $\Gamma_+(t_0)$;
\item The constant $c$ in the cord-arc condition 2; and
\item A positive lower bound for the Rayleigh-Taylor function in 3.
\end{itemize}
\end{description}
\end{thm}
Proof: In section \ref{seccion18}, we showed that the Muskat equation on a time-varying domain $\Omega(t)$ is
equivalent to a generalized Muskat equation on a time-varying strip. In section \ref{seccion19} and \ref{seccion20},
we proved a local existence theorem for generalized Muskat. The present result follows at once.

\subsection{Continuing Back to Negative Time}\label{seccion5.6}

Let $\underline{z}=(\underline{z}_1(\zeta,t),\,\underline{z}_2(\zeta,t))$, and $A$, $\tau$, $\lambda$, $\kappa$, $h(x,t)$, $\hbar(x,t)$, $\Omega(t)$, $\Gamma_\pm(t)$ be as in section 3. Recall that $\hbar(x,\tau^2)< h(x,\tau^2)$ for all $x$.

We suppose we are given initial data $z^0(\zeta)=(z_1^0(\zeta),\, z_2^0(\zeta))$ on $\Omega(\tau)^\text{closure}$, holomorphic in $\Omega(\tau)$, and real-valued for real $\zeta$.

We assume that the functions $z_\mu^0(\zeta)-\un_\mu(\zeta,\tau)$ ($\mu=1,\,2$) are $2\pi-$periodic, and that their $H^4(\Omega(\tau))-$norms are less than $10^{-2}\lambda$.

Starting from the initial data $z^0$ at time $t=\tau$, we solve the Muskat problem backwards in time, until we reach a time $t_{\text{least}}$ at which either the assumptions of theorem \ref{mean} break down, or else we reach $t_\text{least}=\tau^2$. We can produce such a Muskat solution, thanks to a routine argument that combines the local existence result theorem  \ref{5XX} with the a-priori estimates given by theorem \ref{mean}.

For $t\in[t_\text{least},\,\tau]\subset[\tau^2,\tau]$, we have
\begin{equation*}
\frac{d}{dt}||z(\cdot,t)-\un(\cdot,t)||^2_{H^4(\Omega(t))}\geq-C(A)\lambda^2.
\end{equation*}

Since also
\begin{equation*}
||z(\cdot,\tau)-\un(\cdot,\tau)||^2_{H^4(\Omega(\tau))}\leq 10^{-4}\lambda^2,
\end{equation*}
it follows that
\begin{equation*}
||z(\cdot,t_{\text{least}})-\un(\cdot,t_{\text{least}})||^2_{H^4(\Omega(t_{\text{least}}))}\leq 10^{-4}\lambda^2+C(A)\tau\lambda^2\leq 10^{-3}\lambda^2.
\end{equation*}
Consequently, the assumptions of theorem \ref{mean} cannot break down at $t=t_{least}$, so we have $t_{\text{least}}=\tau^2$.

Thus, we have solved the Muskat equation back to time $t=\tau^2$. Our Muskat solution $z(\zeta,t)$ satisfies:
 \begin{enumerate}
 \item It is  holomorphic on
 \begin{equation*}
 \Omega(t)=\{|\Im \zeta|< h(\Re \zeta,t)\}
 \end{equation*}
for each fixed $t\in [\tau^2,\tau]$.
\item We have
\begin{equation*}
||z(\cdot,t)-\un(\cdot,t)||_{H^4(\Omega(t))}\leq \lambda
\end{equation*}
for each $t\in[\tau^2,\tau]$.
\item In particular, for time $t=\tau^2$, the Muskat solution $z(\zeta,\tau^2)$ is holomorphic on the smaller domain
\begin{equation*}
\Omega(t)=\{|\Im \zeta|<\hbar(\Re \zeta,t)\}
\end{equation*}
for $t=\tau^2$ (not the same as 1, due to a defect in our notation, noted in section \ref{seccion33}).
\item Estimate 2 for the domain 1 implies easily that
\begin{equation*}
||z(\cdot,\tau^2)-\un(\cdot,\tau^2)||_{H^4(\Omega(\tau^2))}\leq C_1\lambda
\end{equation*}
for the domain 3, as one sees from the Cauchy integral formula and the $L^2-$boundedness of the Hilbert
transform.  Here, $C_1$ is a universal constant.
\end{enumerate}

For the next few paragraphs, we define $\Omega(t)$ by 3.

Now, starting from the initial data $z(\cdot,\tau^2)$ at time $t=\tau^2$, we solve the Muskat equation backwards in time, until we reach a time $t_{\text{least}}\in [-\tau^2,\,\tau^2]$ at which either the assumptions of theorem \ref{mean} break down, or else we reach $t_{\text{least}}=-\tau^2$. (When we apply theorem \ref{mean}, we use $2C_1\lambda$ in place of $\lambda$; see 4.)

Again, a routine argument using theorem \ref{5XX} and theorem \ref{mean} allows us to produce such a Muskat solution.

For $t\in [t_{\text{least}},\,\tau^2]$, we learn from theorem \ref{mean} that
\begin{equation*}
\frac{d}{dt}||z(\cdot,t)-\un(\cdot,t)||^2_{H^4(\Omega(t))}\geq -C(A)\tau^{-1}\lambda^2.
\end{equation*}
Hence, by 4, we have
\begin{equation*}
||z(\cdot,t_{\text{least}})-\un(\cdot,t_{\text{least}})||^2_{H^4(\Omega(t_{\text{least}}))}
\leq C_1^2\lambda^2+C(A)\tau^{-1}(C_1\lambda)^2\tau^2\leq 2C_1^2\lambda^2.
\end{equation*}
(Notice that the large constant $C(A)\tau^{-1}\lambda^2$ in the conclusion of theorem \ref{mean} is compensated by the short time interval of length $2\tau^2$).

As before, it now follows easily that the assumptions of theorem \ref{mean} (with $2C_1\lambda$ in place of $\lambda$) cannot break down at $t_\text{least}$. Thus, we get all the way back to time $t_{least}=-\tau^2$.

We have therefore proven the following result.

\begin{lemma}\label{5YY}
Let $z^0(\zeta)=(z^0_1(\zeta),\,z^0_2(\zeta))$ be given, with $z_\mu^0(\zeta)-\un_\mu(\zeta,\tau)$ $2\pi-$periodic, holomorphic in $\Omega(\tau)$, and having norm at most $\lambda$ in $H^4(\Omega(\tau))$.

Then there exists a Muskat solution $z(\zeta,t)=(z_1(\zeta,t),\,z_2(\zeta,t))$ defined for $t\in[-\tau^2,\,\tau]$, holomorphic in $\Omega(t)$ for each $t\in[-\tau^2,\,\tau]$, equal to $z^0(\zeta)$ at $t=\tau$, and satisfying
\begin{equation*}
||z_\mu(\cdot,t)-\un_\mu(\cdot,t)||_{H^4(\Omega(t))}\leq C\lambda
\end{equation*}
for each $t\in[-\tau^2,\tau]$, where $C$ is a universal constant.
\end{lemma}

\subsection{Proof of the Main Theorem \ref{conclusion}}\label{ultima}
In the previous sections, we worked with fixed parameters $A$, $\tau$, $\lambda$, $\kappa>0$
satisfying the conditions imposed in section \ref{seccion33}. In this section, we will let $\kappa\to 0^+$. We write
$h_\kappa(x,t)$ to denote the function called $h(x,t)$ before; and we write $\Omega_\kappa(t)$ to denote the region called $\Omega(t)$ in the previous sections.

Let $f_\kappa$ be the analytic function on $\Omega_{\kappa}(\tau)$ defined by setting
\begin{align*}
\pa_x^4 f_{\kappa}(x)=&\log\left(\sin^2\left(\frac{x}{2}\right)
+\sinh^2\left(\frac{\kappa}{2}\right)\right)-\frac{1}{2\pi}\int_{-\pi}^{\pi}\log\left(\sin^2\left(\frac{s}{2}\right)
+\sinh^2\left(\frac{\kappa}{2}\right)\right)ds,\\
\int_{-\pi}^{\pi}f_\kappa(x)dx=&0,
\end{align*}
for $x\in\T$.
One sees easily that $f_\kappa(\zeta)$ can be taken to be analytic outside the union of the slits
\begin{equation*}
\{\Re \zeta=2\pi m,\, \pm \Im \zeta\geq \kappa\}\quad (m\in \Z);
\end{equation*}
thus, $f_\kappa(\zeta)$ is analytic on $\Omega_\kappa(\tau)$ as claimed.

Moreover, the norm of $f_\kappa(\zeta)$ in $H^4(\Omega_\kappa(\tau))$ is bounded by a universal constant $C$, since the function $\log|x|$ belongs to $L^2([-\pi,\pi])$.

Lemma \ref{5YY} from the preceding section gives us a Muskat solution $$z_\kappa(\zeta,t)=(z_{1,\,\kappa}(\zeta,t),z_{2,\,\kappa}(\zeta,t))$$ defined for $t\in[-\tau^2,\tau]$ and $\zeta\in\Omega_\kappa(t)$, with initial condition
$$z_\kappa(\zeta,\tau)=(\un_1(\zeta,\tau)+\lambda f_{\kappa}(\zeta),\, \un_2(\zeta,\tau)).$$
Moreover, $z_\kappa(\zeta,t)$ is holomorphic on $\Omega_\kappa(t)$ for fixed $t$, and we have
$$||z_\kappa(\cdot,t)-\un(\cdot,t)||_{H^4(\Omega_\kappa(t))}\leq C\lambda$$
for $t\in[-\tau^2,\,\tau]$.

A routine limiting argument now produces a Muskat solution
\begin{equation*}z(\zeta,t)=(z_1(\zeta,t),\,z_2(\zeta,t)),\quad t\in[-\tau^2,\tau],\quad \zeta\in \cap_{\kappa>0}\Omega_\kappa^\text{closure}(t)\end{equation*}

with the following properties
\begin{enumerate}
\item $z_1(\zeta,\tau)=\un_1(\zeta,\tau)+\lambda f_0(\zeta)$, $z_2(\zeta,\tau)=\un_2(\zeta,\tau)$, where $f_0(\zeta)=\lim_{\kappa\to 0^+}f_{\kappa}(\zeta)$;
\item $z(\zeta,t)$ is holomorphic in \begin{equation*}\Omega(t)=\cap_{\kappa>0}\Omega_\kappa(t)\end{equation*}
for each fixed $t\in[-\tau^2,\tau)$;
\item $||z(\cdot,t)-\un(\cdot,t)||_{H^4(\Omega(t))}\leq C\lambda$ for $t\in[-\tau^2,\,\tau]$.
\end{enumerate}
Since $\pa_\zeta^4 f_0(\zeta)$ has a logarithmic singularity at the origin, it follows from 1 that the curve $x\mapsto z(x,t)$ is not $C^4-$smooth.

Recalling that $\pa_\zeta\un_2(0,\tau)\neq 0$, we see that there is no way to reparameterize the curve $x\mapsto z(x,t)$ to make it $C^4$ smooth.
However, for $t\in[-\tau^2,\,\tau)$, we see from 2 that the curve $x\mapsto z(x,t)$ is real-analytic.

Finally, for $t=-\tau^2$, we learn from 3 that
\begin{equation}\label{yaesta}
|\pa_x z_1(x,-\tau^2)-\pa_x \un_1(x,-\tau^2)|\leq C\lambda,
\end{equation}
for all $x\in\T$.

Since, from \eqref{grafo}, $\pa_x \un_1(x,-\tau^2)\geq c_{20}\tau^2$ for all real $x$, we learn form \eqref{yaesta} that
\begin{equation}
\pa_x z_1(x,-\tau^2)>0,
\end{equation}
for all $x\in\T$, since $\lambda$ is smaller than $c\tau^2$. Thus, at time $t=-\tau^2$, our Muskat curve is a graph. The proof of our main theorem \eqref{conclusion} is complete.

\subsection*{{\bf Acknowledgments}}

\smallskip

 AC, DC and FG were partially supported by the grant {\sc MTM2008-03754} of the MCINN (Spain) and
the grant StG-203138CDSIF  of the ERC. CF was partially supported by
NSF grant DMS-0901040.
FG was partially supported by NSF grant DMS-0901810.
We are grateful
for the support of the Fundaci\'on General del CSIC.

\begin{tabular}{ll}
\textbf{Angel Castro} & \textbf{Diego C\'ordoba} \\
{\small D\'eparment de Math\'ematiques et Applications} & {\small Instituto de Ciencias Matem\'aticas}\\
{\small \'Ecole Normale Sup\'erieure} & {\small Consejo Superior de Investigaciones Cient\'ificas}\\
{\small 45, Rue d'Ulm, 75005 Paris} & {\small Serrano 123, 28006 Madrid, Spain}\\
{\small Email: castro@dma.ens.fr} & {\small Email: dcg@icmat.es}\\
   & \\
\textbf{Charles Fefferman} & \textbf{Francisco Gancedo}\\
{\small Department of Mathematics} & {\small Departamento de An\'alisis Matem\'atico}\\
{\small Princeton University} & {\small Universidad de Sevilla}\\
{\small 1102 Fine Hall, Washington Rd, } & {\small  Tarfia, s/n}\\
{\small Princeton, NJ 08544, USA} & {\small Campus Reina Mercedes, 41012 Sevilla }\\
 {\small Email: cf@math.princeton.edu} & {\small Email: fgancedo@us.es}\\
  \\

\end{tabular}


\begin{thebibliography}{99}
\bibitem{turning}A. Castro, D. C\'ordoba, C. Fefferman, F. Gancedo and Mar\'ia L\'opez-Fern\'andez.
Turnning waves and breakdown for incompressible flows. \emph{Proc. Natl. Acad. Sci.} 108, no. 12, 4754-4759 (2011).

\bibitem{CCFGL} A. Castro, D. C\'ordoba, C. Fefferman, F. Gancedo and M. Lopez-Fernandez. Rayleigh-Taylor breakdown for the Muskat problem with applications to water waves. To appear in Annals of Math (2012).

\bibitem{ccgs} P. Constantin, D. C\'ordoba, F. Gancedo and R.M. Strain. On the global existence  for
the Muskat problem. \emph{To appear in JEMS} (2012).

\bibitem{Peter} P. Constantin and M. Pugh. Global solutions for small data to the
Hele-Shaw problem. \emph{Nonlinearity}, 6 (1993), 393 - 415.

\bibitem{AD} A. C\'ordoba, D. C\'ordoba. A pointwise estimate for fractionary derivatives with applications
to partial differential equations. \emph{Proc. Natl. Acad. Sci. USA} 100, 26,
(2003), 15316-15317.

\bibitem{CFG} D. C\'ordoba, D. Faraco and F. Gancedo. Lack of uniqueness for weak solutions of the incompressible porous media equation. \emph{Arch. Rat. Mech. Anal.} 200 (2011), 725-746.

\bibitem{DY} D. C\'ordoba and F. Gancedo. Contour dynamics of incompressible 3-D fluids
in a porous medium with different densities. \emph{Comm. Math.
Phys.} 273, 2,(2007), 445-471.


\bibitem{Esch2} J. Escher and B.-V. Matioc. On the parabolicity of the Muskat problem: Well-posedness, fingering,
and stability results. \emph{Z. Anal. Awend.} 30, 193-218, (2011).

\bibitem{H-S} Hele-Shaw. On the motion of a viscous fluid between two parallel plates. \emph{Trans. Royal
Inst. Nav. Archit.}, London 40, 21, (1898).

\bibitem{Muskat} M. Muskat. Two fluid systems in porous media. The encroachment of water into an oil sand.
\emph{Physics}, 5, (1934), 250-264.

\bibitem{Nirenberg} L. Nirenberg. An abstract form of the nonlinear Cauchy-Kowalewski theorem.
\emph{J. Differential Geometry}, 6, (1972), 561-576.

\bibitem{Nishida} T. Nishida. A note on a theorem of Nirenberg. \emph{J. Differential Geometry}, 12, (1977), 629-633.



\bibitem{S-T} P.G. Saffman and Taylor. The penetration of a fluid into a porous medium or
Hele-Shaw cell containing a more viscous liquid.
\emph{Proc. R. Soc. London, Ser. A} 245, (1958), 312-329.

\bibitem{Shvydkoy} R. Shvydkoy. Convex integration for a class of active scalar equations. J. Amer. Math. Soc. 24 (2011) 1159-1174.

\bibitem{SCH} M. Siegel, R. Caflisch and S. Howison. Global
Existence, Singular Solutions, and Ill-Posedness for the Muskat
Problem. \emph{Comm. Pure and Appl. Math.}, 57, (2004), 1374-1411.

\bibitem{Eli} E. M. Stein. Harmonic Analysis. Princeton University Press. Princeton, New Jersey (1993)


\bibitem{Laszlo} L. Szekelyhidi Jr. Relaxation of the incompressible porous media equation. arXiv:1102.2597.

\bibitem{Otto1} F. Otto. Evolution of microstructure in unstable porous media flow: a relaxational approach.  \emph{Comm. Pure and Appl. Math.}, 52,  7(1999), 873-915.


\bibitem{Otto2} F. Otto. Evolution of microstructure: an example. In ergodic theory, analysis and efficient simulation of dynamical systems. \emph{Springer}, Berlin, 2001, pp. 501-522.

\bibitem{Yi2} F. Yi. Global classical solution of Muskat free boundary problem, \emph{J. Math. Anal. Appl.}, 288 (2003), 442-461.

\end{thebibliography}
\end{document}